\documentclass[a4paper,11pt]{amsart}

\usepackage{IntegralDifferentials}

\title{Models and Integral Differentials of Hyperelliptic Curves}
\author{Simone Muselli}
\address{University of Bristol, Bristol, UK.}

\begin{document}

\maketitle

\begin{abstract}
    Let $C: y^2=f(x)$ be a hyperelliptic curve of genus $g\geq 1$, defined over a complete discretely valued field $K$, with ring of integers $O_K$. Under certain conditions on $C$, mild when residue characteristic is not $2$, we explicitly construct the minimal regular model with normal crossings $\mathcal{C}/O_K$ of $C$. In the same setting we determine a basis of integral differentials of $C$, that is an $O_K$-basis for the global sections of the relative dualising sheaf $\omega_{\mathcal{C}/O_K}$.
\end{abstract}

\section{Introduction}\label{IntroductionSection}
The purpose of this paper is to construct regular models of hyperelliptic curves and to describe a basis of integral differentials attached to them. Moreover, we want these constructions explicit and easy to compute. 

\subsection{Overview}\label{MainResultsIntroductionSubsection}

To describe the arithmetic of curves over global fields, for example in the study of the Birch \& Swinnerton-Dyer conjecture, it is essential to understand regular models and integral differentials over all primes, including those with very bad reduction.
Constructing regular models of curves over discrete valuation rings is not an easy problem, even in the hyperelliptic curve case. In fact, there is no practical algorithm able to determine a model, unless the genus of the curve is $1$ or we have some tameness or nondegeneracy hypothesis.

One possible approach to tackle this problem is giving a full classification of possible regular models in a fixed genus, as done by the Kodaira--N\'{e}ron (\cite{Kod}, \cite{Ner}) and Namikawa--Ueno (\cite{NU}, \cite{Liu94}) classifications for curves of genera $1$ and $2$, respectively. However, this strategy seems impractical in general, since the number of models grows fast with the genus. Recently, new approaches based on clusters \cite{D2M2}, Newton polytopes \cite{Dok}, and MacLane valuations \cite{OW}, have been developed (see \S\ref{RelatedWorksSection} for more detail). 

On one side, clusters define nice and clear invariants from which one can extract information on the local arithmetic of hyperelliptic curves. 
Such invariants turn out to be particularly useful from a Galois theoretical point of view. However, for describing regular models, restrictions on the reduction type of the curve and on the residue characteristic of its base field (\cite{D2M2}, \cite{FN}) need to be imposed. On the other side, Newton polytopes and MacLane valuations have a potential to solve the problem in general, but the respective constructions are more algorithmic and so do not give the result in closed form. Furthermore, they often depend on the chosen equation rather than on the curve itself.

In this paper, we present a new approach that preserves both positive aspects from the above and provides a link between the two sides. We describe a model from simple invariants defined from what we call \textit{rational cluster picture} (Definition \ref{RationalClusterPictureIntroductionDefinition}). This object modifies the theory in \cite{D2M2} and appears to be more suitable for our purpose (see \S\ref{RationalClusterPictureSubsection}). 
In fact, the rational cluster picture also carries intrinsic connections with the other presented approaches, as it is closely related to Newton polygons and to degree $1$ MacLane valuations (see \cite{FGMN}).
When these valuations are enough to describe a regular model we say that the curve has an \textit{almost rational cluster picture} (Definition \ref{AlmostRationalIntroductionDefinition}; see also \ref{NewtonPolygonAlmostRationaCharacterisationCorollary}, \ref{AlmostRationalDiskCharacterisationProposition}). It turns out that the approach even works in residue characteristic 2, under an extra assumption that the curve is \textit{$y$-regular} (Definition \ref{yRegularIntroductionDefinition}). Our main result is:


\emph{Let $K$ be a complete\footnote{The assumption on the completeness of $K$ is not restrictive since regular models do not change under completion of the base field.} discretely valued field with $\mathrm{char}(K)\neq 2$,
and let $K^{nr}$ be its maximal unramified extension. Let $C/K$ be a hyperelliptic curve, having an almost rational cluster picture over $K^{nr}$. If the residue characteristic of $K$ is $2$, assume that $C_{K^{nr}}$ is $y$-regular. Then via the rational cluster picture we determine:
\begin{enumerate}[label=(\roman*)]
    \item the minimal regular model with normal crossings $\mathcal{C}^\mathrm{min}$,\label{maini}
    \item a basis of integral differentials of $C$.\label{mainii}
\end{enumerate}}


This result applies to a wide class of curves, covering wild cases and base fields with even residue characteristic. For example, if $g=2$, then $107$ out of $120$ Namikawa-Ueno types (\cite{NU}) arise from hyperelliptic curves satisfying the conditions of our theorem. In addition, the author believes it has a potential to solve the problem in general.
Heuristically speaking, the rational clusters invariants are expected to extend to general MacLane valuations. This approach could eventually lead to a full characterisation of minimal models with normal crossings of hyperelliptic curves (over any discretely valued field). 

\subsection{Main results}
We will now present (a simplified version of) the main results of this paper. We will then illustrate them with an explicit example in \S\ref{ExampleSubsection}.

Let $K$ be a complete discretely valued field of residue characteristic $p$, with normalised discrete valuation $v$ and ring of integers $O_K$. We require $\mathrm{char}(K)$ to be not $2$, but we allow $p=2$ and $p=0$. In this subsection we will assume for simplicity that $K=K^{nr}$. Extend the valuation $v$ to an algebraic closure $\bar K$ of $K$. 
Let $C/K$ be a hyperelliptic curve, i.e. a geometrically connected smooth projective curve, double cover of $\P^1_K$. Let $g$ be the genus of $C$. Assume $g\geq 1$. Fix a Weierstrass equation 
\[C:y^2=f(x).\] Let $\roots$ be the set of roots of $f$ in $\bar K$. Thus
\[f(x)=c_f\prod_{r\in\roots}(x-r).\]
For any $r,r'\in\roots$, with $r\neq r'$, denote by $\mathcal{D}_{r,r'}$ the smallest $v$-adic disc containing $r$ and $r'$.

\begin{defn}[Definition \ref{AlmostRationalDefinition}]\label{AlmostRationalIntroductionDefinition}
We say that $C$ has an \textit{almost rational cluster picture} if for any roots $r,r'\in\roots$ with $r\neq r'$, either
\begin{enumerate}
    \item [(a)] $\mathcal{D}_{r,r'}\cap K\neq\varnothing$, or
    \item [(b)] $p>0$ and $|\mathcal{D}_{r,r'}\cap\roots|\leq |v(r-w)|_p$ for some $w\in K$,
\end{enumerate}
where $|\cdot|_p$ denotes the canonical $p$-adic absolute value on $\Q$.
\end{defn}

\begin{defn}
A \textit{rational cluster} is a non-empty subset $\s\subset\roots$ of the form $\mathcal{D}\cap\roots$, where $\mathcal{D}$ is a $v$-adic disc $\mathcal{D}=\{x\in\bar K\mid v(x-w)\geq \rho\}$ for some $w\in K$ and $\rho\in\Q$. 
We denote by $\Sigma_K$ the set of rational clusters.
\end{defn}





In the following definition we introduce most of the notation and quantities, associated with rational clusters, needed in order to state our main theorems. 

\begin{defn}
For any $\s\in\Sigma_K$ we say:
\begin{center}
\begin{tabular}{|l@{ if }l|}
\hline
$\s$ proper, & $|\s|>1$\cr
$\s'$ is a child of $\s$, & $\s'\in\Sigma_K$ and $\s'\subsetneq\s$ is a maximal subcluster \cr
$\s$ minimal, & $\s$ has no proper children\cr
$\s$ \"{u}bereven, & $\s=\bigcup_{\s'\text{ child of }\s}\s'$ and $|\s'|$ even for all children $\s'$ of $\s$\cr
\hline
\end{tabular}
\end{center}
Moreover, we write $\s'<\s$, or $\s=P(\s')$, for a child $\s'$ of $\s$, and $r\wedge\s$ for the smallest rational cluster containing the root $r\in\roots$ and $\s$. 

Let $\mathring\Sigma_K$ be the set of proper rational clusters. For any $\s\in\mathring\Sigma_K$, define its \textit{radius} \[\rho_\s=\max_{w\in K}\min_{r\in\s} v(r-w)\] and the following quantities:

\begin{center}
\begin{tabular}{|l@{$=\>\,$}l|}
\hline
            $b_\s$            & denominator of $\rho_\s$\cr
            $\epsilon_\s$     & $v(c_f) + \sum_{r\in\roots} \rho_{r\wedge \s}$\cr
            $D_\s$            & $1$ if $b_\s\epsilon_\s$ odd, $2$ if $b_\s\epsilon_\s$ even\cr
            $m_\s$            & $(3-D_\s)b_\s$\cr
            $p_\s$            & $1$ if $|\s|$ is odd, $2$ if $|\s|$ is even\cr
            $s_\s$            & $\frac 12(|\s|\rho_\s+p_\s\rho_\s-\epsilon_\s)$\cr
            $\gamma_\s$       & $2$ if $|\s|$ is even and $\epsilon_\s\!-\!|\s|\rho_\s$ is odd, $1$ otherwise\cr
            $p_\s^0$          & $1$ if $\s$ is minimal and $\s\cap K\neq\varnothing$, $2$ otherwise\cr

            $s_\s^0$          & $-\epsilon_\s/2+\rho_\s$\cr
            $\gamma_\s^0$     & 2 if $p_\s^0=2$ and $\epsilon_\s$ is odd, 1 otherwise\cr
\hline
\end{tabular} 
\end{center}
\end{defn}

\begin{defn}[Definition \ref{yRegularDefinition}]\label{yRegularIntroductionDefinition}
We say that the hyperelliptic curve $C$ is \textit{$y$-regular} if either $p\neq 2$ or $D_\s=1$ for any $\s\in\mathring\Sigma_K$.
\end{defn}

\begin{defn}
Let $\s\in\mathring\Sigma_K$ and let $c\in\{0,\dots,b_\s-1\}$ such that $c\rho_\s-\frac{1}{b_\s}\in\Z$. Define
\[\tilde\s=\{\s'\in\Sigma_K\cup\{\varnothing\}\mid \s'<\s\text{ and }\tfrac{|\s'|}{b_\s}-c\epsilon_\s\notin 2\Z\},\]
where $\varnothing<\s$ if $\s$ is minimal and $p_\s^0=2$.

The \textit{genus $g(\s)$} of a rational cluster $\s\in\mathring\Sigma_K$ is defined as follows:
\begin{itemize}
    \item If $D_\s=1$, then $g(\s)=0$.
    \item If $D_\s=2$, then $2g(\s)+1$ or $2g(\s)+2$ equals
    $\smaller{\dfrac{|\s|-\sum_{\s'<\s}|\s'|}{b_\s}}+|\tilde \s|$.
\end{itemize}
\end{defn}

\begin{nt}
Let $\alpha,a,b\in\Z$, with $\alpha>0$, $a>b$, and fix $\frac{n_i}{d_i}\in\Q$ so that
    \[\alpha a=\frac{n_0}{d_0}>\frac{n_1}{d_1}>\ldots>\frac{n_r}{d_r}>\frac{n_{r+1}}{d_{r+1}}=\alpha b,\quad\text{with\scalebox{0.9}{$\quad\begin{vmatrix}n_i\!\!\!&n_{i+1}\cr d_i\!\!\!&d_{i+1}\cr\end{vmatrix}=1$}},\]
    and $r$ minimal. We write $\P^1(\alpha,a,b)$ for a chain of $\P^1$s (Notation \ref{ChainNotation}) of length $r$ and multiplicities $\alpha d_i,\dots,\alpha d_r$. Denote by $\P^1(\alpha,a)$ the chain $\P^1(\alpha,a,\lfloor\alpha a-1\rfloor/\alpha)$.
\end{nt}

The following theorem describes the special fibre of a regular model of $C$ with strict normal crossings.\footnote{In this paper a `normal crossings' divisor is not a `strict normal crossings' divisor in general (see e.g.\ \cite[Remark 9.1.7]{Liu}).} It follows from a more general result constructing a proper flat model of $C$ unconditionally (Theorem \ref{ConstructionProperModelGeneralCaseTheorem}). For the special fibre $\mathcal{C}^\mathrm{min}_s$ of the minimal regular model with normal crossings, the reader can refer to Theorem \ref{MinimalRegularSNCModelTheorem}, where we also describe a defining equation for all components of $\mathcal{C}_s^\mathrm{min}$ and discuss the Galois action (for general $K$). Finally, note that all these models are constructed in \S\ref{ConstructionModelsSection} by giving an explicit open affine cover (see \S\ref{ChartsSubsection}-\ref{GlueingSubsection}). 

\begin{thm}[Regular SNC model]\label{ModelIntroductionTheorem}
Suppose $C$ is $y$-regular and has almost rational cluster picture.
Then we can explicitly construct a regular model with strict normal crossings $\mathcal{C}/O_{K}$ of $C$ (\S\ref{ChartsSubsection}-\ref{GlueingSubsection}). Its special fibre $\mathcal{C}_s/k$ is given as follows.

\begin{enumerate}[label=(\arabic*)]
\item Every $\s\in\mathring\Sigma_K$ gives a $1$-dimensional closed subscheme $\Gamma_\s$ of multiplicity 
$m_\s$. If $\s$ is \"{u}bereven and $\epsilon_\s$ is even, then $\Gamma_\s$ is the disjoint union of $\Gamma_\s^{-}\simeq\P^1$ and $\Gamma_\s^{+}\simeq\P^1$, otherwise $\Gamma_\s$ is a smooth geometrically integral curve of genus $g(\s)$ (write $\Gamma_\s^{-}=\Gamma_\s^{+}=\Gamma_\s$ in this case). 

\item Every $\s\in\mathring\Sigma_K$ with $D_\s=1$ gives $(|\s|-\sum_{\s'\in\mathring\Sigma_K,\,\s'<\s}|\s'|+p_\s^0-2)/b_\s$ open-ended $\P^1$s of multiplicity $b_\s$ from $\Gamma_\s$.

\item Finally, for any $\s\in\mathring\Sigma_K$ draw the following chains of $\P^1$s:
\cellspacetoplimit4pt
\cellspacebottomlimit4pt
\begin{center}
\begin{tabular}{|Sc|Sc|Sc|Sc|}
\hline
Conditions & Chain & From & To\cr
\Xhline{1 pt}
$\s$ minimal & $\P^1(\gamma_\s^0,-s_\s^0)$ & $\Gamma_\s^-$ & open-ended\cr 
\hline
$\s$ minimal, $p_\s^0/\gamma_\s^0=2$ & $\P^1(\gamma_\s^0,-s_\s^0)$ & $\Gamma_\s^+$ & open-ended\cr
\hline
$\s\neq\roots$ & $\P^1(\gamma_\s,s_\s,s_\s-p_\s\cdot\frac{\rho_\s-\rho_{P(\s)}}{2})$ & $\Gamma_\s^-$ & $\Gamma_{P(\s)}^-$\cr
\hline
$\s\neq\roots$, $p_\s/\gamma_\s=2$ & $\P^1(\gamma_\s,s_\s,s_\s-p_\s\cdot\frac{\rho_\s-\rho_{P(\s)}}{2})$ & $\Gamma_\s^+$ & $\Gamma_{P(\s)}^+$\cr
\hline
$\s=\roots$ & $\P^1(\gamma_\s,s_\s)$ & $\Gamma_\s^-$ & open-ended\cr
\hline
$\s=\roots$, $p_\s/\gamma_\s=2$ & $\P^1(\gamma_\s,s_\s)$ & $\Gamma_\s^+$ & open-ended \cr
\hline
\end{tabular}
\end{center}

\end{enumerate}
\end{thm}

\begin{defn}
For any $\s\in\mathring\Sigma_K$, an element $w_\s\in K$ is called \textit{rational centre} of $\s$ if $\min_{r\in\s}v(r-w_\s)=\rho_\s$.
\end{defn}

If $\s'<\s$ and $w_{\s'}$ is a rational centre of $\s'$, then $w_{\s'}$ is also a rational centre of $\s$. For any minimal rational cluster $\s'$ fix a rational centre $w_{\s'}$. For any $\s\in\mathring\Sigma_K$ fix $w_\s=w_{\s'}$ for some minimal rational cluster $\s'\subseteq\s$.

The following result gives a basis of integral differentials when $K=K^{nr}$. In Theorem \ref{DifferentialsGeneraclCaseTheorem} we extend it to the case $K\neq K^{nr}$.

\begin{thm}[Theorem \ref{DifferentialTheorem}]\label{DifferetialsIntroductionTheorem}
Suppose $C$ is $y$-regular and has almost rational cluster picture. For $i=0,\dots,g-1$, inductively 
\begin{enumerate}[label=(\roman*)]
    \item define $e_i:=\displaystyle\max_{\t\in\mathring\Sigma_K}\bigg\{\sfrac{\epsilon_{\t}}{2}-\rho_\t-\smaller{\sum_{j=0}^{i-1}}\rho_{\s_j\wedge\t}\bigg\}$;
    \item choose clusters $\s_i\in\mathring\Sigma_K$ so that $e_i=\frac{\epsilon_{\s_i}}{2}-\sum_{j=0}^i\rho_{\s_j\wedge\s_i}$. If $\s$ and $\s'$ are two possible choices for $\s_i$ satisfying $\s'\subsetneq\s$, then choose $\s_i=\s$.
\end{enumerate}
Then a basis of integral differentials is given by
\[\mu_i=\pi^{\lfloor e_i\rfloor}\prod_{j=0}^{i-1}(x-w_{\s_j})\sfrac{dx}{2y},\qquad i=0,\dots,g-1.\]
\end{thm}

Note that given $e_i$ as in the previous theorem, the sum $\sum_{i=0}^{g-1}\lfloor e_i \rfloor$ is the quantity, often denoted by $v(\nicefrac{\omega^\circ}{\omega})$, appearing in the period in the Birch and Swinnerton-Dyer conjecture (for more details see \cite{FLSSSW}, \cite[\S 1.3]{Bom}).

\subsection{Rational cluster picture}\label{RationalClusterPictureSubsection}
In this subsection we define the rational cluster picture and compare it with the \textit{classical} cluster picture defined in \cite{D2M2}. We will show, via a simple example, in which sense the new object we introduce appears to be more suitable for the study of regular models.

\begin{defn}[Definition \ref{RationalClusterPictureDefinition}]\label{RationalClusterPictureIntroductionDefinition}
Let $K$ and $C$ as before. The \textit{rational cluster picture} of $C$ is the collection of its rational clusters $\Sigma_K$ together with their radii.
\end{defn}

\begin{exa}
Let $p$ be any prime number and set $K=\Q_p^{nr}$. Let $E_p/\Q_p^{nr}$ given by $y^2=x^3-p$. Then $E_p$ is an elliptic curve with Kodaira-N\'{e}ron reduction type II. Therefore the minimal regular model (with normal crossings) of $E_p$ does not depend on $p$. This is in accordance with the fact that the rational cluster picture of $E_p$ is the same for all $p$. Indeed, the set of roots of the polynomial $x^3-p$ is $\roots=\{\sqrt[3]{p}, \zeta_3\sqrt[3]{p}, \zeta_3^2\sqrt[3]{p}\}$, where $\zeta_3$ is a primitive $3$-rd of unity. Hence the rational cluster picture of $E_p$ is
\begin{center}
\begin{minipage}{2.4cm}\begin{center}\includegraphics[trim=9.8cm 22.33cm 9.8cm 4.6cm,clip]{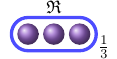}\end{center}\end{minipage} for any $p$,
\end{center}
where we denoted with bullet points the roots in $\roots$, with a surrounding oval the only rational cluster $\roots$, and with the subscript the radius $\rho_\roots$ of $\roots$.

A different behaviour is observed when we consider the cluster picture \cite[Definition 1.26]{D2M2} of $E_p$, collection of its clusters together with their depths.
The cluster picture of $E_p$ is
\begin{center}
\begin{tabular}{|c|c|c|}
\hline
    $p=2$ & $p=3$ & $p>3$ \cr
\hline
 $\substack{\text{cluster picture}\\ \text{not defined}}$ &      \begin{minipage}{2.4cm}\begin{center}\includegraphics[trim=9.8cm 22.33cm 9.8cm 4.6cm,clip]{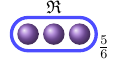}\end{center}\end{minipage} &     \begin{minipage}{2.4cm}\begin{center}\includegraphics[trim=9.8cm 22.33cm 9.8cm 4.6cm,clip]{Clusters/Cluster6.pdf}\end{center}\end{minipage}\\
\hline
\end{tabular}
\end{center}
where the subscripts represent the depth of the cluster $\roots$. It does depend on $p$ and differs from the rational cluster picture when $p=3$. Thus, although the cluster picture is particularly useful for Galois theoretical problems, the rational cluster picture appears to be a more suitable object for the study of regular models of the curve.

Finally, note that $E_p$ has an almost rational cluster picture. For any two distinct roots $r,r'\in\roots$, the smallest $v$-adic disc $D_{r,r'}$ containing them also contains the whole $\roots$. The element $0\in\Q_p^{nr}$ belongs to $D_{r,r'}$ when $p\neq 3$, while $|D_{r,r'}\cap\roots|=3=|v(r)|_p$, if $p=3$.
\end{exa}

The advantages of the rational cluster picture discussed in this subsection can also be observed in the following example where we study a more complex family of hyperelliptic curves having almost rational cluster picture. 

\subsection{Example}\label{ExampleSubsection}
In this subsection we are going to present an example of a family of hyperelliptic curves $C_p$ satisfying the hypothesis of Theorems \ref{ModelIntroductionTheorem} and \ref{DifferetialsIntroductionTheorem}. Via those results we will then describe the special fibre of the minimal regular model and a basis of integral differentials of $C_p$. All the computations involved are explained in detail in Examples \ref{AlmostRationalExample}, \ref{MinimalRegularModelExample} and \ref{IntegralDifferentialsExample}.

For any prime number $p$, let $a\in\Z_p$, $b\in\Z_p^\times$ such that the polynomial $x^2+ax+b$ is not a square modulo $p$. Let $C_p/\Q_p$ be the hyperelliptic curve of genus $4$ given by $y^2=f(x)$, where 
$f(x)=(x^6+ap^4x^3+bp^8)((x-p)^3-p^{11})$.
The curve $C_p/\Q_p^{nr}$ has an almost rational cluster picture and is $y$-regular when $p=2$. Its rational cluster picture is
\begin{center}
    \includegraphics[trim=6cm 22.2cm 6cm 4.5cm,clip]{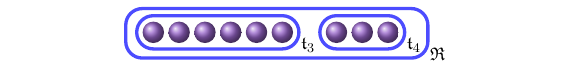}
\end{center}
where $\rho_{\t_3}=\frac{4}{3}$, $\rho_{\t_4}=\frac{11}{3}$, and $\rho_\roots=1$.
From Theorem \ref{ModelIntroductionTheorem} we can construct a regular model with strict normal crossings of $C_p$ with special fibre
\begin{center}
\pbox[c]{20cm}{
\begin{tikzpicture}[xscale=1,yscale=0.9,
  l1/.style={shorten >=-1.3em,shorten <=-0.5em,thick},
  l2/.style={shorten >=-0.3em,shorten <=-0.3em},
  lfnt/.style={font=\tiny},
  leftl/.style={left=-3pt,lfnt},
  rightl/.style={right=-3pt,lfnt},
  mainl/.style={scale=0.8,above left=-0.17em and -1.5em},
  mainleftl/.style={scale=0.8,above right=-0.17em and -1.5em},
  abovel/.style={above=-2.5pt,lfnt},
  facel/.style={scale=0.7,blue,below right=-0.5pt and 6pt},
  faceleftl/.style={scale=0.7,blue,below left=-0.5pt and 6pt},
  redbull/.style={red,label={[red,scale=0.6,above=-0.17]#1}}]
\draw[l1] (-.66,0.00)--(2.53,0.00) node[mainl] {2} node[facel] {$\Gamma_\roots$};
\draw[l2] (0.50,0.00)--node[rightl] {1} (0.50,0.66);
\draw[l1] (1.30,2.00)--(3.33,2.00) node[mainl] {6} node[facel] {$\Gamma_{\t_4}$};
\draw[l2] (1.96,0.00)--node[rightl] {1} (1.96,0.66);
\draw[l2] (1.30,0.66)--node[abovel] {2} (1.96,0.66);
\draw[l2] (1.30,0.66)--node[leftl] {3} (1.30,1.33);
\draw[l2] (1.30,1.33)--node[abovel] {4} (1.96,1.33);
\draw[l2] (1.96,1.33)--node[rightl] {5} (1.96,2.00);
\draw[l2] (1.30,2.00)--node[rightl] {3} (1.30,2.66);
\draw[l2] (2.76,2.00)--node[rightl] {4} (2.76,2.66);
\draw[l2] (2.10,2.66)--node[abovel] {2} (2.76,2.66);
\draw[l1] (-.60,0.66)--(-3.16,0.66) node[mainleftl] {6} node[faceleftl] {$\Gamma_{\t_3}$};
\draw[l2] (-.60,0.00)--node[rightl] {4} (-.60,0.66);
\draw[l2] (-1.20,0.66)--node[rightl] {3} (-1.20,1.33);
\draw[l2] (-2,0.66)--node[rightl] {3} (-2,1.33);
\draw[l2] (-2.70,0.66)--node[rightl] {2} (-2.70,1.33);
\end{tikzpicture}
}
\end{center}
over $\bar \F_p$. Computing the self-intersection of each irreducible component we easily see that this model coincides with the minimal regular model $\mathcal{C}^\mathrm{min}$. Theorem \ref{MinimalRegularSNCModelTheorem} also describes the action of the Galois group $\Gal(\bar \F_p/\F_p)$ on the special fibre $\mathcal{C}_s^\mathrm{min}$ of $\mathcal{C}^\mathrm{min}$. If the roots of $x^2+ax+b\mod p$ are in $\F_p$ then the absolute Galois group acts trivially on each component, otherwise it swaps the $2$ irreducible components of multiplicity $3$ intersecting $\Gamma_{\t_3}$.

From Theorem \ref{DifferetialsIntroductionTheorem} it follows that, for any $p$, a basis of integral differentials of $C_p/\Q_p^{nr}$ is given by
\[\mu_0=p^4\cdot\sfrac{dx}{2y},\quad\mu_1=p^3(x-p)\cdot\sfrac{dx}{2y},\quad\mu_2=p(x-p)x\cdot\sfrac{dx}{2y},\quad\mu_3=(x-p)x^2\cdot\sfrac{dx}{2y}.\]
In fact, this is also a basis of integral differentials of $C_p/\Q_p$ since they are all defined over $\Q_p$ (see Proposition \ref{DualisingSheafProposition}).


Below we will present related works of other authors concerning regular models and integral differentials of hyperelliptic curves. Note that the example presented here is not covered by \cite{D2M2} and \cite{Dok} since the curve $C_p$ is not semistable and not $\Delta_v$-regular. In fact, if $p=3$ the curve $C_p$ does not even have tamely potential semistable reduction. The results in \cite{FN} assume $p>2$ and $C_p$ with tamely potential semistable reduction, hence they can not be used when $p=2,3$. Finally, there is no classification for genus $4$ curves.

\subsection{Related works of other authors}\label{RelatedWorksSection}
Let $K$ be a
discretely valued field with residue field $k$ of characteristic $p$ and let $C/K$ be a hyperelliptic curve of genus $g$. 

In genus $1$, when $k$ is perfect, thanks to Tate's algorithm, one can describe the minimal regular model and the space of integral differentials of an elliptic curve $C$ (see for example \cite[IV.8.2]{Sil2}, \cite[Theorem 9.4.35]{Liu}).

If $K=\C(t)$ and $C$ has genus $2$, then Namikawa and Ueno \cite{NU} and Liu \cite{Liu94bis} give a full classification of the possible configurations of the special fibre of the minimal regular model of $C$. 


If $p\neq 2$, then Liu and Lorenzini show in \cite{LL} that regular models of $C$ can be seen as double cover of well-chosen regular models of $\P^1_K$. Since the latter can be found by using the MacLane valuations (\cite{Mac}) approach in \cite{OW}, this argument gives a way to describe any regular model of a hyperelliptic curve. At the moment there is no known closed form description of a regular model based on this approach and it has not been generalised to the $p=2$ case.

If $p>2$, $k$ finite, and $C$ is semistable, then in \cite{D2M2} the authors explicitly construct a minimal regular model in terms of the cluster picture of $C$. 
Under the same assumptions, Kunzweiler \cite{Kun} gives a basis of integral differentials rephrasing \cite[Proposition 5.5]{Kau} in terms of the cluster invariants introduced in \cite{D2M2}. These results can be recovered from Theorem \ref{MinimalRegularSNCModelTheorem} (see Corollary \ref{SemistableMinimalRegularModelCorollary}) and Theorem \ref{DifferentialTheorem}.

If $p>2$ and $C$ is semistable over some tamely ramified extension $L/K$, then Faraggi and Nowell \cite{FN} find the special fibre of the minimal regular model of $C$ with strict normal crossings taking the quotient of the stable model of $C_L$ and resolving the (tame) singularities. However, since they do not describe the charts of the model, their result does not immediately yield all arithmetic invariants, such as a basis of integral differentials.

The last work we want to recall represents an important ingredient of the strategy we will use in this paper (described more precisely in the next subsection). T.\ Dokchitser in \cite{Dok} shows that the toric resolution of $C$ gives a regular model in case of $\Delta_v$-regularity (\cite[Definition 3.9]{Dok}). This result, used also in \cite{FN}, holds for general curves and in any residue characteristic. In his paper, Dokchitser also describes a basis of integral differentials since his model is given as open cover of affine schemes. In Corollary \ref{RegularityAfterTranslationCorollary} and Theorem \ref{DifferentialsNestedTheorem}, we will rephrase his results for hyperelliptic curves by using rational cluster picture invariants from \S \ref{ClustersSection}.

\subsection{Strategy and outline of the paper}\label{StrategyOutlineSubsection}
In \cite{Dok}, Dokchitser not only describes a regular model of $C$ in case of $\Delta_v$-regularity, but also constructs a proper flat model $\mathcal{C}_\Delta$ without any assumptions on $C$. Assume $C$ is $y$-regular and has an almost rational cluster picture over $K^{nr}$ with rational centres $w_1,\dots, w_m\in K^{nr}$. Our approach to construct the minimal regular model with normal crossings of $C$ is composed by the following steps:
\begin{itemize}
    \item Consider the $x$-translated hyperelliptic curves $C^{w_h}/K^{nr}:y^2=f(x+w_h)$, for $h=1,\dots,m$. For each $h$, \cite[Theorem 3.14]{Dok} constructs a proper flat model $\mathcal{C}_\Delta^{w_h}$, possibly singular.
    \item We glue regular open subschemes of these models along common opens, and show that the result is a proper flat regular model $\mathcal{C}$ of $C_{K^{nr}}$ with strict normal crossings.
    \item We give a complete description of what components of the special fibre of $\mathcal{C}$ have to be blown down to obtain the minimal model with normal crossings $\mathcal{C}^\mathrm{min}$ of $C_{K^{nr}}$. 
    \item Finally, we describe the action of the absolute Galois group $G_k$ of $k$ on the special fibre of $\mathcal{C}^\mathrm{min}$.
\end{itemize}
We will explicitly describe both the models $\mathcal{C}_\Delta^{w_h}$ and $\mathcal{C}$.
This allows us to study the global sections of its relative dualising sheaf $\omega_{\mathcal{C}/O_K}(\mathcal{C})$.

In \S \ref{NewtonPolygonSection}, we present some results on Newton polygons used in the following sections. In \S \ref{ClustersSection}, we recall the basic objects and notation of \cite{D2M2} and define the rational cluster picture. Moreover, we relate it with the notions given in \S \ref{NewtonPolygonSection}. 
This comparison allows us to rephrase the objects in \cite{Dok} in terms of rational clusters invariants in \S \ref{DescriptionModelsSection}. In the same section we also state the theorems which describe the special fibres of a proper flat model (Theorem \ref{ConstructionProperModelGeneralCaseTheorem})
and of the minimal regular model with normal crossings (Theorem \ref{MinimalRegularSNCModelTheorem}) of $C$. The construction of these models, from which the two theorems above follow, is presented in \S \ref{ConstructionModelsSection}. Finally, in \S \ref{IntegralDifferentialsSection}, Theorems \ref{DifferentialTheorem} and \ref{DifferentialsGeneraclCaseTheorem} describe a basis of integral differentials of $C$, in terms of rational clusters invariants defined in \S \ref{ClustersSection}.

\subsection*{Acknowledgements} The author would like to thank his supervisor Tim Dokchitser for the very useful conversations, corrections and general advice.

\newpage
\subsection{Notation}
In the following, we present the main notation used for fields, hyperelliptic curves and Newton polytopes.
\begin{center}
\begin{tabular}{l@{$\>\quad$}l}
$K,v$ & complete field with normalised discrete valuation $v$\cr
$O_k,\pi,k,p$ & ring of integers, uniformiser, residue field, $\mathrm{char}(k)$\cr
$\bar K,\bar k$ & fixed algebraic closure of $K$, residue field of $\bar K$\cr
$K^\mathrm{s},k^\mathrm{s}$ & separable closure of $K$ in $\bar K$, residue field of $K^\mathrm{s}$\cr
$K^{nr}$ & maximal unramified extension of $K$ in $K^\mathrm{s}$\cr
$F$ & extension of $K$ in $\bar K$, unramified in \S\ref{DescriptionModelsSection}\cr
$G_K, G_k$ & absolute Galois groups $\Gal(K^\mathrm{s}/K), \Gal(k^\mathrm{s}/k)$ \cr 
$f(x)$ & $=\sum a_i x^i$, polynomial in $K[x]$, separable from \S\ref{ClustersSection}\cr
$\Np(f)$ & Newton polygon of $f$, lower convex hull of $\{(i,v(a_i))\mid i\}$\cr
$f|_L,\overline{f|_L}$ & restriction and reduction of $f$ to an edge $L$ of $\Np(f)$ (\ref{RestrictionReductionDefinition})\cr
$g(x,y)$ & $=y^2-f(x)$, polynomial in $K[x,y]$ defining $C$ \cr
$C$ & hyperelliptic curve defined over $K$ by $g(x,y)=0$\cr
$f_w(x), f_h(x)$ & $=f(x+w), f(x+w_h)$, for a given rational centre $w_h$\cr
$g_w(x,y), g_h(x,y)$ & $=y^2-f_w(x), y^2-f_h(x)$\cr
$C^w$ & $\simeq C$, hyperelliptic curve given by $g_w(x,y)=0$\cr
$\Delta^w, \Delta_v^w$ & Newton polytopes attached to $C^w$ as in \cite[\S1.1]{Dok}\cr
$\mathcal{C}_\Delta^{w}$ & proper model of $C^w$ constructed from $\Delta_v^w$ by \cite[3.14]{Dok}\cr
$F_\t^w,L_\t^w,V_\t^w,V_0^w$ & $v$-faces and $v$-edges of $\Delta^w$ (\ref{FacesEdgesClustersNotation})
\end{tabular}
\end{center}

For a separable polynomial $f\in k[x]$ or a hyperelliptic curve $C/K:y^2=f(x)$ as above, the following is the main notation for clusters.
\begin{center}
\begin{tabular}{l@{$\>\quad$}l}
$c_f,\roots$ & leading coefficient and set of roots of $f$\cr
$\Sigma_f,\Sigma_C$ & cluster picture, the set of clusters of $f$,$C$ (\ref{ClusterPictureDefinition})\cr
$\s\in\Sigma_C$ & cluster, $\s=\mathcal{D}\cap\roots$, for a $v$-adic disc $\mathcal{D}$ (\ref{ClusterDepthDefinition})\cr
 $G_\s, K_\s, k_\s$ & $G_\s=\mathrm{Stab}_{G_K}(\s)$; $K_\s=\lb K^\mathrm{s}\rb^{G_\s}$; $k_\s$ residue field of $K_\s$ \cr
$d_\s$ & $=\min_{r,r'\in\s}v(r-r')$ is the depth of a cluster $\s$ (\ref{ClusterDepthDefinition})\cr
$\s'<\s=P(\s')$ & $\s'$ is a child of $\s$ and $\s$ is the parent of $\s'$ (\ref{ParentChildWedgeDefinition})\cr
$\s\wedge\t$ & smallest cluster containing $\s$ and $\t$ (\ref{ParentChildWedgeDefinition})\cr
$\rho_\s$     & $=\max_{w\in F}\min_{r\in\s} v(r-w)$, radius of $\s\in\Sigma_{C_F}$ (\ref{RadiusRationalCentreDefinition}, \ref{QuantitiesForTheoremsOnModelsDefinition})\cr
$b_\s$            & denominator of $\rho_\s$ (\ref{QuantitiesForTheoremsOnModelsDefinition})\cr
$w_\s$      & rational centre of $\s$ (\ref{RadiusRationalCentreDefinition})\cr
$\epsilon_\s$     & $=v(c_f) + \sum_{r\in\roots} \rho_{r\wedge \s}$ (\ref{EpsilonDefinition}, \ref{QuantitiesForTheoremsOnModelsDefinition})\cr
$\Sigma_f^\mathrm{rat},\Sigma_C^\mathrm{rat}$ & rational cluster picture (\ref{RationalClusterPictureDefinition})\cr
$\s\in\Sigma_C^\mathrm{rat}$ & rational cluster (\ref{RationalClusterPictureDefinition})\cr
$\Sigma_F$ & $=\Sigma_{C_F}^\mathrm{rat}$, for some extension $F/K$ (\ref{QuantitiesForTheoremsOnModelsDefinition})\cr
$\Sigma_f^z,\Sigma_C^z$ & cluster picture centred at $z$ (\ref{CentredClusterPictureDefinition})\cr
$\s\in\Sigma_C^z$ & cluster centred at $z$ (\ref{CentredClusterDefinition})\cr
$\rho_\s^z,\epsilon_\s^z$ & $\rho_\s^z=\min_{r\in\s}v(r-z)$, $\epsilon_\s^z=v(c_f)+\sum_{r\in\roots}\rho_{r\wedge\s}^z$ (\ref{RadiusWithRespectToCentreDefinition})\cr
$\Sigma^W$, $\mathring\Sigma$ & $\Sigma^W=\bigcup_{w\in W}\Sigma_C^w$, $\mathring\Sigma\subset\Sigma_{K^{nr}}$ non-removable clusters (\ref{RemovableContractibleDefinition})\cr
$w_{hl}$ & $=w_h-w_l$ for fixed rational centres $w_h,w_l$ (\S\ref{ChartsSubsection})\cr
$u_{hl},\rho_{hl}$ & $u_{hl}\in O_K^\times$, $\rho_{hl}\in \Z$ such that $w_{hl}=u_{hl}\pi^{\rho_{hl}}$; $u_{hh}=0$ (\S\ref{ChartsSubsection})\cr
$D_\s,m_\s$            & $D_\s=1$ if $b_\s\epsilon_\s$ odd, $2$ if $b_\s\epsilon_\s$ even; $m_\s=(3-D_\s)b_\s$ (\ref{QuantitiesForTheoremsOnModelsDefinition})\cr
$p_\s$            & $=1$ if $|\s|$ is odd, $2$ if $|\s|$ is even (\ref{QuantitiesForTheoremsOnModelsDefinition})\cr
 $\gamma_\s$       & $=2$ if $|\s|$ is even and $\epsilon_\s\!-\!|\s|\rho_\s$ is odd, $1$ otherwise (\ref{QuantitiesForTheoremsOnModelsDefinition})\cr
 $p_\s^0$          & $=1$ if $\s$ is minimal and $\s\cap K_\s\neq\varnothing$, $2$ otherwise (\ref{QuantitiesForTheoremsOnModelsDefinition})\cr
$\gamma_\s^0$     & $=2$ if $p_\s^0=2$ and $\epsilon_\s$ is odd, 1 otherwise (\ref{QuantitiesForTheoremsOnModelsDefinition})\cr
 $s_\s$, $s_\s^0$           & $s_\s=\frac 12(|\s|\rho_\s+p_\s\rho_\s-\epsilon_\s)$, $s_\s^0=-\epsilon_\s/2+\rho_\s$ (\ref{QuantitiesForTheoremsOnModelsDefinition})\cr
 $\ch{g_\s}, \ch{g_\s^0}, \ch{f_\s^W}, \ch{f_\s}, \tilde f_\s$ & polynomials in one variable over $k_\s$ (\ref{SchemesXsDefinition}, \ref{PolynomialsDefinition}) 
\end{tabular}
\end{center}


\section{Newton polygon}\label{NewtonPolygonSection}
Let $K$ be a complete field with a discrete valuation $v$, ring of integers $O_K$, uniformiser $\pi$, and residue field $ k$ of characteristic $p$.
We fix $\bar K$, an algebraic closure of $K$, of residue field $\bar k$, and we denote by $K^\mathrm{s}$ the separable closure of $K$ in $\bar K$, and by $ k^\mathrm{s}$ its residue field. Note that $ k^\mathrm{s}$ is the separable closure of $k$ in $\bar k$. Extend the valuation $v$ to $\bar K$. Write $G_K$, $G_k$ for the Galois groups $\Gal(K^\mathrm{s}/K)$, $\Gal( k^\mathrm{s}/ k)$, respectively. Finally, denote by $K^{nr}$ the maximal unramified extension of $K$ in $K^\mathrm{s}$.
\begin{nt}
Let $O_{\bar K}=\{a\in\bar K\mid v(a)\geq 0\}$. Throughout this paper, given an element $a\in O_{\bar K}$, we will write $a\mod\pi$ for the reduction of $a$ in $\bar k$. Similarly, given a polynomial $h\in O_{\bar K}[x_1,\dots,x_n]$, namely $h=\sum a_{i_1,\dots,i_n}\cdot x_1^{i_1}\cdots x_n^{i_n}$, we will write $h\mod\pi$ for the polynomial $\sum (a_{i_1,\dots,i_n}\mod\pi)\cdot x_1^{i_1}\cdots x_n^{i_n}\in\bar k[x_1,\dots,x_n]$.
\end{nt}

Let $f\in K[x]$ be a non-zero polynomial of degree $d$, say
\[ f(x)=\sum_{i=0}^da_ix^i.\]
The \textit{Newton polygon} of $f$, denoted $\NP f$, is 
\[\NP f = \mbox{lower convex hull}~\left\{(i,v(a_i))\mid\, i=0,\dots, d, \, a_i\neq 0\right\}\subset \R^2.\]

We recall the following well-known result (see for example \cite[II.6.3,6.4]{Neu}).

\begin{thm}
\label{NewtonpolygonFactorisation(Theorem)}
Let $i_0<\ldots<i_s=d$ be the set of indices in $\{0,\dots,d\}$ such that the points $(i_0,v(a_{i_0})),\dots,(i_s,v(a_{i_s}))$ are the vertices of $\NP f$. For any $j=1,\dots,s$, denote by $\rho_j$ the slope of the edge of $\NP f$ which links the points $(i_{j-1},v(a_{i_{j-1}}))$ and $(i_j,v(a_{i_j}))$. Then $f$ has a unique factorisation over $K$ as a product
\[f=a_d \cdot g_0 \cdot g_1\cdots g_s,\]
where $g_0=x^{i_0}$ and, for all $j=1,\dots, s$, 
\begin{itemize}
    \item the polynomials $g_j\in K[x]$ are monic of degree $d_j=i_j-i_{j-1}$,
    \item all the roots of $g_j$ have valuation $-\rho_j$ in $\bar K$.
\end{itemize}
In particular, $\NP{g_j}$ is a segment of slope $-\rho_j$.
\end{thm}


\begin{cor}\label{ValuationofRoots(Corollary)}
With the notation of Theorem \ref{NewtonpolygonFactorisation(Theorem)}, the polynomial $f$ has exactly $d_j$ roots of valuation $-\rho_j$ for all $j=1,\dots,s$.
\end{cor}

\begin{cor}
If $f=\sum a_ix^i$ is irreducible of degree $d$ and $a_0\neq 0$, then $\NP f$ is a segment linking the points $(0,v(a_0))$ and $(d,v(a_d))$.
\end{cor}

\begin{defn}[Restriction and reduction]\label{RestrictionReductionDefinition}
Let $f=\sum_{i=0}^da_ix^i\in K[x]$ and consider an edge $L$ of its Newton polygon $\NP f$. Let $(i_1,v(a_{i_1})), (i_2,v(a_{i_2}))$, $i_1<i_2$ be the two endpoints of $L$. Denote by $\rho$ the slope of $L$ and by $n$ the denominator of $\rho$.
Define the \textit{restriction} of $f$ to $L$ as
\[f|_L:=\sum_{i= 0}^{(i_2-i_1)/n}a_{ni+i_1}x^{i}\in K[x].\]

Moreover we define the \textit{reduction} of $f$ with respect to $L$ to be the polynomial
\[\ch{f|_L}:=\pi^{-c}f|_L(\pi^{-n\rho} x)\mbox{ mod }\pi\in  k[x],\]
where $c=v(a_{i_1})=v(a_{i_2})+(i_1-i_2)\rho.$
\end{defn}

\begin{rem}
These definitions coincide with the ones given in \cite[Definitions 3.4, 3.5]{Dok} when the number of variables is $1$ (for suitable choices of basis of the lattices used in the definitions).
\end{rem}

Until the end of the section let $f\in K[x]$, consider a factorisation $f=a_d\cdot g_0\cdot g_1\cdots g_s$ as in Theorem \ref{NewtonpolygonFactorisation(Theorem)}. Denote by $L_j$ the edge of slope $-\rho_j$ of $\NP f$, for any $j=1\dots s$.

\begin{rem}\label{Remarkfor(Lemma1)}
By the lower convexity of $\NP f$, for all $j=1,\dots,s$, note that $\ch{f|_{L_j}}=\bar c_j\cdot \ch{g_j|_{\NP{g_j}}}$ for some $\bar c_j\in k^\times$. In particular they define the same $ k$-scheme in $\mathbb{G}_{m, k}$. More precisely, for any $j=1,\dots,s$, let
\[u_j=a_d\cdot \prod_{i=j+1}^{s}g_i(0).\]
Then $\bar c_j=u_j/\pi^{v(u_j)}\mod \pi$.

\end{rem}

\begin{defn}
We say that $f$ is \textit{$\Np$-regular} if the $ k$-scheme
\[X_{L_j}:\{\ch{f|_{L_j}}=0\}\subset\mathbb{G}_{m, k}\]
is smooth for all $j=1,\dots,s$.
\end{defn}

\begin{lem}\label{regularityfactorisation(Lemma)}
The polynomial $f=a_d\cdot g_0\cdot g_1\cdots g_s$ is $\Np$-regular if and only if $g_j$ is $\Np$-regular for every $j=1,\dots,s$.
\proof
The lemma follows from Remark \ref{Remarkfor(Lemma1)}.
\endproof
\end{lem}

We conclude this section with two examples.

\begin{exa}
Let $f=x^{11}+9x^7-3x^6+9x^5+81x-27\in\Q_3[x]$. Then the Newton polygon of $f$ is
\begin{center}
\begin{tikzpicture}[scale=0.8]
    \draw[very thin,color=gray] (-0.1,-0.1) grid (11.3,4.3);    
    \draw[->] (-0.2,0) -- (11.4,0) node[right] {$i$};
    \draw[->] (0,-0.2) -- (0,4.4) node[above] {$v(a_i)$};
    \tkzDefPoint(0,3){A}
    \tkzDefPoint(1,4){B}
    \tkzDefPoint(6,1){C}
    \tkzDefPoint(7,2){D}
    \tkzDefPoint(5,2){E}
    \tkzDefPoint(11,0){F}
    \tkzLabelPoint[right,below](F){$(11,0)$}
    \tkzLabelPoint[right,below](C){$(6,1)$}
    \tkzLabelPoint[left](A){$(0,3)$}
    \foreach \n in {A,B,C,D,E,F}
    \node at (\n)[circle,fill,inner sep=1.5pt]{};
    \tkzDrawSegment[black!60!black](A,C)
    \tkzDrawSegment[black!60!black](C,F)
    \tkzLabelSegment[sloped](A,C){$\rho_1=-\frac{1}{3}$}
    \tkzLabelSegment[sloped, below](A,C){$L_1$}
    \tkzLabelSegment[sloped](C,F){$\rho_2=-\frac{1}{5}$}
    \tkzLabelSegment[sloped,below](C,F){$L_2$}
\end{tikzpicture}
\end{center}
Corollary \ref{ValuationofRoots(Corollary)} implies that $f$ has $6$ roots of valuation $\frac{1}{3}$ and $5$ roots of valuation $\frac{1}{5}$. Furthermore, the two polynomials $g_1$ and $g_2$ in the factorisation $f=g_1\cdot g_2$ of Theorem \ref{NewtonpolygonFactorisation(Theorem)} turn out to be 
\[g_1=x^6+9,\qquad g_2=x^5+9x-3.\]
Finally, 
\[f|_{L_1}=-3x^2-27=-3\cdot g_1|_{\NP{g_1}},\qquad f|_{L_2}=x-3= g_2|_{\NP{g_2}};\]
and
\[\ch{f|_{L_1}}=-x^2-1=-(x^2+1)=-\ch{g_1|_{\NP{g_1}}},\qquad \ch{f|_{L_2}}=x-1=\ch{g_2|_{\NP{g_2}}}\qquad\mbox{in }\F_3[x].\]
Thus $f$ is $\Np$-regular.
\end{exa}

\begin{exa}
We now show an example of a polynomial that is not $\Np$-regular.
Let $f=x^9+12x^6+36x^3+81\in\Q_3[x]$. Then the Newton polygon of $f$ is
\begin{center}
\begin{tikzpicture}[scale=0.8]
    \draw[very thin,color=gray] (-0.1,-0.1) grid (9.3,4.3);
    \draw[->] (-0.2,0) -- (9.4,0) node[right] {$i$};
    \draw[->] (0,-0.2) -- (0,4.4) node[above] {$v(a_i)$};
    \tkzDefPoint(0,4){A}
    \tkzDefPoint(3,2){B}
    \tkzDefPoint(6,1){C}
    \tkzDefPoint(9,0){D}
    \tkzLabelPoint[right,below](D){$(9,0)$}
    \tkzLabelPoint[right,below](B){$(3,2)$}
    \tkzLabelPoint[left](A){$(0,4)$}
    \foreach \n in {A,B,C,D}
    \node at (\n)[circle,fill,inner sep=1.5pt]{};
    \tkzDrawSegment[black!60!black](A,B)
    \tkzDrawSegment[black!60!black](B,D)
    \tkzLabelSegment[sloped](A,B){$\rho_1=-\frac{2}{3}$}
    \tkzLabelSegment[sloped, below](A,B){$L_1$}
    \tkzLabelSegment[sloped](B,D){$\rho_2=-\frac{1}{3}$}
    \tkzLabelSegment[sloped,below](B,D){$L_2$}
\end{tikzpicture}
\end{center}
Corollary \ref{ValuationofRoots(Corollary)} implies that $f$ has $3$ roots of valuation $\frac{2}{3}$ and $6$ roots of valuation $\frac{1}{3}$. 
Furthermore, the two polynomials $g_1$ and $g_2$ in the factorisation $f=g_1\cdot g_2$ of Theorem \ref{NewtonpolygonFactorisation(Theorem)} 
are
\[g_1=x^3+9,\qquad g_2=x^6+3x^3+9.\]
Finally, 
\[f|_{L_1}=36x+81\qquad f|_{L_2}=x^2+12x+36,\]\[ g_1|_{\NP{g_1}}=x+9,\qquad g_2|_{\NP{g_2}}=x^2+3x+9;\]
and
\[\ch{f|_{L_1}}=x+1=\ch{g_1|_{\NP{g_1}}},\qquad \ch{f|_{L_2}}=(x+2)^2=\ch{g_2|_{\NP{g_2}}}\qquad\mbox{in }\F_3[x].\]
Then $f$ is not $\Np$-regular. In fact, according to Lemma \ref{regularityfactorisation(Lemma)}, $g_2$ is not $\Np$-regular.
\end{exa}

\section{Rational clusters}\label{ClustersSection}
From now on, let $f\in K[x]$ be a separable polynomial and denote by $\roots$ the set of its roots in $K^\mathrm{s}$ and by $c_f$ its leading coefficient. Then
\[f(x)=c_f\prod_{r\in\roots}(x-r).\]

\begin{defn}[{\cite[Definition 1.1]{D2M2}}]\label{ClusterDepthDefinition}
A \textit{cluster} (for $f$) is a non-empty subset $\s\subseteq\roots$ of the form $\mathcal{D}\cap\roots$, where $\mathcal{D}$ is a $v$-adic disc $\mathcal{D}=\{x\in\bar K\mid v(x-z)\geq d\}$ for some $z\in\bar K$ and $d\in\Q$. If $|\s|>1$ we say that $\s$ is \textit{proper} and define its \textit{depth} $d_\s$ to be
\[d_\s=\min_{r,r'\in\s}v(r-r').\]
Note that every proper cluster is cut out by a disc of the form \[\mathcal{D}=\{x\in\bar K\mid v(x-r)\geq d_\s\}\] for any $r\in\s$. 
\end{defn}

\begin{defn}[{\cite[Definition 1.26]{D2M2}}]\label{ClusterPictureDefinition}
The \textit{cluster picture} of $f$ is the collection of its clusters, together with their depths.

We denote by $\Sigma_f$ the set of all clusters of $f$ and by $\mathring{\Sigma}_f$ the subset of $\Sigma_f$ of proper clusters.
\end{defn}

\begin{defn}[{\cite[Definition 1.3]{D2M2}}]\label{ParentChildWedgeDefinition}
If $\s'\subsetneq\s$ is maximal subcluster, then we say that $\s'$ is a \textit{child} of $\s$ and $\s$ is the \textit{parent} of $\s'$, and we write $\s'<\s$. Since every cluster $\s\neq \roots$ has one and only one parent we write $P(\s)$ to refer to the unique parent of $\s$.

We say that a proper cluster $\s$ is \textit{minimal} if it does not have any proper child.
    
For two clusters (or roots) $\s_1, \s_2$, we write $\s_1\wedge\s_2$ for the smallest cluster that contains them.
\end{defn}

\begin{defn}[{\cite[Definition 1.4]{D2M2}}]\label{OddEvenUberevenClusterDefinition}
A cluster $\s$ is \textit{odd/even} if its size is odd/even. If $|\s|=2$, then we say $\s$ is a \textit{twin}. A cluster $\s$ is \textit{\"{u}bereven} if it has only even children.
\end{defn}

\begin{defn}[{\cite[Definition 1.9]{D2M2}}]\label{CentreDefinition}
A \textit{centre} $z_\s$ of a proper cluster $\s$ is any element $z_\s\in K^\mathrm{s}$ such that $\s=\mathcal{D}\cap\roots$, where
\[\mathcal{D}=\{x\in\bar K\mid v(x-z_\s)\geq d_\s\}.\]
Equivalently, $v(r-z_\s)\geq d_\s$ for all $r\in\s$. The \textit{centre} of a non-proper cluster $\s=\{r\}$ is $r$.
\end{defn}

\begin{defn}[{\cite[Definition 1.6]{D2M2}}]
For a proper cluster $\s$ set
\[\nu_\s:=v(c_f)+\sum_{r\in\roots}d_{r\wedge\s}.\]
\end{defn}

\begin{defn}
We say that $\Sigma_f$ is \textit{nested} if one of the following equivalent conditions is satisfied:
\begin{enumerate}
    \item [(i)] there exists $z\in  K^\mathrm{s}$ such that $z$ is a centre for all proper clusters $\s\in\Sigma_f$;
    \item [(ii)] there is only one minimal cluster in $\Sigma_f$;
    \item [(iii)] every non-minimal proper cluster has exactly one proper child.
\end{enumerate}
\end{defn}

\begin{defn}\label{RadiusRationalCentreDefinition}
A \textit{rational centre} of a cluster $\s$ is any element $w_\s\in K$ such that
\[\min_{r\in\s}v(r-w_\s)=\max_{w\in K}\min_{r\in\s}v(r-w).\]
If $\s=\{r\}$, with $r\in K$, then $w_\s=r$.

If $w_\s$ is a rational centre of a proper cluster $\s$, we define the \textit{radius} of $\s$ to be
\[\rho_\s=\min_{r\in\s}v(r-w_\s).\]
\end{defn}

\begin{defn}\label{RationalClusterPictureDefinition}
A \textit{rational cluster} is a cluster cut out by a $v$-adic disc of the form $\mathcal{D}=\{x\in\bar K\mid v(x-w)\geq d\}$ with $w\in K$ and $d\in\Q$.

The \textit{rational cluster picture} is the collection of all rational clusters of $f$ together with their radii.

We denote by $\Sigma_f^\mathrm{rat}\subseteq\Sigma_f$ the set of rational clusters and by $\mathring{\Sigma}_f^\mathrm{rat}$ the subset of $\Sigma_f^\mathrm{rat}$ of proper rational clusters.
\end{defn}

\begin{lem}
Let $\s$ be a proper cluster. Then $d_\s\geq\rho_\s$.
\proof
First we want to show that
\[\min_{r,r'\in\s}v(r-r')=\max_{z\in K^\mathrm{s}}\min_{r\in \s}v(r-z).\]
Clearly $\min_{r,r'\in\s}v(r-r')\leq\max_{z\in K^\mathrm{s}}\min_{r\in \s}v(r-z)$. Let $z_\s\in K^\mathrm{s}$ such that
\[\max_{z\in K^\mathrm{s}}\min_{r\in \s}v(r-z)=\min_{r\in \s}v(r- z_\s).\]
Then, for any $r,r'\in\s$, one has
\[v(r-r')\geq \min\{v(r-z_\s), v(r'-z_\s)\}\geq \min_{r\in \s}v(r- z_\s),\]
and so
\[\min_{r,r'\in\s}v(r-r')\geq\max_{z\in K^\mathrm{s}}\min_{r\in \s}v(r-z),\]
as required.
From
\[d_\s=\min_{r,r'\in\s}v(r-r')=\max_{z\in K^\mathrm{s}}\min_{r\in \s}v(r-z)\geq\max_{w\in K}\min_{r\in \s}v(r-w)=\rho_\s,\]
the lemma follows.
\endproof
\end{lem}

\begin{defn}\label{RationalisationDefinition}
Given a proper cluster $\s\in\Sigma_f$, we define the \textit{rationalisation} $\s^\mathrm{rat}$ of $\s$ to be the smallest rational cluster containing $\s$. By definition
\[\s^\mathrm{rat}=\roots\cap\{x\in\bar K\mid v(x-w_\s)\geq\rho_\s\},\]
where $w_\s$ is a rational centre of $\s$ and $\rho_\s$ is its radius.
\end{defn}

\begin{lem}\label{RadiusTimesSizeIntegerLemma}
Let $\s\in\Sigma_f^\mathrm{rat}$ be a proper cluster with rational centre $w_\s$. Let $\s'\in\Sigma_C^\mathrm{rat}$ be the child of $\s$ with rational centre $w_\s$ (let $\s'=\varnothing$ if it does not exist). Then $(|\s|-|\s'|)\rho_\s\in\Z$.
\proof
As $\s\in\Sigma_f^\mathrm{rat}$, one has $\s=\s^\mathrm{rat}$. Let $b_\s$ be the denominator of $\rho_\s$. Then $b_\s$ divides the degree of the minimal polynomial of $r$, for any $r\in\s$ satisfying $v(w_\s-r)=\rho_\s$. Then $(|\s|-|\s'|)\rho_\s\in\Z$, where
\[\s'=\roots\cap\{x\in\bar K\mid v(x-w_\s)>\rho_\s\},\]
as required.
\endproof
\end{lem}

\begin{rem}\label{TameNestedClusterPictureRemark}
If a proper cluster $\s\in\Sigma_f$ satisfies $d_\s=\rho_\s$, then a rational centre $w_\s\in K$ of its is also a centre. Hence $\s$ is a rational cluster and, in particular, is $G_K$-invariant. On the other hand, if a proper cluster $\s\in\Sigma_f$ is $G_K$-invariant and $K(\s)/K$ is tamely ramified, then $\s$ has a centre $z_\s\in K$ by \cite[Lemma B.1]{D2M2}. Thus $\rho_\s=d_\s$ and $\s\in\Sigma_f^\mathrm{rat}$.  
\end{rem}

\begin{lem}\label{rationalCentresRadiusChildParentLemma}
Let $\s$ be a proper cluster with rational centre $w_\s$ and let $\t\in\Sigma_f$ satisfying $\t\supseteq \s$. Then $w_\s$ is a rational centre of $\t$ and $\rho_{\t}\leq\rho_\s$. Furthermore, if $\s$ is a rational cluster and  $\t\supsetneq \s$, then $\rho_\t<\rho_\s$.
\proof
It suffices to prove the lemma for $\t=P(\s)$. Hence we first want to show that $\min_{r\in P(\s)}v(r-w_\s)=\rho_{P(\s)}$ and $\rho_{P(\s)}\leq\rho_\s$. Note that
\[\min_{r\in P(\s)}v(r-w_\s)\leq\max_{w\in K}\min_{r\in P(\s)}v(r-w)=\rho_{P(\s)}.\]
Moreover,
\[\rho_{P(\s)}=\max_{w\in K}\min_{r\in P(\s)}v(r-w)\leq\max_{w\in K}\min_{r\in\s}v(r-w)=\rho_\s.\]
If $r\in\s$ then $v(w_\s-r)\geq\rho_\s$, by definition of $\rho_\s$. On the other hand, if $r\in P(\s)\smallsetminus\s$ then fixing $r'\in\s$ we have
\[v(r-w_\s)=v(r-r'+r'-w_\s)\geq\min\{v(r-r'),v(r'-w_\s)\}\geq\min\{d_{P(\s)},\rho_\s\}\geq\rho_{P(\s)},\]
by the previous lemma. Thus $\min_{r\in P(\s)}v(r-w_\s)=\rho_{P(\s)}$, as required.

Now suppose $\s\in\Sigma_f^\mathrm{rat}$ with $\t\supsetneq \s$. From Definition \ref{RadiusRationalCentreDefinition}, it follows that
\[\{x\in \bar K\mid v(x-w_\s)\geq\rho_\s\}\cap \roots=\s\subsetneq\t\subseteq\{x\in \bar K\mid v(x-w_\s)\geq\rho_\t\}\cap \roots,\]
as $w_\s$ is a rational centre of $\t$. Thus $\rho_\t<\rho_\s$.
\endproof
\end{lem}

\begin{lem}\label{NoRationalClusterNoRationalChildLemma}
Every cluster $\s$ with $\rho_\s<d_\s$ has no rational subcluster $\s'\subsetneq\s$. 
\proof
Suppose by contradiction there exists $\s'\in\Sigma_C^\mathrm{rat}$, $\s'\subsetneq\s$, and fix a rational centre $w_{\s'}$ of $\s'$. Then $w_{\s'}$ is a rational centre of $\s$ by the previous lemma.  If $|\s'|=1$, then $w_{\s'}$ is also a centre of $\s$ and this contradicts $\rho_\s<d_\s$; so assume $\s'$ proper. Let $r'\in\s'$ such that $v(r'-w_{\s'})=\rho_{\s'}$ and $r\in\s$ such that $v(r-w_{\s'})=\rho_\s$. But then $d_\s\leq v(r-w_{\s'}+w_{\s'}-r')=\rho_\s$ again by Lemma \ref{rationalCentresRadiusChildParentLemma}.
\endproof
\end{lem}

In particular, the lemma above shows that if $\s\in\Sigma_f$ and $\s'\in\Sigma_f^\mathrm{rat}$ is a maximal rational subcluster of $\s$, with $\s'\subsetneq\s$, then $\s'$ is a child of $\s$. Moreover, the parent of a rational cluster is rational.

\begin{defn}\label{RationallyminimalDefinition}
We say that a proper rational cluster $\s\in\Sigma_f^\mathrm{rat}$ is \textit{(rationally) minimal} if it does not have any proper rational child.
\end{defn}


\begin{lem}\label{DifferentialsLemma}
Let $\s,\s'\in\Sigma_f^\mathrm{rat}$ such that $\s'\nsubseteq\s$. If $w_\s$ is a rational centre of $\s$ then \[\min_{r\in\s'}v(r-w_\s)=\rho_{\s\wedge\s'}.\]
\proof
By Lemma \ref{rationalCentresRadiusChildParentLemma} we have \[\min_{r\in\s\wedge\s'}v(r-w_\s)=\rho_{\s\wedge\s'}.\]
Therefore $\min_{r\in\s'}v(w_\s-r)\geq\rho_{\s\wedge\s'}$. Suppose by contradiction that \[\min_{r\in\s'}v(r-w_\s)=:\rho>\rho_{\s\wedge\s'}.\] 
It follows from Lemma \ref{rationalCentresRadiusChildParentLemma} that \[\min_{r\in\s}v(r-w_\s)=\rho_\s>\rho_{\s\wedge\s'}\] as $\s'\nsubseteq\s$. But then there exists $ \tilde r\in(\s\wedge\s')\smallsetminus(\s\cup\s')$ such that $v(\tilde r-w_\s)=\rho_{\s\wedge\s'}$. Consider the rational cluster
\[\t:=\roots\cap\left\{x\in\bar K\mid v(x-w_\s)\geq\min\{\rho,\rho_\s\}\right\}\in\Sigma_f^\mathrm{rat}.\]
Then $\s,\s'\subseteq\t$, but since $\tilde r\notin \t$ we have $\s\wedge\s'\nsubseteq\t$ that contradicts the minimality of $\s\wedge\s'$.
\endproof
\end{lem}

\begin{lem}\label{TwoChildIntegralRadiusLemma}
Let $\t\in\Sigma_f$ with at least two children in $\Sigma_f^\mathrm{rat}$. Then $d_\t=\rho_\t\in\Z$ and $\t\in\Sigma_f^\mathrm{rat}$. More precisely, if $\s,\s'\in\Sigma_f^\mathrm{rat}$ such that $\s\subsetneq\s\wedge\s'\supsetneq\s'$, then
\[\rho_{\s\wedge\s'}=v(w_\s-w_{\s'})=d_{\s\wedge\s'},\]
where $w_\s$ and $w_{\s'}$ are rational centres of $\s$ and $\s'$ respectively.
\proof
Clearly it suffices to prove the second statement as $v(w_\s-w_{\s'})\in\Z$. For our assumptions $\s'\not\subseteq\s$. Then by Lemma \ref{DifferentialsLemma} there exists $r\in\s'$ so that $v(r-w_\s)=\rho_{\s\wedge\s'}$. Thus, 
\[v(w_\s-w_{\s'})=\min\{v(w_\s-r),v(r-w_{\s'})\}=\rho_{\s\wedge\s'},\]
as $v(r-w_{\s'})\geq\rho_{\s'}>\rho_{\s\wedge\s'}$ by Lemma \ref{rationalCentresRadiusChildParentLemma}. 
Finally, $d_{\s\wedge\s'}=\rho_{\s\wedge\s'}$ follows from Lemma \ref{NoRationalClusterNoRationalChildLemma}.
\endproof
\end{lem}

\begin{defn}\label{EpsilonDefinition}
For a proper cluster $\s$ set
\[\epsilon_\s:=v(c_f)+\sum_{r\in\roots}\rho_{r\wedge\s}.\]
\end{defn}

\begin{exa}\label{Examplef=x^11-3x^6+9x^5-27}
Let $f=x^{11}-3x^6+9x^5-27\in\Q_3[x]$. The set of roots of $f$ is \[\roots=\{\sqrt[3]{3}, \zeta_3\sqrt[3]{3},\zeta_3^2\sqrt[3]{3},-\sqrt[3]{3}, -\zeta_3\sqrt[3]{3},-\zeta_3^2\sqrt[3]{3},\sqrt[5]{3}, \zeta_5\sqrt[5]{3}, \zeta_5^2\sqrt[5]{3}, \zeta_5^3\sqrt[5]{3}, \zeta_5^4\sqrt[5]{3}\},\]
where $\zeta_q$ is a primitive $q$-th root of unity for $q=3,5$. Then the proper clusters of $f$ are
\[\s_1=\{\sqrt[3]{3}, \zeta_3\sqrt[3]{3},\zeta_3^2\sqrt[3]{3}\},\quad\s_2=\{-\sqrt[3]{3}, -\zeta_3\sqrt[3]{3},-\zeta_3^2\sqrt[3]{3}\},\quad\s_3=\s_1\cup\s_2,\quad\roots\]
with $d_{\s_1}=d_{\s_2}=\frac{5}{6}$, $d_{\s_3}=\frac{1}{3}$ and $d_\roots=\frac{1}{5}$.
The graphic representation of the cluster picture of $f$ is then
   \begin{center}
    \includegraphics[trim=6cm 21.8cm 6cm 4.5cm,clip]{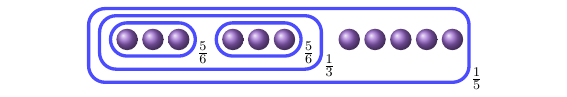}
    \end{center}
where the subscripts of clusters (represented as circles) are their depths.

Furthermore, note that $0$ is a rational centre for all proper clusters and we have $\rho_{\s_1}=\rho_{\s_2}=\rho_{\s_3}=\frac{1}{3}$ and $\rho_\roots=\frac{1}{5}$.

Finally, for every cluster $\s$ we can also compute $\nu_\s$ and $\epsilon_\s$, that are
\[\nu_{\s_1}=\nu_{\s_2}=\sfrac{9}{2},\quad\nu_{\s_3}=\epsilon_{\s_1}=\epsilon_{\s_2}=\epsilon_{\s_3}=3,\quad\nu_{\roots}=\epsilon_{\roots}=\sfrac{11}{5}.\]
\end{exa}
\begin{exa}
Let $f=x^9+12x^6+36x^3+81\in\Q_3[x]$ and fix an isomorphism $\ch \Q_3\simeq \C$. Then the set of roots of $f$ is \[\roots=\{\sqrt[3]{3^2}, \zeta_3\sqrt[3]{3^2},\zeta_3^2\sqrt[3]{3^2},\zeta_9\sqrt[3]{3}, \zeta_9^2\sqrt[3]{3},\zeta_9^4\sqrt[3]{3},\zeta_9^5\sqrt[3]{3},\zeta_9^7\sqrt[3]{3},\zeta_9^8\sqrt[3]{3}\},\]
where $\zeta_q=e^{2\pi i/q}$ is a primitive $q$-th root of unity for $q=3,9$. Then the proper clusters of $f$ are
\[\arraycolsep=1.4pt\def\arraystretch{1.4}\begin{array}{c}
\s_1=\{\sqrt[3]{3^2}, \zeta_3\sqrt[3]{3^2},\zeta_3^2\sqrt[3]{3^2}\},\quad\s_2=\{\zeta_9\sqrt[3]{3}, \zeta_9^4\sqrt[3]{3}, \zeta_9^7\sqrt[3]{3}\},\\ \s_3=\{\zeta_9^2\sqrt[3]{3}, \zeta_9^5\sqrt[3]{3}, \zeta_9^8\sqrt[3]{3}\},\quad\s_4=\s_2\cup\s_3,\quad\roots
\end{array}\]
with $d_{\s_1}=\frac{7}{6}$, $d_{\s_2}=d_{\s_3}=\frac{5}{6}$, $d_{\s_4}=\frac{1}{2}$, and $d_\roots=\frac{1}{3}$.
The cluster picture of $f$ is then
   \begin{center}
    \includegraphics[trim=6cm 21.8cm 6cm 4.5cm,clip]{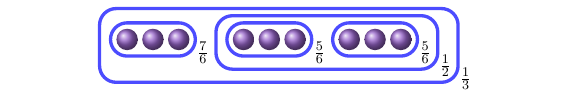}
    \end{center}
It is easy to see that $0$ is a rational centre for all proper clusters and that $\rho_{\s_1}=\frac{2}{3}$, $\rho_{\s_2}=\rho_{\s_3}=\rho_{\s_4}=\rho_\roots=\frac{1}{3}$. Finally,
\[\smaller{\nu_{\s_1}=\frac{11}{2},\quad\nu_{\s_2}=\nu_{\s_3}=5,\quad\nu_{\s_4}=4,\quad\nu_{\roots}=3;\qquad\epsilon_{\s_1}=4,\quad\epsilon_{\s_2}=\epsilon_{\s_3}=\epsilon_{\s_4}=\epsilon_\roots=3.}\]
\end{exa}

The goal of this section is to describe the \Np-regularity of $f\in K[x]$ in terms of conditions on its cluster picture.


\begin{nt}
If $p>0$, we denote by $|\cdot|_p$ the standard $p$-adic absolute value attached to $\Q$, i.e.\ $|a|_p=p^{-v_p(a)}$ for all $a\in\Q$. If $p=0$, then we write $|\cdot|_p$ for the function on $\Q$ identically equal to $1$, i.e.\ $|a|_p=1$ for all $a\in\Q$.
\end{nt}

\begin{lem}\label{Clusters-NewtonPolygonRegularity(MainLemma)}
Suppose that $x\nmid f$ and that $\NP f$ is a segment $L$ of slope $-\rho$. Let $n$ be the denominator of $\rho$. Then $f$ is \Np-regular if and only if all proper clusters $\s\in\Sigma_f$ with $|\s|>|\rho|_p$
satisfy $d_\s=\rho$.

More precisely:
\begin{enumerate}[label=(\roman*)]
    \item If $\s\in\mathring\Sigma_f$ with $|\s|>|\rho|_p$ but $d_\s>\rho$, then $\ch {f|_L}$ has a non-zero multiple root $\bar u=\frac{r^n}{\pi^{n\rho}} \mod \pi$, for some (any) $r\in\s$. \label{ClustersNewtonRegularityLemmaParti}
    \item The multiplicity of a root $\bar u\in\bar k^\times$ of $\ch {f|_L}$ equals $|\s^0|/n$, where \label{ClustersNewtonRegularityLemmaPartii} \[\s^0=\left\{r\in\roots\mid\bar u=\tfrac{r^n}{\pi^{n\rho}}\mod \pi\right\}.\] 
    \item All multiple roots of $\ch {f|_L}$ come from clusters $\s$ as described in \ref{ClustersNewtonRegularityLemmaParti}. \label{ClustersNewtonRegularityLemmaPartiii}
\end{enumerate}
\proof
Let $q$ be the highest power of $p$ dividing $n$ (set $q=1$ if $p=0$). Let $m=n/q$ so that $p\nmid m$. Let $\roots=\{r_i\mid i=1,\dots,D\}$ be the (multi-)set of roots of $f$, where $D:=\deg f$. Fix some choice of $\sqrt[n]{\pi}$ and define $\bar u_i\in\bar k^\times$ as $\bar u_i=r_i/\pi^\rho\mod\pi$, for all $i=1,\dots, D$. Firstly, note that there exists a proper cluster $\s$ with $|\s|>|\rho|_p$ and $d_\s>\rho$ if and only if there exists a subset $I\subseteq\{1,\dots,D\}$ of size $|I|>q$ such that $\bar u_{i_1}=\bar u_{i_2}$ for all $i_1,i_2\in I$.
Indeed, given $\s$, then $I=\{i\in\{1,\dots,D\}\mid r_i\in \s\}$, while given $I$, then $\s=\{r_i\mid \bar u_i=\bar u_{i_0},\text{ for any }i_0\in I\}$.
Secondly, recall that $f$ is not \Np-regular if and only if $\ch {f|_L}$ has a multiple root in $\bar  k^\times$. Therefore we will prove that $\ch {f|_L}$ has a non-zero multiple root if and only if there exists a subset $I\subseteq\{1,\dots,D\}$ with size $|I|>q$ and such that $\bar u_{i_1}=\bar u_{i_2}$ for all $i_1,i_2\in I$.



Note that for the lower convexity of $\NP{f}=L$, we have
\[\ch{f|_L}(x^n)=\pi^{-(v(c_f)+D\rho)}f(\pi^\rho x)\mod\pi.\]
Hence $\{\bar u_i \mid i=1,\dots,D\}$ is the multiset of roots of $\ch {f|_L}(x^n)$. Then there exists an $n$-to-$1$ map
\[\bar\phi:\mapcenter \{\bar u_i\}--{\{\bar w_j\}},\bar u_i--\bar u_i^m.,\]
where $\{\bar w_j\mid j=1,\dots,D/n\}$ is the multiset of roots of $\ch {f|_L}$. Note that $\bar w_j\neq 0$ for all $j=1,\dots,D/n$, so all roots of $\ch {f|_L}$ are non-zero.

Now, suppose that $f$ is not \Np-regular. We want to show that there exists a subset $I\subset\{1,\dots,D\}$ with $|I|>q$ such that $\bar u_{i_1}=\bar u_{i_2}$ for all $i_1,i_2\in I$. Since $f$ is not \Np-regular, its reduction $\ch {f|_L}$ has a (non-zero) multiple root. Then there exist $j_1,j_2\in\{1,\dots,D/n\}$ so that $\bar w_{j_1}=\bar w_{j_2}=:\bar w$. Hence, by the definition of $\bar\phi$, for some (any) $\bar u\in\bar\phi^{-1}(\bar w)$, there are at least $2 q$ $\bar u_i$'s with $\bar u_i=\bar u$. Let $I$ denote the set of their indices. Then $|I|\geq 2q>q$ and $\bar u_{i_1}=\bar u_{i_2}$ for all $i_1,i_2\in I$, as required.

On the other hand, suppose that there exists a subset $I\subset\{1,\dots,D\}$ with $|I|>q$ and such that $\bar u_{i_1}=\bar u_{i_2}$ for all $i_1,i_2\in I$. We want to show that $\ch {f|_L}$ has a multiple root, that is there exist two indices $j_1,j_2\in \{1,\dots,D/n\}$ such that $\bar w_{j_1}=\bar w_{j_2}$. Suppose not and let $j\in\{1,\dots,D/n\}$ such that $\bar w_j=\bar u_i^m=\bar\phi(\bar u_i)$ for some (all) $i\in I$. 
Then the polynomial $x^n-\bar w_j=(x^m-\bar w_j)^{q}\in\bar  k[x]$, factor of $\ch{f|_L}(x^n)$, should have a root of order $|I|>q$. This would imply $x^m-\bar w_j$ is inseparable, a contradiction as $p\nmid m$.

The parts \ref{ClustersNewtonRegularityLemmaParti}, \ref{ClustersNewtonRegularityLemmaPartii} and \ref{ClustersNewtonRegularityLemmaPartiii} of the lemma follow from above:

\noindent \ref{ClustersNewtonRegularityLemmaParti} Given a proper cluster $\s\in\Sigma_f$ with $|\s|>|\rho|_p$ and $d_\s>\rho$, we showed that $\ch{f|_L}$ has a non-zero multiple root $\bar w_j=\bar u_i^n=\nicefrac{r_i^n}{\pi^{n\rho}}\mod\pi$, where $r_i$ is any root in $\s$.

\noindent \ref{ClustersNewtonRegularityLemmaPartii} By the definition of $\bar \phi$, given $\bar w\in\bar k$, the number of $\bar w_j$'s such that $\bar w_j=\bar w$ equals $|\s^0|/n$, where $\s^0=\{r_i\mid \bar u_i^n=\bar w\}$.

\noindent \ref{ClustersNewtonRegularityLemmaPartiii} Given a (non-zero) multiple root $\bar w$ of $\ch{f|_L}$ we showed that there exists $I\subseteq\{1,\dots,D\}$, with $|I|>q$ and $\bar u_{i_1}=\bar u_{i_2}$ for any $i_1,i_2\in I$, such that $\bar u_i^n=\bar w$ for all $i\in I$. The set $\s=\{r_i\mid \bar u_i=\bar u_{i_0},\text{ for any }i_0\in I\}$ is a proper cluster as in \ref{ClustersNewtonRegularityLemmaParti}. 
\endproof
\end{lem}

\begin{thm}\label{RegularityClustersTheorem}
Let $w\in K$ and $f_w(x)=f(x+w)$. For all clusters $\s\in\Sigma_f$ define $\lambda_\s=\min_{r\in\s}v(r-w)$, and let $b$ be the denominator of $\lambda_\s$. Then $f_w$ is \Np-regular if and only if all proper clusters $\s\in\Sigma_f$ with $|\s|>|\lambda_\s|_p$ have $d_\s=\lambda_\s$.

More precisely:
\begin{enumerate}[label=(\roman*)]
    \item Let $\s\in\mathring\Sigma_f$ with $|\s|>|\lambda_\s|_p$ but $d_\s>\lambda_\s$, and let $r\in\s$ with $v(r-w)=\lambda_\s$. Then $\ch{f_w|_L}$ has a non-zero multiple root $\bar u=\frac{(r-w)^b}{\pi^{b\lambda_\s}}\mod \pi$, where $L$ is the edge of $\NP{f_w}$ of slope $-\lambda_\s$.\label{RegularityClustersThmParti} 
    \item Let $L$ be an edge of $\Np(f_w)$ of slope $-\lambda$. Let $l$ be the denominator of $\lambda$. The multiplicity of a root $\bar u\in\bar k^\times$ of $\ch{f_w|_L}$ equals $|\s^0|/l$, where \label{RegularityClustersThmPartii}
    \[\s^0=\big\{r\in\roots\mid v(r-w)=\lambda\quad\text{and}\quad\bar u=\tfrac{(r-w)^l}{\pi^{l\lambda}}\mod \pi\big\}.\] 
    \item For every edge $L$ of $\NP{f_w}$, the multiple roots of $\ch{f_w|_L}$ come from proper clusters $\s$ for $f$ as described in \ref{RegularityClustersThmParti}. \label{RegularityClustersThmPartiii}
\end{enumerate}
\proof
Let $\roots_w$ be the set of roots of $f_w$. Note that we have a natural bijection $\roots\rightarrow\roots_w$, $r\mapsto r-w$, which induces a bijective function $\psi:\Sigma_{f}\rightarrow\Sigma_{f_w}$, 
sending \[\s=\roots\cap\{x\in\bar K\mid v(x-z)>d\}\quad\mapsto\quad\psi(\s)=\roots_w\cap\{x\in\bar K\mid v(x+w-z)>d\}.\]
In particular, if $\s\in\Sigma_{f}$, $|\s|=|\psi(\s)|$, $d_\s=d_{\psi(\s)}$ and
\[\lambda_\s=\min_{r\in\s}v(r-w)=\min_{r\in\psi(\s)}v(r).\]
Hence it suffices to show the theorem for $w=0$.

Assume $w=0$. Let $f=c_f\cdot g_0\cdot g_1\dots g_t$ be a factorisation of Theorem \ref{NewtonpolygonFactorisation(Theorem)}. Note that if $t=0$, then either $f\in K$ or $f\in Kx$. In both cases, $f$ is clearly \Np-regular and has no proper clusters. Then assume $t>0$ and let $-\rho_i$ be the slope of $\NP{g_i}$ for any $i=1,\dots,t$. Denote by $\roots$ the set of roots of $f$ and by $\roots_i$ the set of roots of $g_i$ for $i=0,\dots,t$. Note that the $\roots_i$'s are pairwise disjoint. From Remark \ref{Remarkfor(Lemma1)}, for every edge $L$ of $\NP f$ there exists $i$ such that $\ch{f|_L}=\bar c_i\cdot\ch{g_i|_{\NP{g_i}}}$ for some $\bar c_i\in k^\times$. Hence, by Lemma \ref{regularityfactorisation(Lemma)} and Lemma \ref{Clusters-NewtonPolygonRegularity(MainLemma)}, we need to prove that there exists a proper cluster $\s\in\Sigma_f$ such that $|\s|>|\lambda_\s|_p$ and $d_\s>\lambda_\s$ if and only if for some $i=1,\dots,t$ there exists a proper cluster $\s_i\in\Sigma_{g_i}$ such that $|\s_i|>|\lambda_{\s_i}|_p=|\rho_i|_p$ and $d_{\s_i}>\lambda_{\s_i}=\rho_i$. We will show that one can choose $\s=\s_i$.

First, note that if $\s$ is a proper cluster , then $\s\not\subseteq\roots_0$, as $|\roots_0|\leq 1$. Furthermore, if $\s\in\Sigma_f$ contains roots of different valuations, that is $\s\nsubseteq\roots_i$ for all $i$, then 
\[d_\s=\min_{r,r'\in\s}v(r-r')=\min_{r\in\s}v(r)=\lambda_\s=\min\{\rho_i\mid \roots_i\cap\s\neq\varnothing\}.\] 

Now suppose there exists a proper cluster $\s\in\Sigma_f$ such that $|\s|>|\lambda_\s|_p$ and $d_\s>\lambda_\s$. For the observation above, the inequality $d_\s>\lambda_\s$ implies that $\s\subseteq\roots_i$ for some $i=1,\dots,t$. Let $\mathcal{D}$ be the $v$-adic disc such that $\s=\mathcal{D}\cap \roots$. Since $\s\subseteq\roots_i$, one has $\s=\mathcal{D}\cap \roots_i$ which means that $\s\in\Sigma_{g_i}$, as required.

Finally suppose that for some $i=1,\dots,s$, there exists a proper cluster $\s_i\in\Sigma_{g_i}$ such that $|\s_i|>|\rho_i|_p$ and $d_{\s_i}>\rho_i$. Let $r_i\in\s_i$. Then \[\s_i=\{x\in\bar K \mid v(x-r_i)\geq d_{\s_i}\}\cap\roots_i.\]
Consider the cluster $\s:=\{x\in\bar K \mid v(x-r_i)\geq d_{\s_i}\}\cap\roots$ of $f$. Clearly $\s_i\subseteq\s$. Therefore
\[\lambda_{\s_i}=\min_{r\in\s_i}v(r)\geq\min_{r\in\s}v(r)=\lambda_{\s},\]
which implies
\[d_{\s}=d_{\s_i}>\rho_i=\lambda_{\s_i}\geq\lambda_{\s},\]
where $d_{\s}=d_{\s_i}$ by construction. Again from the observation above the inequality $d_{\s}>\lambda_{\s}$ implies that $\s$ is contained in $\roots_j$ for some $j$. As $\s\cap\roots_i\supseteq\s_i\cap\roots_i=\s_i$, we must have $\s\subseteq\roots_i$. Thus $\s=\s_i$, that concludes the proof.
\endproof
\end{thm}

\begin{cor}\label{RegularityAfterTranslationCorollary}
Let $f\in K[x]$ be a separable polynomial. Let $w\in K$ and $f_w(x)=f(x+w)$. Then $f_w$ is \Np-regular if and only if all proper clusters $\s\in\Sigma_f$ have rational centre $w$ and those with $|\s|>|\rho_\s|_p$ satisfy $d_\s=\rho_\s$.
\proof
If $f_w$ is \Np-regular, then, from the previous theorem, all proper clusters $\s\in\Sigma_f$ with $|\s|>|\lambda_\s|_p$ have $d_\s=\lambda_\s$, where $\lambda_\s=\min_{r\in\s}v(r-w)$. First let $\s\in\Sigma_f$ proper and assume $|\s|>|\lambda_\s|_p$. Then
\[d_\s=\lambda_\s=\min_{r\in\s}v(r-w)\leq \max_{z\in K}\min_{r\in\s}v(r-z)=\rho_\s\leq d_\s,\]
so $d_\s=\lambda_\s=\rho_\s$, and $w$ is a rational centre of $\s$. Now assume $|\s|\leq|\lambda_\s|_p$. In particular, $p>0$ and $\lambda_\s\notin\Z$, and so \[\min_{r\in\s}v(r-w)=\lambda_\s\neq v(w-w_\s),\]
where $w_\s$ is a rational centre of $\s$. Let $r\in\s$ such that $v(r-w)=\lambda_\s$.
Then
\[\rho_\s\leq v(r-w+w-w_\s)=\min\{\lambda_\s,v(w-w_\s)\}\leq\lambda_\s.\]
Clearly
\[\rho_\s=\max_{z\in K}\min_{r\in\s}v(r-z)\geq \min_{r\in\s}v(r-w)=\lambda_\s,\]
that implies $\rho_\s=\lambda_\s=\min_{r\in\s}v(r-w)$. Hence $w$ is a rational centre of $\s$.

On the other hand, suppose that all proper clusters $\s\in\Sigma_f$ have rational centre $w\in K$ and those with $|\s|>|\rho_\s|_p$ satisfy $d_\s=\rho_\s$. Then $\rho_\s=\min_{r\in\s}v(r-w)$ for any $\s\in\Sigma_f$. Thus $f_w$ is $\Np$-regular again by Theorem \ref{RegularityClustersTheorem}.
\endproof
\end{cor}

\begin{defn}\label{AlmostRationalDefinition}
We say that $f$ has an \textit{almost rational cluster picture} if all proper clusters $\s\in\Sigma_f$ with $|\s|>|\rho_\s|_p$ have $d_\s=\rho_\s$.
\end{defn}

In the following we give different characterisations of the previous definition.

\begin{cor}\label{AlmostRationalTameCaseCorollary}
Suppose that $K(\roots)/K$ is a tamely ramified extension. Then $f$ has an almost rational cluster picture if and only if every proper cluster $\s\in\Sigma_f$ is $G_K$-invariant.
\proof
Since $K(\roots)/K$ is tamely ramified, every cluster $\s\in\Sigma_f$ has $|\rho_\s|_p\leq 1$. Therefore the corollary follows from Remark \ref{TameNestedClusterPictureRemark}. 
\endproof
\end{cor}

\begin{cor}\label{RegularityAfterTranslationTamecaseCorollary}
Suppose that $K(\roots)/K$ is a tamely ramified extension. Then $f_w$ is \Np-regular for some $w\in K$ if and only if $\Sigma_f$ is nested.
\proof
First note that every cluster $\s\in\Sigma_f$ has $|\rho_\s|_p\leq 1$, as $K(\roots)/K$ is tamely ramified. Therefore from Corollary \ref{RegularityAfterTranslationCorollary}, we need to prove that $\Sigma_f$ is nested if and only if all clusters $\s\in\Sigma_f$ have $d_\s=\rho_\s$ and rational centre $w$, for some $w\in K$. But this follows from Remark \ref{TameNestedClusterPictureRemark}.
%
%
\endproof
\end{cor}

\begin{cor}\label{NewtonPolygonAlmostRationaCharacterisationCorollary}
The polynomial $f$ has an almost rational cluster picture if and only if for every $r\in\roots\setminus K$, there exists $w\in K$ so that $r_w^b:=\frac{(r-w)^b}{\pi^{b\cdot v(r-w)}}\mod\pi$ is a simple root of $f_w|_L$, where $b$ is the denominator of $v(r-w)$, $f_w(x)=f(x+w)$ and $L$ is the edge of $\NP{f_w}$ of slope $-v(r-w)$.
\proof
Fix $\tilde r\in\roots\setminus K$ and let $\s$ be the smallest proper cluster containing $\tilde r$. Let $w_\s$ be a rational centre of $\s$. Note that $v(\tilde r-w_\s)=\rho_\s=\min_{r\in\s} v(r-w_\s)$, for the choice of $\s$, as $\tilde r\notin K$. 
Moreover, for any proper cluster $\t$ containing $\tilde r$, we have $\s\subseteq\t$. In particular, $w_\s$ is a rational centre of all such clusters. Let $L$ be the edge of $\Np(f_{w_\s})$ of slope $-\rho_\s$. Theorem \ref{RegularityClustersTheorem} shows that $\tilde r_{w_\s}^{b_\s}$ is a multiple root of $f_{w_\s}|_L$ if and only if there exists $\t\in\Sigma_f$ such that $\tilde r\in\t$, $|\t|>|\rho_\t|_p$ and $d_\t>\rho_\t$. Therefore if $f$ has an almost rational cluster picture, then $\tilde r_{w_\s}^{b_\s}$ is a simple root.

Suppose there exists $\t\in\Sigma_f$ such that $|\t|>|\rho_\t|_p$ and $d_\t>\rho_\t$. Then $\t\cap K=\varnothing$. By Theorem \ref{RegularityClustersTheorem}, it remains to show that for any $w\in K$, we have $|\t|>|\lambda_\t|_p$ and $d_\t>\lambda_\t$, where $\lambda_\t=\min_{r\in\t}v(r-w)$. First note $d_\t>\rho_\t\geq\lambda_\t$. Moreover, in the proof of Corollary \ref{RegularityAfterTranslationCorollary}, we saw that $|\t|\leq|\lambda_\t|_p$ implies $\rho_\t=\lambda_\t$, which contradicts $|\t|>|\rho_\t|_p$.
\endproof
\end{cor}

\begin{lem}\label{UberevenLemma}
Suppose $f$ has an almost rational cluster picture. Let $\s\in\Sigma_f$ proper. If
$d_\s>\rho_\s$, then $p>0$ and $|\s|$ is a $p$-power.
In particular, if $w_\s$ is a rational centre of $\s$, for any $r\in\s$, the elements $r-w_\s$ are all the roots of a polynomial with coefficients in $K^\mathrm{s}$, and constant term $c$ such that $|v(c)|_p\geq 1$.
\proof
Let $\s\in\Sigma_f$ proper, with $d_\s>\rho_\s$.
Since $f$ has an almost rational cluster picture, we must have $|\s|\leq|\rho_{\s}|_p$. Since $\s$ is proper, $p>0$.
Let $b_\s$ be the denominator of $\rho_\s$. Then $v_p(b_\s)>1$.
Fix a rational centre $w_\s$ of $\s$ and a root $r\in\s$ such that $v(r-w_\s)=\rho_\s$. Consider $\s'=\{x\in\roots\mid v(x-r)>\rho_\s\}$. Then $\s\subseteq\s'\leq\s^\mathrm{rat}$ and $|\s'|\leq|\rho_\s|_p$ (as $d_{\s'}>\rho_\s=\rho_{\s'}$).
Let $I_w$ be the wild inertia subgroup of $G_K$. As $v(r-w_\s)=\rho_\s$ there exist $\sigma_1={id},\sigma_2,\dots,\sigma_{|\rho_\s|_p}\in I_w$ such that $\sigma_i(r)\neq\sigma_j(r)$ if $i\neq j$. Moreover, $v(\sigma_i(r)-r)>\rho_\s$ from the definition of $I_w$. Therefore $\sigma_i(r)\in\s'$ for all $i$ and so $|\rho_\s|_p\leq|\s'|$. Thus $|\s'|=|\rho_\s|_p$ and $\s\subseteq\s'=\{\sigma_i(r)\mid i=1,\dots,|\rho_\s|_p\}$. Finally, as $\s'$ contains only conjugates of $r\in\s$, the cluster $\s'$ is union of orbits of $\s$. In particular, $|\s|\mid|\s'|=|\rho_\s|_p$, and so $|\s|$ is a $p$-power. The rest of the lemma follows.
\endproof
\end{lem}

\begin{prop}\label{AlmostRationalDiskCharacterisationProposition}
The polynomial $f$ has an almost rational cluster picture if and only if for every proper cluster $\s\in\Sigma_f$ one of the following is satisfied:
\begin{enumerate}[label=(\alph*)]
    \item the smallest disc containing $\s$ also contains a rational point;\label{AlmostRationalCharacterisation(a)}
    \item $p>0$ and after a translation by an element of $K$, the elements in $\s$ are all the roots of a polynomial with coefficients in $K^\mathrm{s}$ of $p$-power degree and constant term $c$ such that $|v(c)|_p\geq 1$. \label{AlmostRationalCharacterisation(b)}
\end{enumerate}
\proof
First of all note that point \ref{AlmostRationalCharacterisation(a)} is equivalent to requiring $d_\s=\rho_\s$. Therefore by Lemma \ref{UberevenLemma} it only remains to show that if $d_\s>\rho_\s$ and \ref{AlmostRationalCharacterisation(b)} is satisfied, then $|\s|\leq |\rho_\s|_p$. Let $F\in K^\mathrm{s}[x]$ be the polynomial in \ref{AlmostRationalCharacterisation(b)} and let $w\in K$ such that $r-w$, for $r\in\s$, are all the roots of $F$. We have $\rho_s\geq \min_{r\in\s}v(r-w)$. Fix $r\in\s$ such that $\rho_\s\geq v(r-w)=:\rho$. Since $d_\s>\rho_\s\geq  v(r-w)$, we have $v(r'-w)=v(r-w)=\rho$ for any $r'\in\s$.  Then
\[|\s|=\deg F = |1/\deg F|_p\leq|v(c)/\deg F|_p=|\rho|_p.\]
Let $w_\s$ be a rational centre of $\s$. Suppose by contradiction that $\rho_\s>\rho$. Let $r_\s\in\s$ such that $v(r_\s-w_\s)=\rho_\s$. Hence
\[v(w-w_\s)=v(w-r_\s+r_\s-w_\s)=\min\{\rho,\rho_\s\}=\rho.\]
But then $\rho\in\Z$, which contradicts $|\s|\leq|\rho|_p$.
\endproof
\end{prop}

\begin{exa}\label{AlmostRationalExample}
Let $p$ be a prime number and let $a\in\Z_p$, $b\in\Z_p^\times$ such that the polynomial $x^2+ax+b$ is not a square modulo $p$. Let $f\in\Q_p[x]$ given by $f(x)=(x^6+ap^4x^3+bp^8)((x-p)^3-p^{11})$. For any prime $p$ the rational cluster picture of $f$ is
\begin{center}
    \includegraphics[trim=6cm 22.2cm 6cm 4.5cm,clip]{Clusters/Cluster5.pdf}
\end{center}
where $\rho_{\t_3}=\frac{4}{3}$, $\rho_{\t_4}=\frac{11}{3}$, and $\rho_\roots=1$. 

If $p\neq 3$, then the proper clusters of $\Sigma_f$ coincide with the rational clusters above and $d_\s=\rho_\s$ for any $\s=\t_3,\t_4,\roots$. In particular, $f$ has an almost rational cluster picture when $p\neq 3$. 

Suppose $p=3$. Then the cluster picture of $f$ is 
\begin{center}
    \includegraphics[trim=6cm 22cm 6cm 4.5cm,clip]{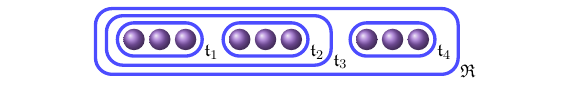}
\end{center}
where $d_{\t_1}=d_{\t_2}=\frac{11}{6}$, $d_{\t_3}=\rho_{\t_1}=\rho_{\t_2}=\frac{4}{3}$, $d_{\t_4}=\frac{25}{6}$ and $d_\roots=1$. Thus $f$ has an almost rational cluster picture for all $p$.
\end{exa}

We conclude this section by showing that the \textit{cluster picture centred at $0$} completely determines the Newton polygon of $f$.

\begin{defn}\label{CentredClusterDefinition}
Let $z\in \bar K$. A \textit{cluster centred at $z$} is a cluster cut out by a $v$-adic disc of the form $\mathcal{D}=\{x\in\bar K\mid v(x-z)\geq d\}$ for some $d\in\Q$.
\end{defn}

\begin{defn}\label{CentredClusterPictureDefinition}
Let $z\in \bar K$. Define $\Sigma_f^z$ to be the set of all clusters centred at $z$. Write $\mathring\Sigma_f^z$ for the set $\Sigma_f^z\smallsetminus\{\{z\}\}$.
Note that $\Sigma_f^z$ is \textit{nested}, i.e.\ every cluster $\s\in\Sigma_f^z$ has at most one child in $\Sigma_f^z$.
\end{defn}

\begin{defn}\label{RadiusWithRespectToCentreDefinition}
Let $z\in \bar K$, and let $\s\in\Sigma_f\setminus\{\{z\}\}$. The \textit{radius of $\s$ with respect to the centre $z$} is
\[\rho_\s^z=\min_{r\in\s}v(r-z).\]
The \textit{cluster picture centred at $z$} of $f$ is the collection of all clusters in $\Sigma_f^z$ together with their radii with respect to $z$.
Finally set
\[\epsilon_\s^z:=v(c_f)+\sum_{r\in\roots}\rho_{r\wedge\s}^z.\]
\end{defn}

\begin{rem}
From the definitions above, if $\s$ is a cluster centred at $z\in K^\mathrm{s}$, then $\s=\roots\cap\{x\in\bar K\mid v(x-z)\geq \rho_\s^z\}$. But this does not mean $z$ is a centre for $\s$, that is false in general. For example, $\roots$ is clearly a cluster centred at any $z\in K^\mathrm{s}$, but there are elements of $K^\mathrm{s}$ which are not centres of $\roots$, e.g.\ any $z\in K^\mathrm{s}$ with valuation $v(z)<\min_{r\in\roots}v(r)$.
\end{rem}

\begin{rem}\label{RationalRadiusCentredRadiusRemark}
Let $\s\in\Sigma_f$ be a proper cluster with centre $z$ and rational centre $w$. Then $\s\in\Sigma_f^{z}$, $d_\s=\rho^z_\s$, $\nu_\s=\epsilon^z_\s$, $\rho_\s=\rho_\s^w$, and $\epsilon_\s=\epsilon_\s^w$.
Furthermore, $\s\in\Sigma_f^\mathrm{rat}$ if and only if $\s\in\Sigma_f^{w}$. 
\end{rem}

\begin{lem}\label{DescriptionNewtonPolygonClustersLemma}
Let $w\in K$ and let $f_w(x)=f(x+w)$. Then there is a $1$-to-$1$ correspondence between the clusters in $\mathring\Sigma_f^w$ and the edges of $\NP {f_w}$. More explicitly, let $\s_1\subset\dots\subset \s_n=\roots$ be the clusters in $\mathring\Sigma_f^w$ and let $\s_0=\{w\}$ if $\{w\}\in\Sigma_f^w$ or $\s_0=\varnothing$ otherwise. Then $\NP {f_w}$ has vertices $Q_i$, $i=0,\dots,n$, where
\begin{itemize}
    \item $Q_n=(|\roots|,\epsilon_\roots^w-|\roots|\rho_\roots^w)=(\deg f,v(c_f))$,
    \item $Q_i=(|\s_i|,\epsilon_{\s_i}^w-|\s_i|\rho_{\s_i}^w)=(|\s_i|,\epsilon_{\s_{i+1}}^w-|\s_i|\rho_{\s_{i+1}}^w)$,\quad for $i=1,\dots,n-1$,
    \item $Q_0=(|\s_0|,\epsilon_{\s_{1}}^w-|\s_0|\rho_{\s_{1}}^w)$.
\end{itemize}
and edges $L_i$, $i=1,\dots,n$, of slope $-\rho_{\s_i}^w$ linking $Q_{i-1}$ and $Q_i$.

Furthermore, for any $i=1,\dots,n$ we have 
\[\overline{f_w|_{L_i}}(x^{b_i})= \tfrac{u}{\pi^{v(u)}}\tprod_{r\in\s_i\setminus \s_{i-1}} (x+\tfrac{w-r}{\pi^{\rho_i}})\mod \pi,\qquad u=c_f\tprod_{r\in\roots\setminus\s} (w-r),\]
where $\rho_i=\rho_{\s_i}^w$, and $b_i$ is the denominator of $\rho_i$.
\proof
Without loss of generality we can assume $w=0$ so that $f_w=f$. First note that the coordinates of $Q_n$ are trivial. Now consider a factorisation $f=c_f\cdot g_0\cdot g_1\cdots g_s$ of Theorem \ref{NewtonpolygonFactorisation(Theorem)}. 
Recall the polynomials $g_j$ are monic and $g_0\mid x$.
Let $\roots_j$ be the set of roots of $g_j$. It follows from the definition of cluster centred at $0$ that
\[n=s,\quad\mbox{and}\quad \s_i=\bigcup_{j=0}^i\roots_j\quad\text{for all }i=0,\dots,n.\]
Therefore $\s_0=\roots_0$ and $\roots_i=\s_i\setminus\s_{i-1}$ for any $i=1,\dots,n$.

Let $i=1,\dots,n-1$. Then the $x$-coordinate of $Q_i$ follows as \[|\s_i|=\sum_{j=0}^i|\roots_j|=\sum_{j=0}^i\deg g_j=\deg\prod_{j=0}^ig_j.\] The $y$-coordinate of $Q_i$ equals the sum of $v(c_f)$ and the valuation of the constant term of $\prod_{j=i+1}^n g_j$, so
\[Q_i=\bigg(|\s_i|,v(c_f)+\sum_{j=i+1}^n|\roots_j|v(r_j)\bigg),\]
where $r_j$ is any root in $\roots_j$. But since $\s_i=\bigcup_{j=0}^i\roots_j$, we have $v(r_j)=\rho_{\s_j}^0$. Therefore
\[v(c_f)+\sum_{j=i+1}^n|\roots_j|v(r_j)=v(c_f)+\sum_{j=i+1}^n(|\s_j|-|\s_{j-1}|)\rho_{\s_j}^0=\epsilon_{\s_i}^0-|\s_i|\rho_{\s_i}^0.\]
Moreover, \[\epsilon_{\s_i}^0-|\s_i|\rho_{\s_i}^0=\epsilon_{\s_{i+1}}^0-|\s_i|\rho_{\s_{i+1}}^0\] from the easy computation 
$\epsilon_{\s_i}^0-\epsilon_{\s_{i+1}}^0=|\s_i|\big(\rho_{\s_i}^0-\rho_{\s_{i+1}}^0\big).$
Finally the $x$-coordinate of $Q_0$ is trivial, while its $y$-coordinate equals
\[v(c_f)+\sum_{j=1}^n|\roots_j|v(r_j)=v(c_f)+\sum_{j=1}^n(|\s_j|-|\s_{j-1}|)\rho_{\s_j}^0=\epsilon_{\s_1}^0-|\s_0|\rho_{\s_1}^0,\]
that concludes the first part of the proof as $|\s_0|=|\roots_0|=\deg g_0$.

The computation of $f|_{L_i}$ follows from Remark \ref{Remarkfor(Lemma1)}. Indeed, let $i=1,\dots,n$, and define $\bar c_i=u/\pi^{v(u)}\mod\pi$, where $u=c_f\tprod_{j=i+1}^ng_j(0)$. Then
$\ch{f|_{L_i}}(x^{b_i})=\bar c_i\cdot \ch{g_i|_{\NP{g_i}}}(x^{b_i})$, where $b_i$ is the denominator of $\rho_{\s_i}^0$. But \[\ch{g_i|_{\NP{g_i}}}(x^{b_i})=g_i\big(\pi^{\rho_{\s_i}^0}x\big)/\pi^{\rho_{\s_i}^0\deg g_i}\mod\pi.\] Thus the lemma follows as $\roots_i=\s_i\setminus\s_{i-1}$.
\endproof
\end{lem}

\begin{nt}
Let $\s\in\mathring\Sigma_f^w$. Following the notation of Lemma \ref{DescriptionNewtonPolygonClustersLemma}, let $i\in\{1,\dots,n\}$ be such that $\s=\s_i$. We will write $L_\s^w$ for the edge $L_i$.
\end{nt}

\section{Description of a regular model}\label{DescriptionModelsSection}
From now on, assume $\mathrm{char}(K)\neq 2$ and let $C/K$ be a hyperelliptic curve, i.e.\ a geometrically connected, smooth, projective curve, equipped with a separable morphism $C\rightarrow \P^1_K$ of degree $2$. Let $y^2=f(x)$ be a Weierstrass equation of $C$. Suppose $\deg f>1$. Let $g$ be the genus of $C$.
Accordingly with \cite{D2M2} we define the \textit{cluster picture} of $C$ as the cluster picture of $f$. Analogously, all definitions and notations attached to $f$ given in \S \ref{ClustersSection} (e.g. $\Sigma_f$, $\Sigma_f^\mathrm{rat}$, $\Sigma_f^z$) are given for $C$ in the same way (e.g. $\Sigma_C$, $\Sigma_C^\mathrm{rat}$, $\Sigma_C^z$). 
In particular, we will say that $C$ has an almost rational cluster picture if $f$ does (Definition \ref{AlmostRationalDefinition}).

For the following sections we will use the main definitions, notations and results of 
\cite[\S $3$]{Dok}. In particular, we recall (without stating) the definitions of Newton polytopes $\Delta$ and $\Delta_v$  attached to a polynomial $g\in K[x,y]$, $v$-vertices/edges/faces of $\Delta$, the denominator $\delta_\lambda$ of a $v$-face/edge $\lambda$, the slopes $s_1^\lambda, s_2^\lambda$ of a $v$-edge $\lambda$. 

\begin{nt}
We denote by $\Delta_v^w$ and $\Delta^w$ respectively the polytopes $\Delta_v$ and $\Delta$ attached to the polynomial $g_w(x,y)=y^2-f(x+w)$. The piecewise affine function $v:\Delta^w\rightarrow\R$ determining the bijection $\Delta^w\rightarrow\Delta_v^w$, $P\mapsto(P,v(P))$, will be denoted by $v$ (with a little abuse of notation). For a $v$-face $F$ of $\Delta^w$, denote by $v_F:\Delta^w\rightarrow\R$ the linear function equal to $v$ on $F$. Since the projection $\Delta_v^w\rightarrow\Delta^w$ is a bijection, given a vertex/edge/face $\lambda$ of $\Delta_v^w$ we will denote by the same symbol $\lambda$ the corresponding $v$-vertex/edge/face of $\Delta^w$. Since they are mainly used for indexing, this will not cause confusion.
\end{nt}

\begin{nt}
Given a $v$-edge $\lambda$ of $\Delta^w$, we will denote by $r_\lambda$ the smallest non-negative integer such that we can fix $\frac{n_i}{d_i}\in\Q$, for $i=0,\dots,r_\lambda+1$ so that
    \[s_1^\lambda=\frac{n_0}{d_0}>\frac{n_1}{d_1}>\ldots>\frac{n_{r_\lambda}}{d_{r_\lambda}}>\frac{n_{r_\lambda+1}}{d_{r_\lambda+1}}=s_2^\lambda,\quad\text{with\scalebox{0.9}{$\quad\begin{vmatrix}n_i\!\!\!&n_{i+1}\cr d_i\!\!\!&d_{i+1}\cr\end{vmatrix}=1$}}.\]
\end{nt}

Thanks to Lemma \ref{DescriptionNewtonPolygonClustersLemma} we can explicitly relate the Newton polytope $\Delta_v^w$ of $g_w(x,y)$ and the cluster picture centred at $w$ of $C$.

\begin{lem}\label{DescriptionNewtonPolytopeClustersLemma}
Let $w\in K$. Then there is a $1$-to-$1$ correspondence between the clusters in $\mathring\Sigma_C^w$ and the faces of the Newton polytope $\Delta_v^w$.
More explicitly, let $\s_1\subset\dots\subset \s_n=\roots$ be the clusters in $\mathring\Sigma_C^w$ and let $\s_0=\{w\}$ if $\{w\}\in\Sigma_C^w$ or $\s_0=\varnothing$ otherwise. Then $\Delta_v^w$ has vertices $T$, $Q_i$, $i=0,\dots,n$, where
\begin{itemize}
    \item $T=(0,2,0)$,
    \item $Q_n=(|\roots|,0,v(c_f))$,
    \item $Q_i=(|\s_i|,0,\epsilon_{\s_{i+1}}^w-|\s_i|\rho_{\s_{i+1}}^w)$ for $i=0,\dots,n-1$,
\end{itemize}
and edges $L_i$ ($i=1,\dots,n$), linking $Q_{i-1}$ and $Q_i$, and $V_j$ ($j=0,\dots,n$), linking $Q_j$ and $T$. 
Furthermore, (possible choices for) the slopes of the $v$-edges of $\Delta^w$ are:
\begin{itemize}
    \item \[s_1^{V_n}=\delta_{V_n}\tfrac{-v(c_f)+(|\roots|-2g)\rho_\roots^w}{2}\quad \mbox{and}\quad s_2^{V_n}=\lfloor s_1^{V_n}-1\rfloor;\]
    \item \[
    \renewcommand{\arraystretch}{2}
    \begin{array}{l}
        s_1^{V_i}=\delta_{V_i}\lb-\frac{\epsilon_{\s_{i}}^w}{2}+\lb\left\lfloor\frac{|\s_i|}{2}\right\rfloor+1\rb\rho_{\s_{i}}^w\rb,\\ s_2^{V_i}=\delta_{V_i}\lb-\frac{\epsilon_{\s_{i+1}}^w}{2}+\lb\left\lfloor\frac{|\s_i|}{2}\right\rfloor+1\rb\rho_{\s_{i+1}}^w\rb
    \end{array}\renewcommand{\arraystretch}{1}
\quad\text{for all }i=1,\dots,n-1;\]
    \item \[s_1^{V_0}=\delta_{V_0}\lb\tfrac{\epsilon_{\s_1}^w}{2}-\rho_{\s_1}^w\rb\quad\text{and}\quad s_2^{V_0}=\lfloor s_1^{V_0}-1\rfloor;\]
    \item \[s_1^{L_i}=\delta_{L_i}\lb-\tfrac{\epsilon_{\s_{i}}^w}{2}+\lb\left\lfloor\tfrac{|\s_i|}{2}\right\rfloor+1\rb\rho_{\s_{i}}^w\rb \quad \text{and}\quad s_2^{L_i}=\lfloor s_1^{L_i}-1\rfloor,\] for all $i=1,\dots,n$. 
    In particular, as $\delta_{L_i}$ is the denominator of $\rho_{\s_{i}}^w$, \[r_{L_i}=\begin{cases}1&\mbox{if }\delta_{L_i}\epsilon_{\s_i}^w\mbox{ is odd,}\\0&\mbox{if }\delta_{L_i}\epsilon_{\s_i}^w\mbox{ is even.}\end{cases}\]
\end{itemize}

Finally, for suitable choices of basis of the lattices in \cite[3.4, 3.5]{Dok}, we have
\[\overline{g_w|_{L_i}}(x^{b_i})= -\tfrac{u}{\pi^{v(u)}}\tprod_{r\in\s_i\setminus \s_{i-1}} (x+\tfrac{w-r}{\pi^{\rho_i}})\mod \pi,\qquad u=c_f\tprod_{r\in\roots\setminus\s_i} (w-r),\]
for any $i=1,\dots,n$, where $\rho_i=\rho_{\s_i}^w$, and $b_i$ is the denominator of $\rho_i$;
\[\overline{g_w|_{V_j}}(y)=y^{|\bar V_j(\Z)_\Z|-1}- \tfrac{u}{\pi^{v(u)}}\mod \pi,\qquad u=c_f\tprod_{r\in\roots\setminus\s_j} (w-r),\]
for any $j=0,\dots,n$, where $|\bar V_j(\Z)_\Z|$ is the number of integer points $P$ on the $v$-edge $V_j$ with $v(P)\in\Z$, endpoints included.
\proof
The structure of $\Delta_v^w$ follows from Lemma \ref{DescriptionNewtonPolygonClustersLemma}. For the computation of the slopes, we only need to individuate, for all the $v$-edges, the two points $P_0$ and $P_1$ of \cite[Definition 3.12]{Dok}. It is easy to see that the followings are admissible choices.
\begin{itemize}
    \item For $V_i$ and $L_i$ ($i=1,\dots,n$), choose $P_0=(|\s_i|,0)$ and $P_1=\lb\left\lfloor\frac{|\s_i|-1}{2}\right\rfloor,1\rb$.
    \item For $V_0$, choose $P_0=(0,2)$ and $P_1=(1,1)$;
\end{itemize}
The second part of the lemma then follows from the first one.  The computations of the reductions also follows from Lemma \ref{DescriptionNewtonPolygonClustersLemma} by choosing the lattices $Q_{i-1}+(b_i,0)\Z$ for $g_w|_{L_i}$ and $Q_{i}+(-|\s_i|/a,2/a)\Z$ for $g_w|_{V_j}$, where $a=|\bar V_j(\Z)_\Z|-1$.
\endproof
\end{lem}

\begin{nt}\label{FacesEdgesClustersNotation}
Let $C$ be as above and let $w\in K$. For every cluster $\s\in\mathring\Sigma_C^w$ denote by $F_\s^w$ the $v$-face of the Newton polytope $\Delta^w$ of $g_w(x,y)=y^2-f(x+w)$ that corresponds to $\s$.

Following the notation of Lemma \ref{DescriptionNewtonPolytopeClustersLemma}, let $i\in\{1,\dots,n\}$ be such that $\s=\s_i$. We will write $L_\s^w$, $V_\s^w$, $V_0^w$ for the $v$-edges $L_i$, $V_i$, $V_0$, respectively.
\end{nt}

\begin{exa}\label{ExampleC:y^2=x^11-3x^6+9x^5-27}
Let $C$ be the hyperelliptic curve over $\Q_3$ given by the equation $y^2=f(x)$ where $f(x)=x^{11}-3x^6+9x^5-27$ is the polynomial of Example \ref{Examplef=x^11-3x^6+9x^5-27}.

Its cluster picture centred at $0$ is
\begin{center}
    \includegraphics[trim=6cm 22.2cm 6cm 4.5cm,clip]{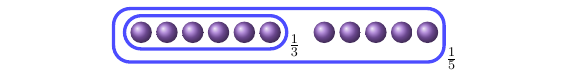}
    \end{center}
    where the subscripts represent the radii with respect to $0$. As we can see, $\Sigma_C^0$ consists of two clusters: $\s_1$ of size $6$, radius $\frac{1}{3}$ and $\epsilon_{\s_1}^0=3$, and $\s_2=\roots$ of size $11$, radius $\frac{1}{5}$ and $\epsilon_{\s_2}^0=\frac{11}{5}$. Therefore the picture of $\Delta^0$ broken into $v$-faces will be
\begin{center}
\begin{tikzpicture}
    \node (Q0) at ( 0, 0) {$Q_0$}; 
    \node (Q1) at ( 6, 0) {$Q_1$};
    \node (Q2) at ( 11,0) {$Q_2$};
    \node (T) at ( 0,2) {$T$};
       \draw (Q0) -- (Q1)node[midway, below] {$L_1$};
       \draw (Q1) -- (Q2)node[midway, below] {$L_2$}; 
       \draw (Q0) -- (T)node[midway, left] {$V_0$};
       \draw (Q1) -- (T)node[midway, below=7pt, left] {$V_1$};
       \draw (Q2) -- (T)node[midway, above=1pt] {$V_2$};
\end{tikzpicture}
\end{center}
where $T=(0,2)$, $Q_0=(0,0)$, $Q_1=(6,0)$, and $Q_2=(11,0)$. Denoting the values of $v$ on vertices, the picture becomes
\begin{center}
\begin{tikzpicture}
\draw[very thin,color=gray] (0.15,0.15) grid (10.85,1.85);
    \node (Q0) at ( 0, 0) {$3$}; 
    \node (Q1) at ( 6, 0) {$1$};
    \node (Q2) at ( 11,0) {$0$};
    \node (T) at ( 0,2) {$0$};
       \draw (Q0) -- (Q1);
       \draw (Q1) -- (Q2); 
       \draw (Q0) -- (T);
       \draw (Q1) -- (T);
       \draw (Q2) -- (T);
\end{tikzpicture}
\end{center}
\end{exa}

To state the theorems which describe the special fibre of the proper flat model $\mathcal{C}$ of $C$ we will construct in \S\ref{ConstructionModelsSection}, we need some definitions.

\begin{defn}\label{QuantitiesForTheoremsOnModelsDefinition}
Let $F/K$ be an unramified extension and let $\Sigma_F=\Sigma_{C_{F}}^\mathrm{rat}$ (i.e.\ set of clusters cut out by discs with centre in $F$). For any proper $\s\in\Sigma_F$ let $G_\s=\mathrm{Stab}_{G_K}(\s)$ and $K_\s=\lb K^\mathrm{s}\rb^{G_\s}$. We define the following quantities:

\begin{center}
\begin{tabular}{|l@{\quad}l@{$=\>\,$}l|}
\hline
\multicolumn{3}{|l|}{$\s\in\Sigma_F$, proper}\cr
\hline
radius       & $\rho_\s$         & $\max_{w\in F}\min_{r\in\s} v(r-w)$\cr
             & $b_\s$            & denominator of $\rho_\s$\cr
             & $\epsilon_\s$     & $v(c_f) + \sum_{r\in\roots} \rho_{r\wedge \s}$\cr
             & $D_\s$            & $1$ if $b_\s\epsilon_\s$ odd, $2$ if $b_\s\epsilon_\s$ even\cr
multiplicity & $m_\s$            & $(3-D_\s)b_\s$\cr
parity       & $p_\s$            & $1$ if $|\s|$ is odd, $2$ if $|\s|$ is even\cr
slope        & $s_\s$            & $\frac 12(|\s|\rho_\s+p_\s\rho_\s-\epsilon_\s)$\cr
             & $\gamma_\s$       & $2$ if $\s$ is even and $\epsilon_\s\!-\!|\s|\rho_\s$ is odd, $1$ otherwise\cr
             & $p_\s^0$          & $1$ if $\s$ is minimal and $\s\cap K_\s\neq\varnothing$, $2$ otherwise\cr
             & $s_\s^0$          & $-\epsilon_\s/2+\rho_\s$\cr
             & $\gamma_\s^0$     & 2 if $p_\s^0=2$ and $\epsilon_\s$ is odd, 1 otherwise\cr
\hline
\end{tabular}
\end{center}
\end{defn}

\begin{lem}\label{QuantitiesIndependentOfUnramifiedBaseChangeLemma}
Keep the notation of the previous definition and let $\s\in\Sigma_K$. Then $\s\in\Sigma_F$ but the quantities in Definition \ref{QuantitiesForTheoremsOnModelsDefinition} do not depend on $F$.
\proof
A cluster $\s\in\Sigma_F$ belongs to $\Sigma_K$ if and only if $\sigma(\s)=\s$ for any $\sigma\in G_K$. Then the result follows from Lemma \ref{RationalCentreTameExtensionLemma}.
\endproof
\end{lem}

\begin{rem}
Lemma \ref{DescriptionNewtonPolytopeClustersLemma} shed some light on the quantities we defined in Definition \ref{QuantitiesForTheoremsOnModelsDefinition}. Let $\s\in\Sigma_F$. Fix a rational centre $w_\s\in F$ of $\s$ such that $w_\s\in K_\s$ if $p_\s^0=1$. Denoting $V=V_\s^{w_\s}$, $L=L_\s^{w_\s}$, and $V_0=V_0^{w_\s}$, we have:
\begin{itemize}
    \item $b_\s=\delta_L$ and $r_L=2-D_\s$.
    \item $\gamma_\s=\delta_V$, $p_\s/\gamma_\s=\bar V(\Z)_\Z-1$ and $s_1^V=\gamma_\s s_\s$. If $V$ is internal, that is $\s\neq\roots$, then $s_2^V=\gamma_\s(s_\s-p_\s\tfrac{\rho_\s-\rho_{P(\s)}}{2})$.
    \item If $\s$ is minimal and so $V_0$ is an edge of $F_\s^{w_\s}$, then $\gamma_\s^0=\delta_{V_0}$, $p_\s^0/\gamma_\s^0=\bar V_0(\Z)_\Z-1$ and $s_1^{V_0}=-\gamma_\s^0 s_\s^0$.
\end{itemize}
\end{rem}

\begin{lem}\label{DsInteriorPointsLemma}
Let $\s\in\Sigma_C^\mathrm{rat}$ with rational centre $w\in K$. Then $D_\s=1$ if and only if $v_{F_\s^w}((a,1))\notin\Z$, for every $a\in\Z$.
\proof
If $D_\s=1$ then $r_{L_\s^w}=1$ by Lemma \ref{DescriptionNewtonPolytopeClustersLemma}, and so $v_{F_\s^w}((a,1))\notin\Z$, for every $a\in\Z$. Now let $c,d\in\Z$ such that $\rho_\s\cdot c+d=1/b_\s$. If $D_\s=2$, then $b_\s\epsilon_\s\in 2\Z$, so \[v_{F_\s^w}(cb_\s\epsilon_\s/2,1)=\frac{v_{F_\s^w}((cb_\s\epsilon_\s,0))}{2}=\frac{\epsilon_\s-(cb_\s\epsilon_\s)\rho_\s}{2}=\frac{db_\s\epsilon_\s}{2}\in \Z,\]
as required.
\endproof
\end{lem}

\begin{defn}\label{yRegularDefinition}
We say that $C$ is \textit{$y$-regular} if $p\nmid D_\s$ for every proper $\s\in\Sigma_C^\mathrm{rat}$, i.e.\ if either $p\neq 2$ or $D_\s=1$ for any proper $\s\in\Sigma_C^\mathrm{rat}$.
\end{defn}

\begin{rem}
Let $F/K$ be an unramified extension. Then from Lemma \ref{QuantitiesIndependentOfUnramifiedBaseChangeLemma}, if $C_F$ is $y$-regular then $C$ is $y$-regular.
\end{rem}

\begin{lem}\label{Regularity-RegularityLemma}
The hyperelliptic curve $C$ is $\Delta_v$-regular if and only if $C$ is $y$-regular and $f$ is \Np-regular.
\proof
Let $g(x,y)=y^2-f(x)$. If $C$ is $y$-regular and $f$ is \Np-regular, then $C$ is $\Delta_v$-regular by Lemma \ref{DescriptionNewtonPolytopeClustersLemma} and Lemma \ref{DsInteriorPointsLemma}. 

Conversely, if $C$ is $\Delta_v$-regular, then $f$ is \Np-regular, and all clusters have rational centre $0$ by Corollary \ref{RegularityAfterTranslationCorollary}. It remains to show that if $p=2$ then $D_\s=1$ for every proper $\s\in\Sigma_C^\mathrm{rat}$. Suppose there exists $\s\in\Sigma_C^\mathrm{rat}$ such that $D_\s=2$.
Consider the variety $\bar X_{F_\s^0}
$ (\cite[Definition 3.7]{Dok}). By Lemma \ref{DsInteriorPointsLemma}, the smoothness of $\bar X_{F_\s^0}
$ implies there exists $\s'\in\Sigma_C^\mathrm{rat}$, such that $|\s|-|\s'|=1$. Hence $\rho_\s\in\Z$ from Lemma \ref{RadiusTimesSizeIntegerLemma}. Therefore $\bar F_\s^0(\Z)=\bar F_\s^0(\Z)_\Z$, by Lemma \ref{DsInteriorPointsLemma}. 
But this gives a contradiction as it forces either $\ch{g|_{V_{\s'}^0}}$ or $\ch{g|_{V_\s^0}}$ to be a square.
\endproof
\end{lem}

\begin{defn}\label{GenussDefinition}
Let $\s\in\Sigma_F$ be a proper cluster and let $c\in\{0,\dots,b_\s-1\}$ such that $c\rho_\s-\frac{1}{b_\s}\in\Z$. Define
\[\tilde\s=\{\s'\in\Sigma_F\cup\{\varnothing\}\mid \s'<\s\text{ and }\tfrac{|\s'|}{b_\s}-c\epsilon_\s\notin 2\Z\},\]
where $\varnothing<\s$ if $\s$ is minimal and $p_\s^0=2$.

The \textit{genus $g(\s)$} of a rational cluster $\s\in\Sigma_F$ is defined as follows:
\begin{itemize}
    \item If $D_\s=1$, then $g(\s)=0$.
    \item If $D_\s=2$, then $2g(\s)+1$ or $2g(\s)+2$ equals
    \[\frac{|\s|-\sum_{\s'\in\Sigma_F,\s'<\s}|\s'|}{b_\s}+|\tilde \s|.\]
\end{itemize}
\end{defn}

\begin{defn}\label{SchemesXsDefinition}
Let $\Sigma_C^\mathrm{min}$ be the set of rationally minimal clusters of $C$ and let $\Sigma\subseteq\Sigma_C^\mathrm{min}$. For each cluster $\s\in\Sigma$, fix a rational centre $w_\s$; if possible, choose $w_\s\in\s$. Let $W$ be the set of these rational centres and define $\Sigma^W=\bigcup_{w\in W}\Sigma_C^w$. For any proper cluster $\s\in\Sigma^W$ fix a rational centre $w_\s\in W$.
Denote $r_\s=\frac{w_\s-r}{\pi^{\rho_\s}}$ for $r\in\roots$. Define reductions $\overline{f_\s^W}(x)\in k[x]$, $\overline{g_\s}\in k[y]$, and for $\s\in\Sigma$ also $\overline{g_\s^0}\in k[y]$ by
\begin{align*}
    \overline{f_\s^W}(x^{b_\s})&= \tfrac{u}{\pi^{v(u)}}\prod_{r\in\s\setminus\bigcup_{\s'<\s}\s'} (x+r_\s)\mod \pi,& u&=c_f\tprod_{r\in\roots\setminus\s} r_\s,\\
      \overline{g_\s}(y) &=  y^{p_\s/\gamma_\s} - \tfrac{u}{\pi^{v(u)}}\mod \pi, &  u&=c_f\tprod_{r\in\roots\setminus\s} r_\s,\\
  \overline{g_\s^0}(y) &=  y^{p_\s^0/\gamma_\s^0} - \tfrac{u}{\pi^{v(u)}}\mod \pi, & u&=c_f\tprod_{r\in\roots\setminus\{w_\s\}} r_\s.
\end{align*}
where the union runs through all $\s'\in\Sigma^W$, $\s'<\s$. 
Finally define the $k$-schemes
\begin{enumerate}
    \item $ X_\s^W:\{\ch{f_\s^W}=0\}\subset\G_{m,k}$;
    \item $ X_\s:\{\ch{g_\s}=0\}\subset\G_{m,k}$;
    \item $ X_\s^0:\{\ch{g_\s^0}=0\}\subset\G_{m,k}$ if $\s\in\Sigma$.
\end{enumerate}
\end{defn}


\begin{nt}\label{GeniricSpecialFibresNotation}
Given a scheme $\mathcal{X}/O_K$ we will denote by $\mathcal{X}_\eta$ its generic fibre $\mathcal{X}\times_{\Spec O_K}\Spec K$, and by $\mathcal{X}_s$ its special fibre $\mathcal{X}\times_{\Spec O_K}\Spec k$.
\end{nt}

\begin{nt}\label{ChainNotation}
If $C=C_1\cup\dots C_r$ is a chain of $\P^1_k$s of length $r$ and multiplicities $m_i\in\Z$ (meeting transversely), then $\infty\in C_i$ is identified with $0\in C_{i+1}$, and $0,\infty\in C$ are respectively $0\in C_{1}$ and $\infty\in C_{r}$. Finally, if $r=0$, then $C=\Spec k$ and $0=\infty$.
\end{nt}

\begin{nt}\label{ChainP1(alpha,a,b)Notation}
Let $\alpha,a,b\in\Z$, with $\alpha>0$, $a>b$, and fix $\frac{n_i}{d_i}\in\Q$ so that
    \[\alpha a=\frac{n_0}{d_0}>\frac{n_1}{d_1}>\ldots>\frac{n_r}{d_r}>\frac{n_{r+1}}{d_{r+1}}=\alpha b,\quad\text{with\scalebox{0.9}{$\quad\begin{vmatrix}n_i\!\!\!&n_{i+1}\cr d_i\!\!\!&d_{i+1}\cr\end{vmatrix}=1$}},\]
    and $r$ minimal. We write $\P^1(\alpha,a,b)$ for a chain of $\P^1_k$s of length $r$ and multiplicities $\alpha d_1,\dots,\alpha d_r$. Furthermore, we denote by $\P^1(\alpha,a)$ the chain $\P^1(\alpha,a,\lfloor\alpha a-1\rfloor/\alpha)$.
\end{nt}

Theorem \ref{ConstructionProperModelGeneralCaseTheorem} and Theorem \ref{MinimalRegularSNCModelTheorem} will follow from \S \ref{ConstructionModelsSection}.

\begin{thm}\label{ConstructionProperModelGeneralCaseTheorem}
Let $C/K$ be a hyperelliptic curve given by a Weierstrass equation $y^2=f(x)$. Suppose $\deg f>1$ and let $\Sigma$, $W$ and $\Sigma^W$ as in Definition \ref{SchemesXsDefinition}. Then there exists a proper flat model $\mathcal{C}/O_K$ (constructed in \S\ref{ConstructionModelsSection}) of $C$ such that its special fibre $\mathcal{C}_s/k$ consists of $1$-dimensional schemes given below in \ref{1},\ref{2},\ref{3},\ref{4},\ref{5}, glued along $0$-dimensional transversal intersections:
\begin{enumerate}[label=(\arabic*)]
    \item Every proper cluster $\s\in\Sigma^W$ gives a $1$-dimensional closed subscheme $\Gamma_\s$ of multiplicity $m_\s$. $\Gamma_\s$ is not integral if and only if $D_\s=2$, $\tilde\s\cap(\Sigma^W\cup\{\varnothing\})=\varnothing$ and $\ch{f_\s^W}$ is a square. When this happens, if $p=2$ then $\Gamma_\s$ is not reduced and $(\Gamma_\s)_\mathrm{red}$ is irreducible of multiplicity $2$ in $\Gamma_\s$, if $p\neq 2$ then $\Gamma_\s$ is reducible, namely $\Gamma_\s=\Gamma_\s^+\cup\Gamma_\s^-$, with $\Gamma_\s^{\pm}=\P^1_k$.\label{1}
    \item Every proper cluster $\s\in\Sigma^W$ with $D_\s=1$ gives the closed subscheme $X_\s^W\times\P^1_k$, of multiplicity $b_\s$, where $X_\s^W\times\{0\}\subset\Gamma_\s$.\label{2}
    \item Every proper cluster $\s\in\Sigma^W$ such that $\s\neq \roots$, 
    gives the closed subscheme $X_\s\times \P^1(\gamma_\s,s_\s,s_\s-p_\s\cdot\frac{\rho_\s-\rho_{P(\s)}}{2})$ where $X_\s\times\{0\}\subset\Gamma_\s$ and $X_\s\times\{\infty\}\subset\Gamma_{P(\s)}$.\label{3}
    \item Every cluster $\s\in\Sigma$ 
    gives the closed subscheme $X_\s^0\times \P^1(\gamma_\s^0,-s_\s^0)$ where $X_\s^0\times\{0\}\subset\Gamma_\s$ (the chains are open-ended). \label{4}
    \item Finally, 
    the cluster $\roots$ gives the closed subscheme $X_\roots\times \P^1(\gamma_\roots,s_\roots)$ where $X_\roots\times\{0\}\subset\Gamma_\s$ (the chains are open-ended). \label{5}
\end{enumerate}
If $\Gamma_\s$ is reducible, the two points in $X_\s\times\{0\}$ (and $X_\s^0\times\{0\}$ if $\s\in\Sigma$) belong to different irreducible components of $\Gamma_\s$. Similarly, if $\s\neq\roots$ and $\Gamma_{P(\s)}$ is reducible, the two points of $X_\s\times\{\infty\}$ belong to different irreducible components of $\Gamma_{P(\s)}$.

Furthermore, if $C$ has an almost rational cluster picture and is $y$-regular, then, by choosing $\Sigma=\Sigma_C^\mathrm{min}$, the model $\mathcal{C}$ is regular with strict normal crossings. In that case, if $\s$ is \"{u}bereven and $\epsilon_\s$ is even, then $\Gamma_\s\simeq X_\s\times\P^1_k$, otherwise $\Gamma_\s$ is irreducible of genus $g(\s)$. 
\end{thm}

\begin{defn} \label{RemovableContractibleDefinition}
Let $\s\in\Sigma_{K^{nr}}$. We say that
\begin{itemize}
    \item $\s$ is \textit{removable} if either $|\s|=1$, or $\s$ has a child $\s'\in\Sigma_{K^{nr}}$ of size $2g+1$ ($\s=\roots$ in this case).
    \item $\s$ is \textit{contractible} if
    \begin{enumerate}
        \item $|\s|=2$ and $\rho_\s\notin\Z$, $\epsilon_\s$ odd, $\rho_{P(\s)}\leq\rho_\s-\frac{1}{2};$ or \label{1Contr}
        \item $\s=\roots$ of size $2g+2$, with a child $\s'\in\Sigma_{K^{nr}}$ of size $2g$, and $\rho_\s\notin\Z$, $v(c_f)$ odd, $\rho_{\s'}\geq\rho_\s+\frac{1}{2}$; or \label{2Contr}
        \item $\s=\roots$ of size $2g+2$, union of its $2$ odd proper children $\s_1,\s_2\in\Sigma_{K^{nr}}$, with $v(c_f)$ odd, $\rho_{\s_i}\geq\rho_\s+1$ for $i=1,2$. \label{3Contr}
    \end{enumerate}
\end{itemize}
\end{defn}

\begin{nt}\label{NonRemovableClustersNotation}
Write $\mathring\Sigma\subseteq\Sigma_{K^{nr}}$ for the subset of non-removable clusters.
\end{nt}

\begin{defn}\label{PolynomialsDefinition}
Choose rational centres $w_\s$ for every $\s\in\mathring\Sigma$, in such a way that 
$w_\s\in\s$ when $p_\s^0=1$, and $\sigma(w_\s)=w_{\sigma(\s)}$ for all $\sigma\in\Gal(K^{nr}/K)$.
Denote $r_\s=\frac{w_\s-r}{\pi^{\rho_\s}}$ for $r\in\roots$ and
define $\overline{g_\s}, \overline{g_\s^0}\in k^\mathrm{s}[y]$ as in Definition \ref{SchemesXsDefinition}, and $\overline{f_\s}(x)\in k^\mathrm{s}[x]$, by
\[
  x^{2-p_\s^0}\overline{f_\s}(x^{b_\s}) = \tfrac{u}{\pi^{v(u)}} \prod_{r\in\s\setminus\bigcup_{\s'<\s}\s'} (x+r_\s)\mod \pi,\quad u=c_f\tprod_{r\in\roots\setminus\s} r_\s,
\]
where the union runs through all $\s'\in\mathring\Sigma$, $\s'<\s$. Let $G_\s=\mathrm{Stab}_{G_K}(\s)$, $K_\s=\lb K^\mathrm{s}\rb^{G_\s}$, and let $k_\s$ be the residue field of $K_\s$. Then $\ch{f_\s}\in k_\s[x]$, $\ch{g_\s}\in k_\s[y]$, and for $\s$ minimal $\ch{g_\s^0}\in k_\s[y]$. 

Let $\s_0\in\mathring\Sigma$ be minimal and contained in $\s$. Denote $\mathring\s=\tilde\s\setminus\{\{r\}<\s\mid r\neq w_{\s_0}\}$. Note that $\mathring\s$ does not depend on the choice of $\s_0$. Define $\tilde{f_\s}\in k_\s[x]$ by
\[\tilde{f_\s}(x)=\prod_{\s'\in\mathring\s}\lb x-\ch{u_{\s',\s}}\rb\cdot\ch{f_\s}(x),\]
where $\ch{u_{\s',\s}}=\frac{w_{\s'}-w_\s}{\pi^{\rho_\s}}\mod\pi$ if $\s'\neq\varnothing$ and $\ch{u_{\s',\s}}=0$ otherwise.
\end{defn}

In the next theorem we describe the special fibre of the minimal regular model of $C$ with normal crossings. We use Definitions/Notations \ref{ClusterDepthDefinition}, \ref{ParentChildWedgeDefinition}, \ref{OddEvenUberevenClusterDefinition}, \ref{ClusterPictureDefinition}, \ref{RadiusRationalCentreDefinition}, \ref{RationalClusterPictureDefinition}, \ref{AlmostRationalDefinition}, \ref{QuantitiesForTheoremsOnModelsDefinition}, \ref{yRegularDefinition}, \ref{GenussDefinition}, \ref{ChainP1(alpha,a,b)Notation}, \ref{RemovableContractibleDefinition}, \ref{NonRemovableClustersNotation}, \ref{PolynomialsDefinition} in the statement. Note that a full description of the model is developed in \S\ref{ConstructionModelsSection}.

\begin{thm}[Minimal regular NC model]\label{MinimalRegularSNCModelTheorem}
Let $C/K: y^2=f(x)$ be a hyperelliptic curve of genus $\geq 1$. Suppose $C_{K^{nr}}$ has an almost rational cluster picture and is $y$-regular.
Then the minimal regular model with normal crossings $\mathcal{C}^\mathrm{min}/O_{K^{nr}}$ of $C$ has special fibre $\mathcal{C}_s^\mathrm{min}/k^\mathrm{s}$ described as follows:
\begin{enumerate}[label=(\arabic*)]
\item 
  Every $\s\in\mathring\Sigma$ gives a $1$-dimensional subscheme $\Gamma_\s$ of multiplicity 
  $m_\s$. If $\s$ is \"{u}bereven and $\epsilon_\s$ is even, then $\Gamma_\s$ is the 
  disjoint union of $\Gamma_\s^{r_{\s,-}}\simeq\P^1$ and $\Gamma_\s^{r_{\s,+}}\simeq\P^1$, otherwise 
  $\Gamma_\s$ is irreducible of genus $g(\s)$ (write $\Gamma_\s^{r_{\s,-}}=\Gamma_\s^{r_{\s,+}}=\Gamma_\s$ 
  in this case). The indices $r_{\s,-}$ and $r_{\s,+}$ are the roots of $\ch{g_\s}$ (where $r_{\s,-}=r_{\s,+}$ if $\deg \ch{g_\s}=1$). \label{reg1}
\item 
  Every $\s\in\mathring\Sigma$ with $D_\s=1$ gives
  open-ended $\P^1$s of multiplicity $b_\s$ from $\Gamma_\s$ indexed by roots of $\ch{f_\s}$.\label{reg2}
\item 
Every non-maximal element $\s\in\mathring\Sigma$ gives chains
$\P^1(\gamma_\s,s_\s,s_\s-p_\s\cdot\frac{\rho_\s-\rho_{P(\s)}}{2})$ from $\Gamma_\s$ to $\Gamma_{P(\s)}$ indexed by roots of $\ch{g_\s}$.
\label{reg3}
\item 
  Every minimal element $\s\in\mathring\Sigma$ gives open-ended chains 
  $\P^1(\gamma_\s^0,-s_\s^0)$ from $\Gamma_\s$ indexed by roots of $\ch{g_\s^0}$.
\label{reg4}
\item 
  The maximal element $\s\in\mathring\Sigma$ gives open-ended chains
  $\P^1(\gamma_\s,s_\s)$ from $\Gamma_\s$ indexed by roots of $\ch{g_\s}$.
\label{reg5}
\item 
Finally, blow down all $\Gamma_\s$ where $\s$ is a contractible cluster.
\label{reg6}
\end{enumerate}
In \ref{reg3} and \ref{reg5}, a chain indexed by $r$ goes from $\Gamma_\s^r$. 
In \ref{reg3} the chain indexed by $r_{\s,-}$ goes to $\Gamma_{P(\s)}^{r_{P(\s),-}}$ and the chain indexed by $r_{\s,+}$ goes to $\Gamma_{P(\s)}^{r_{P(\s),+}}$.

Before blowing down in \ref{reg6}, the components given in \ref{reg1}--\ref{reg5} describe the special fibre of a regular model of $C_{K^{nr}}$ with strict normal crossings.

The Galois group $G_k$ acts naturally, i.e.\ for every $\sigma\in G_k$, $\sigma(\Gamma_\s^r)=\Gamma_{\sigma(\s)}^{\sigma(r)}$, and similarly, on the chains.

If $\Gamma_\s$ is irreducible, then its function field is isomorphic to $k^\mathrm{s}(x)[y]$ with the relation $y^{D_\s}=\tilde{f_\s}(x)$.
\end{thm}

\begin{rem}
Note that if $\Gamma_\s$ or $\Gamma_{P(\s)}$ is reducible then $p_\s/\gamma_\s=2$.
\end{rem}

\begin{exa}\label{MinimalRegularModelExample}
Let $p$ be a prime number and let $a\in\Z_p$, $b\in\Z_p^\times$ such that the polynomial $x^2+ax+b$ is not a square modulo $p$. Let $C$ be the hyperelliptic curve over $\Q_p$ of genus $4$ given by the equation $y^2=f(x)$, where $f(x)=(x^6+ap^4x^3+bp^8)((x-p)^3-p^{11})$. 
In Example \ref{AlmostRationalExample}, we described the rational cluster picture of $C$ and proved that $C$ has an almost rational cluster picture. Recall that $\Sigma_C^\mathrm{rat}$ consists of $3$ clusters $\t_3,\t_4,\roots$ of size $6,3,9$ respectively such that $\t_3<\roots$ and $\t_4<\roots$. In particular, note that $\Sigma_{\Q_p^{nr}}=\Sigma_C^\mathrm{rat}$, and no cluster of $\Sigma_{\Q_p^{nr}}$ is removable, so $\mathring \Sigma=\Sigma_C^\mathrm{rat}$. The minimal elements of $\mathring\Sigma$ are then $\t_3$ and $\t_4$. 

We want to describe the special fibre of the minimal regular model with normal crossings $\mathcal{C}^\mathrm{min}$ of $C$. Compute the quantities in Definitions \ref{QuantitiesForTheoremsOnModelsDefinition} and \ref{GenussDefinition}, and the polynomials $\ch{f_\s}, \ch{g_\s}, \ch{g_\s^0}$ of Definition \ref{PolynomialsDefinition}, for any cluster in $\mathring\Sigma$:

\cellspacetoplimit4pt
\cellspacebottomlimit4pt
\begin{center}
\begin{tabular}{|Sc!{\vline width 1pt}Sc|Sc|Sc|Sc|Sc|Sc|Sc|Sc|Sc|Sc|Sc|Sc|Sc|Sc|Sc|}
\hline
      \!\!\! &\!\!\! $\rho_\s$ \!\!\! &\!\!\!  $b_\s$  \!\!\! &\!\!\!  $\epsilon_\s$   \!\!\! &\!\!\! $D_\s$ 
       \!\!\! &\!\!\!  $m_\s$    \!\!\! &\!\!\!  $p_\s$   \!\!\! &\!\!\!  $s_\s$  \!\!\! &\!\!\! $\gamma_\s$  
       \!\!\! &\!\!\! $p_\s^0$   \!\!\! &\!\!\! $s_\s^0$  \!\!\! &\!\!\! $\gamma_\s^0$ \!\!\! &\!\!\! $g(\s)$ \!\!\! &\!\!\! $\ch{f_{\s}}(x)$ \!\!\! &\!\!\! $\ch{g_\s}(y)$ \!\!\! &$\ch{g_\s^0}(y)$\\ \Xhline{1 pt}
    $\t_3$ \!\!\! &\!\!\! $\tfrac{4}{3}$ \!\!\! &\!\!\! $3$ \!\!\! &\!\!\! $11$ \!\!\! &\!\!\! $1$ \!\!\! &\!\!\! $6$ \!\!\! &\!\!\! $2$ \!\!\! &\!\!\!
        $-\tfrac{1}{6}$ \!\!\! &\!\!\! $2$ \!\!\! &\!\!\! $2$ \!\!\! &\!\!\! $-\tfrac{25}{6}$ \!\!\! &\!\!\! $2$ \!\!\! &\!\!\! $0$ \!\!\! &\!\!\! $x^2+\bar ax+\bar b$ \!\!\! &\!\!\! $y+1$ \!\!\! &\!\!\! $y-1$\\ \hline
    $\t_4$ \!\!\! &\!\!\! $\tfrac{11}{3}$ \!\!\! &\!\!\! $3$ \!\!\! &\!\!\! $17$ \!\!\! &\!\!\! $1$ \!\!\! &\!\!\! $6$ \!\!\! &\!\!\! $1$ \!\!\! &\!\!\!
        $-\tfrac{7}{6}$ \!\!\! &\!\!\! $1$ \!\!\! &\!\!\! $2$ \!\!\! &\!\!\! $-\tfrac{29}{6}$ \!\!\! &\!\!\! $2$ \!\!\! &\!\!\! $0$ \!\!\! &\!\!\! $x-1$ \!\!\! &\!\!\! $y-1$ \!\!\! &\!\!\! $y+1$ \\ \hline
    $\roots$ \!\!\! &\!\!\! $1$ \!\!\! &\!\!\! $1$ \!\!\! &\!\!\! $9$ \!\!\! &\!\!\! $1$ \!\!\! &\!\!\! $2$ \!\!\! &\!\!\! $1$ \!\!\! &\!\!\! $\tfrac{1}{2}$
        \!\!\! &\!\!\! $1$ \!\!\! &\!\!\! $2$ \!\!\! &\!\!\!\!\!\! &\!\!\!\!\!\! &\!\!\! $0$ \!\!\! &\!\!\! $1$ \!\!\! &\!\!\! $y-1$ \!\!\! &\!\!\! \cr
\hline
\end{tabular}
\end{center}
where $\bar a, \bar b$ are the reductions of $a,b$ modulo $p$. Then $C$ is also $y$-regular for any $p$.
Following the steps of Theorem \ref{MinimalRegularSNCModelTheorem} the special fibre of $\mathcal{C}^\mathrm{min}$ over $\bar \F_p$ can be described as follows:

\begin{enumerate}[leftmargin=25pt]

\item [\ref{reg1}] The clusters $\t_3,\t_4,\roots$ give $3$ irreducible components $\Gamma_{\t_3},\Gamma_{\t_4},\Gamma_{\roots}$ of genus $0$ of multiplicities $6,6,2$ respectively;

\item [\ref{reg2}] The cluster $\t_3$ gives $2$ open-ended $\P^1$s of multiplicity $3$ from $\Gamma_{\t_3}$, while $\t_4$ gives $1$ open-ended $\P^1$ of multiplicity $3$ from $\Gamma_{\t_4}$.

\item [\ref{reg3}] From $\gamma_{\t_3}s_{\t_3}=-\tfrac{1}{3}>-\tfrac{1}{2}>-1=\gamma_{\t_3}\big(s_{\t_3}-p_{\t_3}\cdot\tfrac{\rho_{\t_3}-\rho_{\roots}}{2}\big),$ 
the cluster $\t_3$ gives $1$ $\P^1$ of multiplicity $4$ from $\Gamma_{\t_3}$ to $\Gamma_{\roots}$. From
\[\hspace{25pt}\gamma_{\t_4}s_{\t_4}=-\tfrac{7}{6}>-\tfrac{6}{5}>-\tfrac{5}{4}>-\tfrac{4}{3}>-\tfrac{3}{2}>-2>-\tfrac{5}{2}=\gamma_{\t_3}\big(s_{\t_4}-p_{\t_4}\cdot\tfrac{\rho_{\t_4}-\rho_{\roots}}{2}\big)\]
the cluster $\t_4$ gives $1$ chain of $\P^1$s of multiplicities $5,4,3,2,1$ from $\Gamma_{\t_4}$ to $\Gamma_\roots$.

\item [\ref{reg4}] From $-\gamma_{\t_3}^0s_{\t_3}^0=\tfrac{25}{3}>8>7$ the cluster $\t_3$ gives $1$ open-ended $\P^1$ of multiplicity $2$ from $\Gamma_{\t_3}$. From $-\gamma_{\t_4}^0s_{\t_4}^0=\tfrac{29}{3}>\tfrac{19}{2}>9>8$, the cluster $\t_4$ gives $1$ open-ended chain of $\P_1$s of multiplicities $4,2$ from $\Gamma_{\t_4}$.

\item [\ref{reg5}] From $\gamma_{\roots}s_{\roots}=\tfrac{1}{2}>0>-1$, the cluster $\roots$ gives $1$ open-ended $\P^1$ of multiplicity $1$ from $\Gamma_\roots$.

\item [\ref{reg6}] There is no contractible cluster, so the components we considered in the steps above describe the special fibre of $\mathcal{C}^\mathrm{min}$ over $\bar \F_p$:

\begin{center}
    \pbox[c]{20cm}{
\begin{tikzpicture}[xscale=1,yscale=0.7,
  l1/.style={shorten >=-1.3em,shorten <=-0.5em,thick},
  l2/.style={shorten >=-0.3em,shorten <=-0.3em},
  lfnt/.style={font=\tiny},
  leftl/.style={left=-3pt,lfnt},
  rightl/.style={right=-3pt,lfnt},
  mainl/.style={scale=0.8,above left=-0.17em and -1.5em},
  mainleftl/.style={scale=0.8,above right=-0.17em and -1.5em},
  abovel/.style={above=-2.5pt,lfnt},
  facel/.style={scale=0.7,blue,below right=-0.5pt and 6pt},
  faceleftl/.style={scale=0.7,blue,below left=-0.5pt and 6pt},
  redbull/.style={red,label={[red,scale=0.6,above=-0.17]#1}}]
\draw[l1] (-.66,0.00)--(2.53,0.00) node[mainl] {2} node[facel] {$\Gamma_\roots$};
\draw[l2] (0.50,0.00)--node[rightl] {1} (0.50,0.66);
\draw[l1] (1.30,2.00)--(3.33,2.00) node[mainl] {6} node[facel] {$\Gamma_{\t_4}$};
\draw[l2] (1.96,0.00)--node[rightl] {1} (1.96,0.66);
\draw[l2] (1.30,0.66)--node[abovel] {2} (1.96,0.66);
\draw[l2] (1.30,0.66)--node[leftl] {3} (1.30,1.33);
\draw[l2] (1.30,1.33)--node[abovel] {4} (1.96,1.33);
\draw[l2] (1.96,1.33)--node[rightl] {5} (1.96,2.00);
\draw[l2] (1.30,2.00)--node[rightl] {3} (1.30,2.66);
\draw[l2] (2.76,2.00)--node[rightl] {4} (2.76,2.66);
\draw[l2] (2.10,2.66)--node[abovel] {2} (2.76,2.66);
\draw[l1] (-.60,0.66)--(-3.16,0.66) node[mainleftl] {6} node[faceleftl] {$\Gamma_{\t_3}$};
\draw[l2] (-.60,0.00)--node[rightl] {4} (-.60,0.66);
\draw[l2] (-1.20,0.66)--node[rightl] {3} (-1.20,1.33);
\draw[l2] (-2,0.66)--node[rightl] {3} (-2,1.33);
\draw[l2] (-2.70,0.66)--node[rightl] {2} (-2.70,1.33);
\end{tikzpicture}
}

\end{center}
\end{enumerate}

Finally, from the Galois action on the roots of the polynomials $\ch{f_\s}, \ch{g_\s}, \ch{g_\s^0}$, for $\s\in\mathring\Sigma$, we get that $G_k$ acts trivially if $x^2+\bar a x+\bar b$ is reducible in $\F_p$, while it swaps the two components of multiplicity $3$ intersecting $\Gamma_{\t_3}$ (coming from \ref{reg2}) otherwise.
\end{exa}


As application of Theorem \ref{MinimalRegularSNCModelTheorem} we suppose $k$ is finite of characteristic $p>2$ and $C$ is semistable of genus $g\geq 2$. In this setting \cite[Theorem 8.5]{D2M2} describes the minimal regular model of $C$ in terms of its cluster picture $\Sigma_C$. We compare that result with the one obtained from Theorem \ref{MinimalRegularSNCModelTheorem} (Corollary \ref{SemistableMinimalRegularModelCorollary}). 

First note that $C_{K^{nr}}$ is $y$-regular as $p\neq 2$. From \cite[Definition 1.7]{D2M2}, if $C$ is semistable then
\begin{enumerate}
    \item the extension $K(\roots)/K$ has ramification degree at most $2$;\label{1Semistable}
    \item every proper cluster is $\Gal(K^\mathrm{s}/K^{nr})$-invariant;\label{2Semistable}
    \item every principal cluster has $d_\s\in\Z$ and $\nu_\s\in 2\Z$.
\end{enumerate}
It follows from Corollary \ref{AlmostRationalTameCaseCorollary} that $C_{K^{nr}}$ has an almost rational cluster picture.

In fact, (\ref{1Semistable}) and (\ref{2Semistable}) imply $\rho_\s=d_\s$ and $\epsilon_\s=\nu_\s$ for any proper cluster $\s$ (Remark \ref{TameNestedClusterPictureRemark}). In particular, $\mathring\Sigma_{C_{K^{nr}}}^\mathrm{rat}=\mathring\Sigma_C$. We will then say that $\s\in\Sigma_C$ is non-removable if $\s$ is proper and non-removable as rational cluster in $\Sigma_{K^{nr}}$.

\begin{lem}\label{SemistableAlmostRationalLemma}
Suppose $k$ finite and $p>2$. Assume $C$ is semistable and let $\s\in\Sigma_C$ be a non-removable cluster. Then $d_\s\in\frac{1}{2}\Z$ and $\nu_\s\in\Z$. Moreover, $\s$ is contractible if and only if $d_\s\notin\Z$ or $\nu_\s\notin 2\Z$.
\proof
Let $\s\in\Sigma_C$ be a non-removable cluster. Since $K(\roots)/K$ has ramification degree at most $2$, then $d_\s\in\tfrac{1}{2}\Z$.

By Theorem \ref{MinimalRegularSNCModelTheorem} the multiplicity of the $1$-dimensional scheme $\Gamma_\s$ is $m_\s$. Furthermore, $\Gamma_\s$ is an irreducible component of the special fibre of the minimal regular model of $C$ if and only if $\s$ is not contractible. Therefore if $\s$ is not contractible, then $m_\s=1$, i.e.\ $D_\s=2$ and $b_\s=1$. It follows that $\nu_\s\in 2\Z$ and $d_\s\in\Z$. Suppose $\s$ contractible. Then either $d_\s\notin\Z$ (and $\nu_\s\in\Z$) or $\s=\roots$ of size $2g+2$, with $2$ odd rational children and $v(c_f)$ odd. We want to show that in the latter case, $\nu_\s$ is odd. By Lemma \ref{TwoChildIntegralRadiusLemma}, $d_\roots\in\Z$. Then $\nu_\roots=v(c_f)+(2g+2)d_\roots$ is odd.
%
%
%
%
\endproof
\end{lem}

Let $\s\in\Sigma_C$ be a non-removable cluster.
By Lemma \ref{SemistableAlmostRationalLemma}, if $\s$ is not contractible, then $2g(\s)+1$ or $2g(\s)+2$ equals the number of odd children of $\s$. In fact, this also holds when $\s$ is contractible since in that case $g(\s)=0$ and $\s$ has at most $2$ odd children. 

\begin{cor}[Minimal regular model (semistable reduction)]\label{SemistableMinimalRegularModelCorollary}
Suppose that $k$ is finite and $p>2$. Let $C/K$ be a semistable hyperelliptic curve of genus $g\geq 2$. The minimal regular model $\mathcal{C}^\mathrm{min}/O_{K^{nr}}$ of $C$ has special fibre $\mathcal{C}_s^\mathrm{min}/k^\mathrm{s}$ described as follows:
\begin{enumerate}
\item Every non-removable cluster $\s\in\Sigma_C$ gives a $1$-dimensional subscheme $\Gamma_\s$. 
If $\s$ is \"{u}bereven, then $\Gamma_\s$ is the disjoint union of $\Gamma_\s^{r_{\s,-}}\simeq\P^1$ and $\Gamma_\s^{r_{\s,+}}\simeq\P^1$, otherwise $\Gamma_\s$ is irreducible of genus $g(\s)$ (write $\Gamma_\s^{r_{\s,-}}=\Gamma_\s^{r_{\s,+}}=\Gamma_\s$ in this case). 
The indices $r_{\s,-}$ and $r_{\s,+}$ are the roots of $\ch{g_\s}$.
\item Every odd proper cluster $\s\in\Sigma_C$, with size $|\s|\leq 2g$, gives a chain of $\P^1$s of length $\big\lfloor\frac{d_\s-d_{P(\s)}-1}{2}\big\rfloor$ from $\Gamma_\s$ to $\Gamma_{P(\s)}$ indexed by the root of $\ch{g_\s}$.
\item Every even proper cluster $\s\in\Sigma_C$, with size $|\s|\leq 2g$, gives a chain of $\P^1$s of length $\left\lfloor d_\s-d_{P(\s)}-\frac{1}{2}\right\rfloor$ from $\Gamma_\s^{r_{\s,-}}$ to $\Gamma_{P(\s)}^{r_{P(\s),-}}$ indexed by $r_{\s,-}$ and a chain of $\P^1$s of same length from $\Gamma_\s^{r_{\s,+}}$ to $\Gamma_{P(\s)}^{r_{P(\s),+}}$ indexed by $r_{\s,+}$.
\item Finally, blow down all $\Gamma_\s$ where $\s$ is a contractible cluster.
\end{enumerate}
All components have multiplicity $1$, and the absolute Galois group $G_k$ acts naturally, as in Theorem \ref{MinimalRegularSNCModelTheorem}.
\proof
Let $\s\in\Sigma_C$ be a non-removable cluster. From Lemma \ref{SemistableAlmostRationalLemma}, if $\s$ is not contractible, then $D_\s=2$, $\gamma_\s s_\s\in\Z$ and $\gamma_\s^0 s_\s^0\in\Z$. Suppose $\s$ contractible. If $|\s|=2$ with $d_\s\notin\Z$ (case (\ref{1Contr}) of Definition \ref{RemovableContractibleDefinition}), then $\gamma_\s^0 s_\s^0\in\Z$ and $\gamma_\s=1$, $s_\s\in\tfrac{1}{2}\Z\setminus\Z$ and so $s_\s-d_\s+d_{P(\s)}\in\Z$, as $P(\s)$ can not be contractible. If $\s=\roots$ (cases (\ref{2Contr}), (\ref{3Contr}) of Definition \ref{RemovableContractibleDefinition}), then $v(c_f)$ is odd, and so $\gamma_\s=2$ and $\gamma_\s s_\s\in\Z$.
Therefore \ref{reg2}, \ref{reg4} and \ref{reg5} of Theorem \ref{MinimalRegularSNCModelTheorem} do not give any components.

Finally, as $\gamma_\s=1$ and $p_\s\frac{d_\s-d_{P(\s)}}{2}\in\frac{1}{2}\Z$ for any proper $\s$ with size $|\s|\leq 2g$ (i.e.\ non-maximal), the length of $\P^1(\gamma_\s,s_\s,s_\s-p_\s\cdot\frac{d_\s-d_{P(\s)}}{2})$ 
is
\[\left\lfloor\gamma_\s s_\s-\gamma_\s\lb s_\s-p_\s\cdot\frac{d_\s-d_{P(\s)}}{2}\rb-\frac{1}{2}\right\rfloor=\left\lfloor p_\s\cdot\frac{d_\s-d_{P(\s)}}{2}-\frac{1}{2}\right\rfloor.\]
The corollary then follows from Theorem \ref{MinimalRegularSNCModelTheorem}.
\endproof
\end{cor}

\section{Construction of the model}\label{ConstructionModelsSection}
We are going to construct a proper flat model $\mathcal{C}/O_K$ of $C$ by glueing models defined in \cite[\S 4]{Dok}. For this reason we will assume the reader has familiarity with the definitions and the results presented in that paper.  Let us start this section by describing the strategy we will follow. 

Let $\Sigma_C^\mathrm{min}$ be the set of rationally minimal clusters of $C$ and let $\Sigma\subseteq\Sigma_C^\mathrm{min}$. For any cluster $\s\in\Sigma$ fix a rational centre $w_\s$ in such a way that $\mathring{\Sigma}_C^{w_\s}$ consists of the proper clusters in $\Sigma_C^{w_\s}$. 
This requirement can be satisfied by choosing $w_\s\in\s$ whenever possible.\footnote{This is the assumption used in Theorem \ref{ConstructionProperModelGeneralCaseTheorem}.} 
Let $W$ be the set of all such rational centres and define $\Sigma^W:=\bigcup_{w\in W}\Sigma_C^{w}$. For every proper cluster $\t\in\Sigma^W$ fix a rational centre $w_\t\in W$ (Lemma \ref{rationalCentresRadiusChildParentLemma}). For every $w\in W$, consider the curve $C^{w}:y^2=f(x+w)$, isomorphic to $C$, and construct the (proper flat) model $\mathcal{C}_\Delta^{w}/O_K$ by \cite[\S 4, Theorem 3.14]{Dok}. We will define an open subscheme $\mathring{\mathcal{C}}_\Delta^{w}$ of $\mathcal{C}_\Delta^{w}$ and we will show that glueing the schemes $\mathring{\mathcal{C}}_\Delta^{w}$, to varying of $w\in W$, along common opens, gives a proper flat model $\mathcal{C}/O_K$ of $C$. Furthermore, if $\Sigma=\Sigma_C^\mathrm{min}$, and $C$ is $y$-regular and has an almost rational cluster picture, then $\mathring{\mathcal{C}}_\Delta^{w}$ is an open regular subscheme of $\mathcal{C}_\Delta^{w}$ and therefore $\mathcal{C}$ is also regular.

\subsection{Charts}\label{ChartsSubsection}
Let $\Sigma=\{\s_1\dots,\s_m\}\subseteq\Sigma_C^\mathrm{min}$ be a set of rationally minimal clusters and let $W=\{w_1,\dots,w_m\}$ be a set of corresponding rational centres, such that $\mathring{\Sigma}_C^{w_h}$ consists of the proper clusters of $\Sigma_C^{w_h}$, for any $h=1,\dots,m$. Define $\Sigma^W:=\bigcup_{h=1}^m\Sigma_C^{w_h}$. 
For any $h,l=1,\dots,m$, $h\neq l$, define $w_{hl}:=w_h-w_l$, and write $w_{hl}=u_{hl}\pi^{\rho_{hl}}$, where $u_{hl}\in O_K^\times$ and $\rho_{hl}\in\Z$. Note that $\rho_{hl}=\rho_{\s_h\wedge\s_l}=\rho_{lh}$, by Lemma \ref{TwoChildIntegralRadiusLemma}. Set $u_{hh}:=0$. Finally, for any $h,l=1,\dots,m$, denote by $\overline{u_{hl}}\in k$ the reduction of $u_{hl}$ modulo $\pi$.

 \begin{defn}\label{MatricesDefinition}
Let $h=1,\dots,m$ and let $\t\in\Sigma_C^{w_h}$ be a proper cluster. Recall the matrices and cones defined in \cite[\S4]{Dok}. We say that a matrix $M$ \textit{is associated to} $\t$ if $M=M_{L_\t^{w_h},i}$ or $M=M_{V_\t^{w_h},j}$ (or $M=M_{V_0^{w_h},j}$ if $\t=\s_h$). For a matrix $M$ associated to $\t$ we denote by $\delta_M$ and $\sigma_M$ respectively 
\begin{itemize}
    \item the denominator $\delta_{L_\t^{w_h}}$ and the cone $\sigma_{L_\t^{w_h},i,i+1}$ if $M=M_{L_\t^{w_h},i}$, 
    \item the denominator $\delta_{V_\t^{w_h}}$ and the cone $\sigma_{V_\t^{w_h},j,j+1}$ if $M=M_{V_\t^{w_h},j}$,
    \item the denominator $\delta_{V_0^{w_h}}$ and the cone $\sigma_{V_0^{w_h},j,j+1}$ if $M=M_{V_0^{w_h},j}$.
\end{itemize}
Finally, define $X_M=\Spec O_K[\sigma_M^\vee\cap\Z^3]$ and write \[X_\Delta^h=\bigcup_{M}X_M,\]
for the toric scheme defined in \cite[\S 4.2]{Dok}.
\end{defn}

The following lemma describes all possible matrices associated to $\t$.

\begin{lem}\label{MatricesLemma}
Let $\t\in\Sigma_C^{w_h}$ be a proper cluster. Consider the $v$-face $F_\t^{w_h}$. Let $P_0,P_1\in\Z^2$ and $n_i,d_i,k_i\in\Z$ be as in \cite[\S 4]{Dok} and define
\[\delta:=\delta_M,\quad\gamma_i:=\frac{n_0}{\delta d_0}-\frac{n_{i}}{\delta d_{i}}\quad\mbox{ and }\quad T_i=\begin{psmallmatrix}
    \frac{1}{\delta} & -d_{i+1}k_i& d_ik_{i+1}\\
    0 & \delta d_{i+1} & 0\\
    0 & 0 & \delta d_i
    \end{psmallmatrix},\footnote{$\gamma_i$ and $T_i$ should be treated as placeholders for their respective definitions. Thus, for example, $\delta d_id_{i+1}(\gamma_{i+1}-\gamma_i)=1$ even when $d_{i+1}=0$.}\]
for each matrix $M$ associated to $\t$.
\begin{itemize}
    \item 
    Let $c$ be the unique element of $\{0,\dots,b_\t-1\}$ such that $\tfrac{1}{b_\t}-\rho_\t\cdot c=d\in\Z$. For all $i=0,\dots,r_{L_\t^{w_h}}$, choose $k_i=cn_i+d\delta d_i(\lfloor\t/2\rfloor+1)$; then
     \[\smaller{M_{L_\t^{w_h},i}=}\begin{psmallmatrix}
    \delta & -c\delta d_i\lb\tfrac{\epsilon_\t}{2}+\gamma_i\rb & c\delta d_{i+1}\lb\tfrac{\epsilon_\t}{2}+\gamma_{i+1}\rb\\
    0 & d_i & -d_{i+1}\\
    -\delta\rho_{\t} & -d\delta d_i\lb\tfrac{\epsilon_\t}{2}+\gamma_i\rb & d\delta d_{i+1}\lb\tfrac{\epsilon_\t}{2}+\gamma_{i+1}\rb
    \end{psmallmatrix}\smaller{,\quad M_{L_\t^{w_h},i}^{-1}=T_i \cdot} \begin{psmallmatrix}
    1 & \left\lfloor\frac{|\t|}{2}\right\rfloor+1 & 0\\  \rho_\t & \frac{\epsilon_\t}{2}+\gamma_{i+1} & 1\\
     \rho_\t & \frac{\epsilon_\t}{2}+\gamma_i & 1
    \end{psmallmatrix}\smaller{,}\]
    where $P_0=(|\t|,0)$, $P_1=(\lfloor\nicefrac{|\t|-1}{2}\rfloor,1)$ and $\delta=\delta_{L_{\t}^{w_h}}=b_\t$.
    \item If $\t$ is odd, then for all $j=0,\dots, r_{V_\t^{w_h}}$, we have
    \[\smaller{M_{V_\t^{w_h},j}=}\begin{psmallmatrix}
    -|\t| & -\frac{|\t|+1}{2} d_j & \frac{|\t|+1}{2}d_{j+1}\\
    2 & d_j &-d_{j+1}\\
    -\epsilon_\t+|\t|\rho_\t & n_j &-n_{j+1}
    \end{psmallmatrix} \smaller{,\quad M_{V_\t^{w_h},j}^{-1}=T_j}\cdot \begin{psmallmatrix}
    1 & \frac{|\t|+1}{2} & 0\\
    \rho_\t -2\cdot\gamma_{j+1} & \frac{\epsilon_\t}{2}-|\t|\cdot\gamma_{j+1}& 1\\
    \rho_\t -2\cdot\gamma_{j} &\frac{\epsilon_\t}{2}-|\t|\cdot\gamma_j& 1
    \end{psmallmatrix}\smaller,\]
    where $P_0=(|\t|,0)$, $P_1=(\lfloor\nicefrac{|\t|-1}{2}\rfloor,1)$, $\delta=\delta_{V_\t^{w_h}}=1$ and $k_j=k_{j+1}=0$.
    \item If $\t$ is even, then for all $j=0,\dots, r_{V_\t^{w_h}}$, we have
    \[\smaller{M_{V_\t^{w_h},j}=}\begin{psmallmatrix}
    -\delta\frac{|\t|}{2} & -\lb\frac{|\t|}{2}+1\rb d_j - k_j\frac{|\t|}{2} & \lb\frac{|\t|}{2}+1\rb d_{j+1}+k_{j+1}\frac{|\t|}{2}\\
    \delta & d_j+k_j &-d_{j+1}-k_{j+1}\\
    -\delta\frac{\epsilon_\t-|\t|\rho_\t}{2} & \frac{n_j}{\delta}-k_j\frac{\epsilon_\t-|\t|\rho_\t}{2} &-\frac{n_{j+1}}{\delta}+k_{j+1}\frac{\epsilon_\t-|\t|\rho_\t}{2}
    \end{psmallmatrix}\smaller{,}\]\[\smaller{M_{V_\t^{w_h},j}^{-1}=T_j\cdot }\begin{psmallmatrix}
    1 & \frac{|\t|}{2}+1 & 0\\
    \rho_\t -\gamma_{j+1} &  \frac{\epsilon_\t}{2}-\frac{|\t|}{2}\gamma_{j+1}& 1\\
    \rho_\t -\gamma_{j} & \frac{\epsilon_\t}{2} -\frac{|\t|}{2}\gamma_{j}& 1
    \end{psmallmatrix}\smaller{,}\]
    where $P_0=(|\t|,0)$, $P_1=(\lfloor\nicefrac{|\t|-1}{2}\rfloor,1)$ and $\delta=\delta_{V_\t^{w_h}}$.
    \item If $f(w_h)=0$, then for all $j=0,\dots, r_{V_0^{w_h}}$, we have
    \[ \smaller{M_{V_0^{w_h},j}=}\begin{psmallmatrix}
    1 & d_j & -d_{j+1}\\
    -2 & -d_j &d_{j+1}\\
    \epsilon_{\s_h}-\rho_{\s_h} & n_j &-n_{j+1}
    \end{psmallmatrix}\smaller{,\quad   M_{V_0^{w_h},j}^{-1}=T_j\cdot }\begin{psmallmatrix}
    -1 & -1 & 0\\
    \rho_{\s_h} +2\cdot\gamma_{j+1} & \frac{\epsilon_{\s_h}}{2}+\gamma_{j+1}& 1\\
    \rho_{\s_h} +2\cdot\gamma_{j} & \frac{\epsilon_{\s_h}}{2}+\gamma_{j}& 1
    \end{psmallmatrix}\smaller,\]
    where $P_0=(0,2)$, $P_1=(1,1)$, $\delta=\delta_{V_0^{w_h}}=1$ and $k_j=k_{j+1}=0$.
    \item If $f(w_h)\neq 0$, then for all $j=0,\dots, r_{V_0^{w_h}}$, we have
    \[\smaller{M_{V_0^w,j}=}\begin{psmallmatrix}
    0 & d_j & - d_{j+1}\\
    -\delta & -d_j-k_j &d_{j+1}+k_{j+1}\\
    \delta\frac{\epsilon_{\s_h}}{2} & \frac{n_j}{\delta}+k_j\frac{\epsilon_{\s_h}}{2} &-\frac{n_{j+1}}{\delta}-k_{j+1}\frac{\epsilon_{\s_h}}{2}
    \end{psmallmatrix}\smaller{,\quad M_{V_0^{w_h},j}^{-1}=T_j\cdot }\begin{psmallmatrix}
    -1 & -1 & 0\\
    \rho_{\s_h} +\gamma_{j+1} &  \frac{\epsilon_{\s_h}}{2}& 1\\
    \rho_{\s_h} +\gamma_{j} &  \frac{\epsilon_{\s_h}}{2}& 1
    \end{psmallmatrix}\smaller,\]
    where $P_0=(0,2)$, $P_1=(1,1)$ and $\delta=\delta_{V_0^{w_h}}$.
\end{itemize}
\proof
We follow the notation of \cite[\S 4]{Dok}. Choose the points $P_0$ and $P_1$ as in the proof of Lemma \ref{DescriptionNewtonPolytopeClustersLemma}. 

First consider the edge $L_\t^{w_h}$ of $F_\t^{w_h}$. From Lemma \ref{DescriptionNewtonPolytopeClustersLemma} we have \[\nu=(1,0,-\rho_t)\quad\mbox{ and }\quad (w_x,w_y)=\lb-\left\lfloor|\t|/2\right\rfloor-1,1\rb.\]
Then $M_{L_\t^{w_h},i}$ and $M_{L_\t^{w_h},i}^{-1}$ follow from \cite[\S 4.3]{Dok} as $k_i\equiv n_i(\delta\rho_\t)^{-1}\mod \delta$ and
\[\sfrac{n_0}{\delta d_0}=\sfrac{1}{\delta}s_1^{L_\t^{w_h}}=v_{F_\t^{w_h}}(P_1)-v_{F_\t^{w_h}}(P_0)=-\sfrac{\epsilon_\t}{2}+\lb\left\lfloor|\t|/2\right\rfloor+1\rb\rho_\t\]

Now assume $\t$ even and consider the edge $V_\t^{w_h}$ of $F_\t^{w_h}$. Since $\t$ is even,
\[V_\t^{w_h}(\Z)=\left\{(|\t|,0),\lb \sfrac{|\t|}{2},1\rb,(0,2)\right\},\quad \nu=\lb-\sfrac{|\t|}{2},1,-\sfrac{\epsilon_t}{2}+\sfrac{|\t|}{2}\rho_\t\rb\]
and $(w_x,w_y)=\lb-\frac{|\t|}{2}-1,1\rb$ as above. Then $M_{V_\t^{w_h},j}$ and $M_{V_\t^{w_h},j}^{-1}$ follow again from \cite[(4.3)]{Dok} as
\[\sfrac{n_0}{\delta d_0}=\sfrac{1}{\delta}s_1^{V_\t^{w_h}}=v_{F_\t^{w_h}}(P_1)-v_{F_\t^{w_h}}(P_0)=-\sfrac{\epsilon_\t}{2}+\lb\sfrac{|\t|}{2}+1\rb\rho_\t.\]
%

Similar arguments and computations yield the remaining matrices.
\endproof
\end{lem}

\begin{rem}
From the lemma above one can explicitly construct the charts of the model $\mathcal{C}_\Delta^{w_h}$. The description of its special fibre $\mathcal{C}_{\Delta,s}^{w_h}$ which follows from \cite[Theorem 3.14]{Dok}, matches the one given in Theorem \ref{ConstructionProperModelGeneralCaseTheorem} in the case $W=\{w_h\}$. 
\end{rem}

\subsection{Open subschemes}\label{OpenSubschemesSubsection}

Let $h=1,\dots,m$ and let $\t\in\Sigma_C^{w_h}$ be a proper cluster. Let $M$ be a matrix associated to $\t$. Write
\[M=\begin{pmatrix}
m_{11} & m_{12} & m_{13}\\
m_{21} & m_{22} & m_{23}\\
m_{31} & m_{32} & m_{33}
\end{pmatrix}\qquad\mbox{and}\qquad M^{-1}=\begin{pmatrix}
\tilde m_{11} &\tilde m_{12} & \tilde m_{13}\\
\tilde m_{21} & \tilde m_{22} & \tilde m_{23}\\
\tilde m_{31} & \tilde m_{32} & \tilde m_{33}
\end{pmatrix}\]
Recall that $X_M=\Spec R$, where
\[R=\frac{O_K[X^{\pm 1},Y,Z]}{\lb \pi-X^{\tilde m_{13}}Y^{\tilde m_{23}}Z^{\tilde m_{33}}\rb}\hookrightarrow\frac{O_K[X^{\pm 1},Y^{\pm 1},Z^{\pm 1}]}{\lb \pi-X^{\tilde m_{13}}Y^{\tilde m_{23}}Z^{\tilde m_{33}}\rb}\stackrel{M}{\simeq} K\left[x^{\pm 1},y^{\pm 1}\right],\]
via the change of variable
\[
    \begin{psmallmatrix}X\\ Y\\ Z\cr\end{psmallmatrix}=\begin{psmallmatrix}x^{m_{11}}y^{m_{21}}\pi^{m_{31}}\cr x^{m_{12}}y^{m_{22}}\pi^{m_{32}}\cr x^{m_{13}}y^{m_{23}}\pi^{m_{33}}\end{psmallmatrix}=\begin{psmallmatrix}x\\ y\\ z\end{psmallmatrix}\bullet M,
\quad
    \begin{psmallmatrix}x\\ y\\ \pi\end{psmallmatrix}=\begin{psmallmatrix}X^{\tilde m_{11}}Y^{\tilde m_{21}}Z^{\tilde m_{31}}\cr X^{\tilde m_{12}}Y^{\tilde m_{22}}Z^{\tilde m_{32}}\cr X^{\tilde m_{13}}Y^{\tilde m_{23}}Z^{\tilde m_{33}}\end{psmallmatrix}=\begin{psmallmatrix}X\\ Y\\ Z\end{psmallmatrix}\bullet M^{-1}.
\]
Let $l\neq h$.
Set
\[T_M^{hl}(X,Y,Z):=\begin{cases}1+u_{hl}X^{\rho_{hl}\tilde m_{13}-\tilde m_{11}}Y^{\rho_{hl}\tilde m_{23}-\tilde m_{21}}Z^{\rho_{hl}\tilde m_{33}-\tilde m_{31}}&\mbox{if }\t\supseteq\s_h\wedge\s_l,\\
u_{hl}^{-1}X^{\tilde m_{11}-\rho_{hl}\tilde m_{13}}Y^{\tilde m_{21}-\rho_{hl}\tilde m_{23}}Z^{\tilde m_{31}-\rho_{hl}\tilde m_{33}}+1&\mbox{if }\t\not\supseteq\s_h\wedge\s_l, \end{cases}\]
element of $R[Y^{-1},Z^{-1}]$.
Note that 
\begin{align*}
    &\mbox{if}\quad\t\supseteq\s_h\wedge\s_l\quad\mbox{then}\quad T_M^{hl}(X,Y,Z)\stackrel{M}{\longmapsto}\frac{x+w_{hl}}{x},\\&\mbox{if}\quad\t\not\supseteq\s_h\wedge\s_l\quad\mbox{then}\quad T_M^{hl}(X,Y,Z)\stackrel{M}{\longmapsto}\frac{x+w_{hl}}{w_{hl}}.
\end{align*}
The following two lemmas prove that $T_M^{hl}(X,Y,Z)\in R$. Therefore, up to units, $T_M^{hl}(X,Y,Z)$ can be seen as the polynomial in $O_K[X^{\pm 1},Y,Z]$ satisfying
\[x-w_{hl}\stackrel{M}{=}X^{n_X} Y^{n_Y} Z^{n_Z}T_M^{hl}(X,Y,Z),\]
with $n_X,n_Y,n_Z\in\Z$, such that $\ord_Y(T_M^{hl})=\ord_Z(T_M^{hl})=0$.

\begin{lem}\label{NonnegativeExponentsMatricesLemma}
Let $h,l=1,\dots,m$, with $h\neq l$, let $t\in \Sigma_C^{w_h}$ be such that $\t\supseteq\s_h\wedge\s_l$ and let $M$ be a matrix associated to $\t$.
Then \[\rho_{hl} \tilde m_{23}-\tilde m_{21}\geq \rho_\t \tilde m_{23}-\tilde m_{21}\geq 0 \quad\mbox{and}\quad \rho_{hl} \tilde m_{33}- \tilde m_{31}\geq\rho_\t \tilde m_{33}-\tilde m_{31}\geq 0.\]
Furthermore if $M=M_{L_\t^{w_h},i}$ then 
\begin{itemize}
    \item $\rho_{hl} \tilde m_{23}-\tilde m_{21}=0$ if and only if $i=r_{L_{\t}^{w_h}}$ or $\t=\s_h\wedge\s_l$,
    \item $\rho_{hl} \tilde m_{33}-\tilde m_{31}=0$ if and only if $\t=\s_h\wedge\s_l$;
\end{itemize}
if $M=M_{V_\t^{w_h},j}$ then
\begin{itemize}
    \item $\rho_{hl} \tilde m_{23}-\tilde m_{21}>0$,
    \item $\rho_{hl} \tilde m_{33}-\tilde m_{31}=0$ if and only if $\t=\s_h\wedge\s_l$ and $j=0$.
\end{itemize}
\proof
This result follows from Lemma \ref{MatricesLemma}, which gives a complete description of $M$ and $M^{-1}$.
We show it when 
%
%
$\t$ is even and $M=M_{V_\t^{w_h},j}$, and leave the other cases for the reader. First of all recall that $\rho_{hl}=\rho_{\s_h\wedge\s_l}$ by Lemma \ref{TwoChildIntegralRadiusLemma}. Then
\[\rho_{hl} \tilde m_{23}-\tilde m_{21}= \delta d_{j+1}\lb\rho_{hl}-\rho_t+\gamma_{j+1}\rb> \delta d_{j+1}\lb\rho_{\s_h\wedge\s_l}-\rho_\t\rb\geq 0,\]
where $\gamma_{j}=\tfrac{n_0}{\delta d_0}-\tfrac{n_{j}}{\delta d_{j}}$ and $\delta=\delta_M$. 
Similarly,
\[\rho_{hl} \tilde m_{33}-\tilde m_{31}= \delta d_{j}\lb\rho_{hl}-\rho_t+\gamma_j\rb\geq \delta d_{j}\lb\rho_{\s_h\wedge\s_l}-\rho_\t\rb\geq 0.\]
In particular $\rho_{hl} \tilde m_{33}-\tilde m_{31}=0$ if and only if $\t=\s_h\wedge\s_l$ and $j=0$.
\endproof
\end{lem}

\begin{lem}\label{negativeexponentsLemma}
Let $\t\in\Sigma_C^{w_h}$ be a proper cluster such that $\t\not\supseteq\s_h\wedge\s_l$, and let $M$ be a matrix associated to $\t$. Then
\[\tilde m_{21}-\rho_{hl} \tilde m_{23}\geq 0\quad\mbox{and}\quad\tilde m_{31}-\rho_{hl} \tilde m_{33}>
0.\]
Furthermore, $\tilde m_{21}-\rho_{hl} \tilde m_{23}=0$ if and only if
\begin{itemize}
    \item $M=M_{L_\t^{w_h},i}$ and $i=r_{L_\t^{w_h}}$, or
    \item $\t<\s_h\wedge\s_l$, $M=M_{V_\t^{w_h},j}$, and $j=r_{V_\t^{w_h}}$.
\end{itemize}
\proof
This result follows again from Lemma \ref{MatricesLemma}. As in the previous lemma, we show it when $\t$ is even and $M=M_{V_\t^{w_h},j}$, and leave the other cases for the reader.


Let 
$r=r_{V_\t^{w_h}}$. Note that $t\neq\roots$. Set $\delta=\delta_M$ and $\gamma_j=\tfrac{n_0}{\delta d_0}-\tfrac{n_{j}}{\delta d_{j}}$. Then
\[\tilde m_{31}-\rho_{hl} \tilde m_{33}= \delta d_{j}\lb\rho_t-\rho_{hl}-\gamma_j\rb> \delta d_{j}\lb\rho_{P(\t)}-\rho_{\s_h\wedge\s_l}\rb\geq 0.\]
since
$\gamma_j<\gamma_{r+1}=\rho_\t-\rho_{P(\t)}.$
Similarly,
\[\tilde m_{21}-\rho_{hl} \tilde m_{23}= \delta d_{j+1}\lb\rho_t-\rho_{hl}-\gamma_{j+1}\rb\geq \delta d_{j+1}\lb\rho_{P(\t)}-\rho_{\s_h\wedge\s_l}\rb\geq 0,\]
In particular $\tilde m_{21}-\rho_{hl} \tilde m_{23}=0$ if and only if $\t<\s_h\wedge\s_l$ and $j=r$.
%
%
\endproof
\end{lem}

Let \[T_M^h(X,Y,Z):=\prod_{l\neq h}T_M^{hl}(X,Y,Z),\]
and define 
\[V_M^h:=\Spec R[T_M^h(X,Y,Z)^{-1}]\subset X_M,\quad\mbox{ and }\quad \mathring{X}_\Delta^h:=\bigcup_{\t,M} V_M^h\subseteq X_\Delta^h,\]
where $\t$ runs through all proper clusters in $\Sigma_C^{w_h}$ and $M$ runs through all matrices associated to $\t$. We can then define the subscheme
\[\mathring{\mathcal{C}}_\Delta^{w_h}:= \mathcal{C}_{\Delta}^{w_h}\cap \mathring{X}_\Delta^h\subset X_\Delta^h,\]
where $\mathcal{C}_{\Delta}^{w_h}/O_K$ is the model of the hyperelliptic curve $C^{w_h}:y^2=f(x+w_h)$ described in \cite[Theorem 3.14]{Dok} (see \cite[\S 4]{Dok} for the construction). Explicitly, let $g_h(x,y):=y^2-f(x+w_h)$ and define $\mathcal{F}_M^h\in O_K[X^{\pm 1},Y,Z]$ such that $\ord_Y(\mathcal{F}_M^h)=\ord_Z(\mathcal{F}_M^h)=0$, with all non-zero coefficients in $O_K^\times$, satisfying
 \[y^2-f(x+w_h)\stackrel{M}{=}Y^{n_{Y,h}}Z^{n_{Z,h}}\mathcal{F}_M^h(X,Y,Z),\]
for unique $n_{Y,h},n_{Z,h}\in\Z$. Consider the subscheme
\[U_M^h:=\Spec\frac{R\left[T_M^h(X,Y,Z)^{-1}\right]}{\lb\mathcal{F}_M^h(X,Y,Z)\rb}\subset V_M^h.\]
Then
\[\mathring{\mathcal{C}}_\Delta^{w_h}=\bigcup_{\t,M} U_M^h\subset \mathring{X}_\Delta^h,\]
where $\t$ runs through all proper clusters in $\Sigma_C^{w_h}$ and $M$ runs through all matrices associated to $\t$, as before.

\subsection{Glueing}\label{GlueingSubsection}

Let $h,l=1,\dots,m$, with $h\neq l$. Consider the isomorphism
\begin{equation}\label{GlueingMap(-1)equation}
    \phi:K\left[x^{\pm 1},y^{\pm 1},\prod_{o\neq l}(x+w_{lo})^{-1}\right]\stackrel{\simeq}{\longrightarrow} K\left[x^{\pm 1},y^{\pm 1},\prod_{o\neq h}(x+w_{ho})^{-1}\right]
\end{equation}
sending $x\mapsto x+w_{hl}$, $y\mapsto y$. If $\t\supseteq\s_h\wedge\s_l$ and $M$ is a matrix associated to $\t$, then $\phi$ gives a map
\[ R[Y^{-1},Z^{-1},T_M^l(X,Y,Z)^{-1}]\xrightarrow{M^{-1}\circ\phi\circ M}R[Y^{-1},Z^{-1},T_M^h(X,Y,Z)^{-1}],\]
which sends
 \[F(X,Y,Z)\mapsto F(X\cdot T_M^{hl}(X,Y,Z)^{m_{11}}, Y\cdot T_M^{hl}(X,Y,Z)^{m_{12}},Z\cdot T_M^{hl}(X,Y,Z)^{m_{13}}).\]
Hence it induces the isomorphisms
 \begin{equation}\label{IsomorphismAmbientSpacesGeneralCaseequation}
     R[T_M^l(X,Y,Z)^{-1}]\stackrel{\simeq}{\longrightarrow} R[T_M^h(X,Y,Z)^{-1}], \qquad V_M^h\stackrel{\simeq}{\longrightarrow}V_M^l.
 \end{equation}
Via these maps we see that $g_h(x,y)=Y^{n_{Y,h}}Z^{n_{Z,h}}\mathcal{F}_M^h(X,Y,Z)$ also equals
\[ Y^{n_{Y,l}}\cdot Z^{n_{Z,l}}\cdot (T_M^{hl})^{n_{Y,l}m_{12}+n_{Z,l}m_{13}}\mathcal{F}_M^l\lb X\cdot (T_M^{hl})^{m_{11}},Y\cdot (T_M^{hl})^{m_{12}},Z\cdot (T_M^{hl})^{m_{13}}\rb,\]
where $T_M^{hl}=T_M^{hl}(X,Y,Z)$. Since neither $Y$ nor $Z$ divide $T_M^{hl}(X,Y,Z)$, we have $n_{Y,h}=n_{Y,l}$, $n_{Z,h}=n_{Z,l}$ and
\[\mathcal{F}_M^h(X,Y,Z)=(T_M^{hl})^{n_{Y,l}m_{12}+n_{Z,l}m_{13}}\mathcal{F}_M^l\lb X\, (T_M^{hl})^{m_{11}},Y\, (T_M^{hl})^{m_{12}},Z\, (T_M^{hl})^{m_{13}}\rb.\] 
Hence (\ref{IsomorphismAmbientSpacesGeneralCaseequation}) induces the isomorphisms
\begin{equation}\label{IsomorphismsModelsGeneralCaseEquation}
    \frac{R\left[T_M^l(X,Y,Z)^{-1}\right]}{\lb\mathcal{F}_M^l(X,Y,Z)\rb}\stackrel{\simeq}{\longrightarrow}\frac{R\left[T_M^h(X,Y,Z)^{-1}\right]}{\lb\mathcal{F}_M^h(X,Y,Z)\rb},\qquad U_M^h\stackrel{\simeq}{\longrightarrow}U_M^l.
\end{equation}

Define the subschemes
\[V^{hl}:=\bigcup_{\t_l,M_l}V_{M_l}^{h}\subseteq \mathring{X}_{\Delta}^h,\qquad U^{hl}:=V^{hl}\cap\mathcal{C}_{\Delta}^{w_h}\subseteq \mathring{\mathcal{C}}_{\Delta}^{w_h},
\]
where $\t_l$ runs through all proper clusters in $\Sigma_C^{w_h}\cap \Sigma_C^{w_l}$ (i.e.\ $\t_l\in\Sigma^W$, $\s_h\wedge\s_l\subseteq\t_l$) and $M_l$ runs through all matrices associated to $\t_l$. 
From (\ref{GlueingMap(-1)equation}), (\ref{IsomorphismAmbientSpacesGeneralCaseequation}) and (\ref{IsomorphismsModelsGeneralCaseEquation}) we have isomorphisms of schemes
\begin{equation}\label{GlueingMapEquation}
    V^{hl}\stackrel{\simeq}{\longrightarrow}V^{lh},\qquad U^{hl}\stackrel{\simeq}{\longrightarrow}U^{lh}.
\end{equation}
Now, $U^{hl}\subset V^{hl}$ are open subschemes respectively of $\mathring{\mathcal{C}}_\Delta^{w_h}\subset \mathring{X}_\Delta^h$ for any $l\neq h$. Glueing the schemes $\mathring{\mathcal{C}}_\Delta^{w_h}\subset \mathring{X}_\Delta^h$, to varying of $h=1,\dots,m$, respectively along the opens $U^{hl}\subset V^{hl}$ via (\ref{GlueingMapEquation}) gives the schemes $\mathcal{C}\subset \mathcal{X}$. We will show that $\mathcal{C}/O_K$ is a proper flat\footnote{Note that the flatness of $\mathcal{C}$ is trivial since it is a local property.} model of $C$.

\subsection{Generic fibre} 
We start studying the generic fibre $\mathcal{C}_\eta$ of $\mathcal{C}$. Since it is the glueing of all $\mathring{\mathcal{C}}_{\Delta,\eta}^{w_h}$ through the glueing maps 
\[U^{hl}_\eta\longrightarrow U^{lh}_\eta\]
induced by (\ref{GlueingMapEquation}), we start focusing on $\mathring{\mathcal{C}}_{\Delta,\eta}^{w_h}$ for $h=1,\dots,m$. In particular, as $\mathring{\mathcal{C}}_\Delta^{w_h}$ is an open subscheme of $\mathcal{C}_\Delta^{w_h}$, we study 
$\mathcal{C}_{\Delta,\eta}^{w_h}\smallsetminus \mathring{\mathcal{C}}_{\Delta,\eta}^{w_h}=C^{w_h}\smallsetminus\mathring{\mathcal{C}}_{\Delta,\eta}^{w_h}$.

\begin{lem}
For any $h=1,\dots,m$,
\[C^{w_h}\smallsetminus\mathring{\mathcal{C}}_{\Delta,\eta}^{w_h}=\Spec\frac{K[x,y]}{\lb g_h(x,y), \prod_{o\neq h}\lb x+w_{ho}\rb\rb}.\]
\proof
For every choice of a proper cluster $\t\in\Sigma_C^{w_h}$, and $M$ associated to $\t$, let \[P_M:=\lb \mathcal{C}_{\Delta,\eta}^{w_h}\smallsetminus\mathring{\mathcal{C}}_{\Delta,\eta}^{w_h}\rb\cap X_M=\Spec \frac{R\otimes_{O_K} K}{\lb\mathcal{F}_M^h(X,Y,Z),T_M^h(X,Y,Z)\rb}.\]
To study $P_M$ we are going to use Lemma \ref{MatricesLemma} and the definition of $T_M^h(X,Y,Z)$.

Suppose first $\t\neq\roots$ and $M=M_{V_\t^{w_h},j}$. Then $\tilde m_{23},\tilde m_{33}> 0$, so
\begin{equation}\label{PMEqualityequation}
    P_M=\Spec \frac{R[Y^{- 1}, Z^{- 1}]}{\lb\mathcal{F}_M^h(X,Y,Z),T_M^h(X,Y,Z)\rb}\stackrel{M}{\simeq} \Spec\frac{K[x^{\pm 1},y^{\pm 1}]}{\lb g_h(x,y), \prod_{o}\lb x+w_{ho}\rb\rb},
\end{equation}
where the product runs over all $o\neq h$. Now let $\t=\roots$ and $M=M_{V_\t^{w_h},j}$. If $j\neq r_{V_\roots^{w_h}}$, then $P_M$ is as in the previous case (since $\tilde m_{23},\tilde m_{33}> 0$). If $j= r_{V_\roots^{w_h}}$, then $\tilde m_{33}> 0$, $\tilde m_{23}=0$, but $\rho_{hl} \tilde m_{23}-\tilde m_{21}>0$ by Lemma \ref{NonnegativeExponentsMatricesLemma}. So from the definition of $T_M^{hl}(X,Y,Z)$ we have once more the equality (\ref{PMEqualityequation}). Similarly, if $\t=\s_h$ and $M=M_{V_0^{w_h},j}$, then $\tilde m_{33}> 0$, and $\tilde m_{21}-\rho_{hl} \tilde m_{23}>0$ by Lemma \ref{negativeexponentsLemma}. Hence we have (\ref{PMEqualityequation}) again.

It remains to study $P_M$ when $M=M_{L_\t^{w_h},i}$. If $i\neq r_{L_\t^{w_h}}$, then $\tilde m_{23},\tilde m_{33}> 0$ and so $P_M$ is as in (\ref{PMEqualityequation}). Let $i=r_{L_\t^{w_h}}$. Then $\tilde m_{33}> 0$ but both $\tilde m_{23}$ and $\rho_{hl} \tilde m_{23}-\tilde m_{21}$ equal $0$. Hence $\tilde m_{23}=\tilde m_{21}=0$, which also implies $m_{21}=m_{23}=0$. Therefore $M$ defines an isomorphism $R[Z^{-1}]\simeq K[x^{\pm 1},y],$
which induces
\[P_M=\Spec \frac{R[Z^{- 1}]}{\lb\mathcal{F}_M^h(X,Y,Z),T_M^h(X,Y,Z)\rb}\stackrel{M}{\simeq} \Spec\frac{K[x^{\pm 1},y]}{\lb g_h(x,y), \prod_{o\neq h}\lb x+w_{ho}\rb\rb}.\]

This concludes the proof.
\endproof
\end{lem}

Regarding $\mathcal{C}_\Delta^{w_h}$ as a model of $C$ via the natural isomorphism $C\xrightarrow{\sim}C^{w_h}$, we get
\[C\smallsetminus\mathring{\mathcal{C}}_{\Delta,\eta}^{w_h}=\Spec\frac{K[x,y]}{\lb y^2-f(x), \prod_{o\neq h}\lb x-w_o\rb\rb}.\]
Thus the generic fibre of $\mathcal{C}$ is isomorphic to $C$.

\subsection{Special fibre}\label{SpecialFibreSubsection}
We now study the structure of the special fibre $\mathcal{C}_s$ of $\mathcal{C}$. As for the generic fibre, we consider \[\mathcal{C}_{\Delta,s}^{w_h}\smallsetminus \mathring{\mathcal{C}}_{\Delta,s}^{w_h},\]
for any $h=1,\dots,m$. For every choice of a proper cluster $\t\in\Sigma_C^{w_h}$, and $M$ associated to $\t$, let \[S_M:=\lb \mathcal{C}_{\Delta,s}^{w_h}\smallsetminus \mathring{\mathcal{C}}_{\Delta,s}^{w_h}\rb\cap X_M=\Spec \frac{O_K[X^{\pm 1},Y,Z]}{\lb\mathcal{F}_M^h(X,Y,Z),T_M^h(X,Y,Z),Y^{\tilde m_{23}}Z^{\tilde m_{33}},\pi\rb}.\]

\begin{lem}\label{TMSpecialFibreCaseLLemma}
Let $M=M_{L,i}$ for $L=L_\t^{w_h}$. Let $l\neq h$. If $\t=\s_l\wedge\s_h$, then $T_M^{hl}(X,Y,Z)=X^{-1}(X+u_{hl})$, otherwise
\begin{enumerate}[label=(\roman*)]
    \item $T_M^{hl}(X,Y,0)=1$ for $i=0,\dots,r_L$; \label{TML(i)}
    \item $T_M^{hl}(X,0,Z)=1$ for $i=0,\dots,r_L-1$. \label{TML(ii)}
\end{enumerate}
\proof
Fix $l\neq h$. If $\t\not\supseteq\s_l\wedge\s_h$, then by Lemma \ref{negativeexponentsLemma}, we have $\tilde m_{21}-\rho_{hl} \tilde m_{23}\geq 0$ and $\tilde m_{31}-\rho_{hl} \tilde m_{33}> 0$. Moreover, if $\tilde m_{21}-\rho_{hl} \tilde m_{23}=0$, then $i=r_L$. Therefore the equalities in \ref{TML(i)} and \ref{TML(ii)} follow directly from the definition of $T_M^{hl}$.

On the other hand, if $\t\supsetneq \s_l\wedge\s_h$, then by Lemma \ref{NonnegativeExponentsMatricesLemma}, we have $\rho_{hl} \tilde m_{23}-\tilde m_{21}\geq 0$ and $\rho_{hl} \tilde m_{33}-\tilde m_{31}> 0$. Moreover, if $\rho_{hl} \tilde m_{23}-\tilde m_{21}=0$, then $i=r_L$. Therefore we have \ref{TML(i)} and \ref{TML(ii)} again. 

Finally, assume $\t=\s_l\wedge\s_h$. Since $\rho_\t=\rho_{hl}\in\Z$, then
$\rho_{hl} \tilde m_{13} - \tilde m_{11}=-1$. Hence \[T_M^{hl}(X,Y,Z)=1+ u_{hl} X^{-1}=X^{-1}\lb X+u_{hl}\rb, \]
 by Lemma \ref{NonnegativeExponentsMatricesLemma}. 
\endproof
\end{lem}

\begin{lem}\label{SpecialFibreCaseLLemma}
Suppose $M=M_{L_\t^{w_h},i}$. Then
\[S_M=\Spec{\dfrac{O_K[X^{\pm 1},Y,Z]}{(\mathcal{F}_M^h(X,Y,Z),\prod_l\lb X+ u_{hl}\rb,Y^{\tilde m_{23}}Z^{\tilde m_{33}},\pi)}}\subset\mathcal{C}_\Delta^{w_h},\]
where the product runs over all $l\neq h$ such that $\t=\s_l\wedge\s_h$.
\proof
Lemma \ref{MatricesLemma} shows that $\tilde m_{33}$ is always different from $0$, while $\tilde m_{23}=0$ if and only if $i=r_{L_\t^{w_h}}$. Thus the result follows from Lemma \ref{TMSpecialFibreCaseLLemma}.
\endproof
\end{lem}

\begin{lem}\label{MultipleRootsfhLemma}
Let $f_h(x)=f(x+w_h)$ and $l\neq h$. Then $\ch{u_{lh}}$ is a multiple root of $\ch{f_h|_L}$ of order $|\t_l|$, where $L=L_{\s_h\wedge\s_l}^{w_h}$ and $\t_l\in\Sigma_C^{w_l}$, $\t_l<\s_h\wedge\s_l$.

Furthermore, if $\Sigma=\{\s_1,\dots,\s_m\}=\Sigma_C^\mathrm{min}$, $C$ has an almost rational cluster picture and $\bar\alpha\in \bar k$ is a multiple root of $\ch{f_h|_L}$ for some edge $L$ of $\NP{f_h}$, then $\bar\alpha=\ch{u_{lh}}$ and $L=L_{\s_h\wedge\s_l}^{w_h}$ for some $l\neq h$.
\proof
For any proper cluster $\s\in\Sigma_f$, let $\lambda_\s=\min_{r\in\s}v(r-w_h)$. 
Let $\s\in\Sigma_C^{w_l}$, with $\s_l\subseteq\s\subsetneq\s_h\wedge\s_l$. Then $w_h$ is not rational centre of $\s$, and for every root $r\in\s$, one has
\[v(r-w_h)=v(r-w_{l}+w_{l}-w_h)=\min\{v(r-w_l), \rho_{hl}\}=\rho_{hl},\]
as $v(r-w_l)\geq \rho_\s>\rho_{hl}$. Therefore $\lambda_\s=\rho_{hl}\in\Z$. In particular, $|\lambda_\s|_p\leq 1$. Furthermore, \[d_{\s}\geq\rho_{\s}>\lambda_{\s}=\rho_{hl}\qquad\text{and}\qquad \sfrac{r-w_h}{\pi^{\rho_{hl}}}\equiv \sfrac{w_{lh}}{\pi^{\rho_{hl}}}\mod\pi,\]
and so Theorem \ref{RegularityClustersTheorem}\ref{RegularityClustersThmParti} implies that $\ch{u_{lh}}=\tfrac{w_{lh}}{\pi^{\rho_{hl}}}\mod\pi$ is a multiple root of $\ch {f_h|_L}$, where $L=L_{\s_h\wedge\s_l}^{w_h}$.

Let $\t_l\in\Sigma_C^{w_l}$, $\t_l<\s_h\wedge\s_l$. Since $\s_l\subseteq\t_l<\s_h\wedge\s_l$ we have
\[\t_l=\left\{r\in\roots\mid\ch{u_{lh}}=\tfrac{r-w_h}{\pi^{\rho_{hl}}}\mod\pi\right\},\] 
as $v(r-w_l)>\rho_{hl}$ if and only if $\ch{u_{lh}}=\frac{r-w_h}{\pi^{\rho_{hl}}}\mod\pi$. Thus the multiplicity of $\ch{u_{lh}}$ is $|\t_l|$ by Theorem \ref{RegularityClustersTheorem}\ref{RegularityClustersThmPartii}.

Now let $\bar \alpha$ be a multiple root of $\ch {f_h|_L}$ for some edge $L$ of $\NP{f_h}$ and let $\s\in\Sigma_{f}$ associated to $\bar \alpha$ by Theorem \ref{RegularityClustersTheorem}\ref{RegularityClustersThmPartiii}. We want to prove that if $C$ has an almost rational cluster picture and $\Sigma=\Sigma_C^\mathrm{min}$, then there exists $l\neq h$ so that $\bar\alpha=\ch{u_{lh}}$. Note first $w_h$ is not a rational centre of $\s$. Indeed, if $w_h$ is a rational centre of $\s$, then 
\begin{equation*}
    |\s|>|\lambda_\s|_p=|\rho_\s|_p,\qquad d_\s>\lambda_\s=\rho_\s,
\end{equation*}
which contradicts the fact that $C$ has an almost rational cluster picture. As $\{\s_1,\dots,\s_m\}=\Sigma_C^\mathrm{min}$, we must have that $w_{l}$ is a rational centre of $\s$, for some $l\neq h$. Then $\s_l\subseteq\s\subsetneq\s_h\wedge\s_l$. Since $\bar\alpha=\tfrac{r-w_h}{\pi^{\lambda_\s}}\mod\pi$ for any $r\in\s$, from above we have $\bar\alpha=\ch{u_{lh}}$. Finally, $L$ is the edge of $\NP{f_h}$ of slope $-\lambda_\s=-\rho_{hl}$. Thus $L=L_{\s_h\wedge\s_l}^{w_h}$.
\endproof
\end{lem}

It remains to compute $S_M$ when $M=M_{V,j}$, where $V=V_\t^{w_h}$ or $V=V_0^{w_h}$.

\begin{lem}\label{TMSpecialFibreCaseVLemma}
Let $M=M_{V,j}$ for $V=V_\t^{w_h}$, or $V=V_0^{w_h}$ if $\t=\s_h$. For any $l\neq h$
we have
\begin{enumerate}[label=(\roman*)]
    \item $T_M^{hl}(X,Y,0)=1$ except when $\t=\s_l\wedge\s_h$ and $j=0$; \label{TMV(i)}
    \item $T_M^{hl}(X,0,Z)=1$ except when $\t<\s_l\wedge\s_h$ and $j=r_V$. \label{TMV(ii)}
\end{enumerate}
\proof
The lemma immediately follows from Lemmas \ref{NonnegativeExponentsMatricesLemma} and \ref{negativeexponentsLemma}.
\endproof
\end{lem}

\begin{lem}\label{SpecialFibreCaseVLemma}
Let $M=M_{V,j}$ with $V=V_\t^{w_h}$, or $V=V_0^{w_h}$ if $\t=\s_h$.
Then $S_M=\varnothing$.
\proof
For any $l\neq h$, we want to prove that
\begin{equation}\label{SpecialFibreEquation}
S_M^{hl}:=\{T_M^{hl}(X,Y,Z)=Y^{\tilde m_{23}}Z^{\tilde m_{33}}=0\}=\varnothing.
\end{equation} 
Lemma \ref{MatricesLemma} shows that $\tilde m_{33}$ is always different from $0$ and that $\tilde m_{23}=0$ if and only if $j=r_V$, and $V=V_\roots^{w_h}$ or $V=V_0^{w_h}$.
Assume that if $\t=\s_l\wedge \s_h$ then $j\neq 0$ and that if $\t<\s_l\wedge\s_h$ then $j\neq r_{V}$. Lemma \ref{TMSpecialFibreCaseVLemma} implies (\ref{SpecialFibreEquation}).

If $\t=\s_l\wedge \s_h$ and $j=0$, then $\rho_{hl} \tilde m_{33}-\tilde m_{31}= 0$ but $\rho_{hl} \tilde m_{23}-\tilde m_{21}> 0$. So
\[S_M^{hl}=\{T_M^{hl}(X,Y,Z)=Z^{\tilde m_{33}}=0\}\subset\Spec R[Y^{-1}].\]
Similarly, if $\t<\s_l\wedge \s_h$ and $j=r_{V}$, then $\tilde m_{21}-\rho_{hl} \tilde m_{23}= 0$, however $\tilde m_{31}-\rho_{hl} \tilde m_{33}> 0$. Then \[S_M^{hl}=\{T_M^{hl}(X,Y,Z)=Y^{\tilde m_{23}}=0\}\subset\Spec R[Z^{-1}].\]
In both cases, $S_M^{hl}\subseteq X_{F}$ as sets, where $F=F_{\s_l\wedge\s_h}^{w_h}$ (\cite[Definition 3.7]{Dok}). Let $L=L_{\s_l\wedge\s_h}^{w_h}$, and let $f_h(x)=f(x+w_h)$ and $g_h(x,y)=y^2-f_h(x)$. By Lemmas \ref{SpecialFibreCaseLLemma} and \ref{MultipleRootsfhLemma}, one has 
\[S_M^{hl}\subseteq X_{F}\cap S_{M_{L,0}}=\varnothing,\]
as $\mathcal{F}_{M_{L,0}}^h(X,Y,0)\mod \pi$ equals $Y^b-X^a\ch{f_h|_L}(X)$, for some $a\in\Z$, $b=1,2$ (see Lemma \ref{DefiningEquationXFLemma} for more details, whose proof is independent of this result).
Thus if $V=V_\t^{w_h}$ and $M=M_{V,j}$, then $S_M=\varnothing$.
%
\endproof
\end{lem}

\subsection{Components}\label{ComponentsSubsection}
Now that we have compared the special fibre of $\mathcal{C}$ with those of the models $\mathcal{C}_{\Delta}^{w_h}$, let us describe closed subschemes that form it. We will first study closed subschemes forming $\mathring{\mathcal{C}}_{\Delta,s}^{w_h}$ and then how they glue in $\mathcal{C}_s$.

Let $f_h(x)=f(x+w_h)$ and $g_h(x,y)=y^2-f_h(x)$. According to \cite[Theorem 3.14]{Dok} the special fibre of $\mathcal{C}_\Delta^{w_h}$ is formed by:
\begin{itemize}
    \item Chains of $\P_k^1$s coming from $v$-edges of $\Delta^{w_h}$.
    \item $1$-dimensional subschemes coming from $v$-faces of $\Delta^{w_h}$.
\end{itemize}
More precisely, each $v$-edge $E$ gives a scheme $X_E\times \P_E$, where $\P_E$ is a chain of $\P_k^1$s and $X_E\subset\G_{m,k}$ is given by $\ch{g_h|_E}=0$. The multiplicities and and the length of $\P_E$ can be completely described by the slopes of $E$. On the other hand, each $v$-face $F$ gives a proper scheme $\bar X_{F}$ containing an open subscheme $X_F\subseteq\G_{m,k}^2$ given by $\ch{g_h|_F}=0$. How the previous schemes intersect to form $\mathcal{C}_{\Delta,s}^{w_h}$ is described by \cite[Theorem 3.14]{Dok}. The reader is pointed to \cite{Dok} for more details. 

\begin{defn}\label{PolynomialsGeneralDefinition}
Let $\t\in\Sigma^W$ be a proper cluster. For any rational centre $w$ of $\t$, let $r_{\t,w}=\tfrac{w-r}{\pi^{\rho_\t}}$, $u_{\t,w}=c_f\tprod_{r\in\roots\setminus\t} r_{\t,w}$ and $u_{\s_h,w_h}^0=c_f\tprod_{r\in\roots\setminus\{w_h\}} r_{\s_h,w_h}$. Define $\ch{f_{\t,w}^W}, \ch{g_{\t,w}}\in k[X]$, and $\ch{g_{\s_h,w_h}^0}\in k[X]$ for any $h=1,\dots,m$, as follows:
\begin{enumerate}[label=(\arabic*)]
    \item Let $u=u_{\t,w}$. Define $\ch{f_{\t,w}^W}$ by \[\overline{f_{\t,w}^W}(X^{b_\t})= \tfrac{u}{\pi^{v(u)}}\prod_{r\in\t\setminus\bigcup_{\s<\t}\s} (X+r_{\t,w})\mod \pi,\] 
    where the union runs through all children $\s$ of $\t$ in $\Sigma^W$. 
    If $\Sigma=\Sigma_C^\mathrm{min}$ denote $\ch{f_{\t,w}^W}$ by $\ch{f_{\t,w}}$.
    \item Let $u=u_{\t,w}$. Define $\overline{g_{\t,w}}(X):=  X^{p_\t/\gamma_\t} - \tfrac{u}{\pi^{v(u)}}\mod \pi$.
    \item Let $u=u_{\s_h,w_h}^0$. Define $\overline{g_{\s_h,w_h}^0}(X):=X^{p_{\s_h}^0/\gamma_{\s_h}^0} - \tfrac{u}{\pi^{v(u)}}\mod \pi$.
\end{enumerate}
\end{defn}

Note that the polynomials defined in Definition \ref{PolynomialsGeneralDefinition} agree with the ones in Definition \ref{SchemesXsDefinition} when $w=w_\t$.

\begin{lem}\label{uwwLemma}
Let $\s,\t\in\Sigma_C^\mathrm{rat}$, with $\s\subsetneq\t$. Let $w',w$ be rational centres of $\s$ and $\t$ respectively, and define $\ch{u_{w'w}}=\tfrac{w'-w}{\pi^{\rho_\t}}\mod\pi$. Then $\ch{u_{w'w}}$ does not depend on the choice of a rational centre $w'$ of $\s$.
\proof
Suppose that $w_1,w_2$ are two rational centres of $\s$. Then $v(w_1-w_2)\geq\rho_\s>\rho_\t$, and so the lemma follows.
\endproof
\end{lem}

\begin{rem}\label{uwwRemark}
Let $\t\in\Sigma_C^{w_h}$. Let $l=1,\dots,m$, $l\neq h$. Then $\t=\s_h\wedge\s_l$ if and only if it has a child $\s\in\Sigma_C^{w_l}\setminus\Sigma_C^{w_h}$. In particular, if this happens, Lemma \ref{uwwLemma} shows that $\overline{u_{lh}}=\frac{w-w_h}{\pi^{\rho_\t}}\mod\pi$ for any rational centre $w$ of $\s$.
\end{rem}

\begin{defn}\label{hattGtDefinition}
Let $\t\in\Sigma_C^{w_h}$ be a proper cluster. Define $\hat\t^W:=\{\s\in\Sigma^W\cup\{\varnothing\}\mid \s<\t\}$, where $\varnothing<\t$ only if $\t$ has no child in $\Sigma^W$. If $\varnothing<\t$ then we will say that $w_h$ is the rational centre of $\varnothing$.

Define $\G_{\t,w_h}:=\G_{m,k}\setminus\bigcup_l\{\ch{u_{lh}}\}$, where the union runs through all $l\neq h$ such that $\s_l\wedge\s_h=\t$. Note that Remark \ref{uwwRemark} shows that $\G_{\t,w_h}=\A_k^1\setminus\bigcup_{\s\in \hat\t^W}\{\ch{u_{w_\s w_h}}\}$,
where $\ch{u_{w_\s w_h}}=\frac{w_\s-w_h}{\pi^{\rho_\t}}\mod\pi$, and $w_\s$ is any rational centre of $\s$.
\end{defn}

Let $\t\in\Sigma_C^{w_h}$ be a proper cluster. Let $V=V_\t^{w_h}$ and $M=M_{V,j}$. In \S\ref{SpecialFibreSubsection} we showed the special fibre of $U^h_M$ equals $X_M\cap\mathcal{C}_{\Delta,s}^{w_h}$.
Therefore the components of $\mathring{\mathcal{C}}_{\Delta,s}^{w_h}$ coming from $V$ are the same of those of $\mathcal{C}_{\Delta,s}^{w_h}$ given by the same $v$-edge. 
Therefore $V$ gives a closed subscheme $X_V\times\P_V$ of $\mathring{\mathcal{C}}_{\Delta,s}^{w_h}$, where $\P_V$ is a chain of $\P_k^1$s and $X_V:\{\overline{g_h|_{V}}=0\}$ over $\G_{m,k}$.
Lemma \ref{DescriptionNewtonPolytopeClustersLemma} implies that $\overline{g_h|_{V}}=\overline{g_{\t,w_h}}$. 

Let $V_0=V_0^{w_h}$ and $M=M_{V_0,j}$. Similarly to above, $X_M\cap \mathring{\mathcal{C}}_{\Delta,s}^{w_h}=X_M\cap\mathcal{C}_{\Delta,s}^{w_h}$ and so $V_0$ gives rise to a closed subscheme $X_{V_0}\times\P_{V_0}$ of $\mathring{\mathcal{C}}_{\Delta,s}^{w_h}$, where $\P_{V_0}$ is a chain of $\P_k^1$s and $X_{V_0}:\{\overline{g_h|_{V_0}}=0\}$ over $\G_{m,k}$. Note that $\overline{g_h|_{V_0}}=\overline{g_{\s_h,w_h}^0}$.

Let $\t\in\Sigma_C^{w_h}$ be a proper cluster. Let $L=L_\t^{w_h}$ and $M=M_{L,i}$. 
By Lemma \ref{SpecialFibreCaseLLemma}, the $v$-edge $L$ gives a subscheme $X_{L}^W\times\P_{L}$ of $\mathring{\mathcal{C}}_{\Delta,s}^{w_h}$, where $\P_{L}$ is a chain of $\P_k^1$s of length $r_L$ and $X_L^W:\{\overline{g_h|_{L}}=0\}$ in $\G_{\t,w_h}$. Note that $r_L=0$ or $1$ by Lemma \ref{DescriptionNewtonPolytopeClustersLemma} and $r_L=1$ if and only if $D_\t=1$.
Let $\t_h\in\Sigma_C^{w_h}$ be the unique child of $\t$ with rational centre $w_h$ or set $\t_h=\varnothing$ if $\t$ has no such child. 
We will show that 
\begin{equation}\label{ghftEquality}
\overline{g_h|_{L}}(X)=-\prod_{\s\in\hat\t^W,\,\s\neq\t_ h}(X+\ch{u_{w_\s w_h}})^{|\s|}\cdot \overline{f_{\t,w_h}^W}(X).
\end{equation}
where $\ch{u_{w_\s w_h}}=\frac{w_\s-w_h}{\pi^{\rho_\t}}\mod\pi$, and $w_\s$ is any rational centre of $\s$.

Suppose $\t\neq\s_h\wedge\s_l$ for any $l\neq h$. Equivalently, all children of $\t$ in $\Sigma^W$ (at most one) belong to $\Sigma_C^{w_h}$. 
Then Lemma \ref{DescriptionNewtonPolytopeClustersLemma} shows that $\overline{g_h|_{L}}=-\overline{f_{\t,w_h}^W}$.
Suppose now that $\t=\s_h\wedge\s_l$ for some $l\neq h$. In this case
$b_\t=1$. 
We have
\[\frac{\ch{g_h|_L}(X)}{\prod_{\s\in\hat\t^W,\s\neq\t_ h}(X+\ch{u_{w_\s w_h}})^{|\s|}}=\bigg(\!\frac{-\tfrac{u}{\pi^{v(u)}}\prod_{r\in\t\setminus\t_h}(X+r_{\t,w_h})}{\prod_{\s\in\hat\t^W,\s\neq\t_h}\prod_{r\in\s}(X+r_{\t,w_h})}\!\!\mod \pi\!\!\bigg)=-\overline{f_{\t,w_h}^W}(X),\]
where $r_{\t,w_h}$ and $u=u_{\t,w_h}$ are as in Definition \ref{PolynomialsGeneralDefinition}. 
Indeed, $\overline{u_{w_\s w_h}}= r_{\t,w_h}\mod \pi$ for every $r\in\s$ as $v(w_\s-r)\geq\rho_{\s}>\rho_\t$, and since $b_\t=1$, Lemma \ref{DescriptionNewtonPolytopeClustersLemma} implies that
 \[\ch{g_h|_L}(x)=-\tfrac{u}{\pi^{v(u)}}\tprod_{r\in\t\setminus\t_h}(x+r_{\t,w_h})\mod \pi.\]
In particular, Remark \ref{uwwLemma} and Lemma \ref{MultipleRootsfhLemma} shows that $(X+\ch{u_{hl}})\nmid \ch{f_{\t,w_h}^W}(X)$, for any $l\neq h$ such that $\s_l\wedge\s_h=\t$. Moreover, $X\nmid \ch{f_{\t,w_h}^W}(X)$ by definition. Therefore the scheme $X_L^W$ is equal to the closed subscheme $X_{\t,w_h}^W\subset\A_{k}^1$ given by $\ch{f_{\t,w_h}^W}=0$.

Let $\t\in\Sigma^W$ be a proper cluster. For any $h=1,\dots,m$ such that $\s_h\subseteq\t$, let $\bar X_{F_\t^{w_h}}$ be the $1$-dimensional closed subscheme of $\mathcal{C}_{\Delta,s}^{w_h}$ given by $F_\t^{w_h}$. Define \[\mathring{X}_{F_\t^{w_h}}:=\bar X_{F_\t^{w_h}}\cap\mathring{\mathcal{C}}_{\Delta}^{w_h}.\] Denote by $\Gamma_\t$ the $1$-dimensional closed subscheme of $\mathcal{C}_s$, result of the glueing of the subschemes $\mathring{X}_{F_\t^{w_h}}$ of $\mathring{\mathcal{C}}_{\Delta,s}^{w_h}$ to varying of $h$ such that $\t\in\Sigma_C^{w_h}$.

\begin{lem}\label{MultiplicityLemma}
Let $\t\in\Sigma_C^{w_h}$ be a proper cluster. The multiplicity of $\Gamma_\t$ in $\mathcal{C}_s$ is $m_\t$.
\proof
Let $L=L_\t^{w_h}$, $M=M_{L,0}$, and let $F=F_\t^{w_h}$. The multiplicity of $\bar X_{F_\t^{w_h}}$, and so of $\mathring X_{F_\t^{w_h}}$ and $\Gamma_\t$, is $\delta_F$. Hence we only need to show that $m_\t=\delta_F$. Let $d_0\in\Z$ as in Lemma \ref{MatricesLemma}. Then $\delta_F=\delta_L d_0$. The result follows as $\delta_L=b_\t$ and $d_0$, denominator of $s_1^L$, equals $3-D_\t$ by Lemma \ref{DescriptionNewtonPolytopeClustersLemma}.
\endproof
\end{lem}


\begin{lem}\label{DefiningEquationXFLemma}
Let $L=L_\t^{w_h}$, $F=F_\t^{w_h}$ and $M=M_{L,0}$. Let $c\in\{0,\dots,b_\t-1\}$ such that $1/b_\t-\rho_\t\cdot c\in\Z$. Then $\mathcal{F}_M^h(X,Y,0)\mod\pi$ equals the polynomial \[\ch{g_h|_F}(X,Y)=Y^{D_\t}-\prod_{\s\in \hat\t^W}(X-\ch{u_{w_\s w_h}})^{\frac{|\s|}{b_\t}-c\epsilon_\t}\ch{f_{\t,w_h}^W}(X),\] 
where $\ch{u_{w_\s w_h}}=\frac{w_\s-w_h}{\pi^{\rho_\t}}\mod\pi$, and $w_\s$ is any rational centre of $\s$.

In particular, $\Gamma_\t^h\subset\G_{\t,w_h}\times\A_k^1$ given by $\ch{g|_F}=0$ is the open subscheme $U_M^h\cap\{Z=0\}$ of $\mathring{X}_F$, and
the points in $S_M$ belong to all irreducible components of $\bar X_{F}$.
\proof
From \cite[\S3.5]{Dok} and the equation of $C^{w_h}$, the polynomial $\mathcal{F}_M^h(X,Y,0)$ reduces modulo $\pi$ to $X^{a_1}Y^b+X^{a_2}\ch{g_h|_{L}}(X)$, for some $b=1,2$ and $a\in\Z$.  Lemma \ref{DsInteriorPointsLemma} shows that $b=D_\t$. By Lemma \ref{DescriptionNewtonPolytopeClustersLemma}, $a_1=2\tilde m_{12}$, $a_2=|\t_h|\tilde m_{11}+(\epsilon_\t-|\t_h|\rho_\t)\tilde m_{13}$, where $\t_h\in\Sigma_C^{w_h}\cup\{\varnothing\}$, $\t_h<\t$. Then $a_1=0$ and $a_2=\tfrac{|\t_h|}{b_\t}-c\epsilon_\t$ by Lemma \ref{MatricesLemma}.

If $\t$ has one or no child, or $D_\t=1$, then $\ch{g_h|_L}=-\ch{f_{\t,w_h}^W}$ by (\ref{ghftEquality}). On the other hand, if $D_\t=2$ and $\t$ has two or more children in $\Sigma_C^\mathrm{rat}$, then $b_\t=1$, and so $c=0$. Therefore the equality (\ref{ghftEquality}) concludes the proof of the first part of the statement also in this case.
Finally, the last part of the lemma follows from Lemma \ref{SpecialFibreCaseLLemma}.
\endproof
\end{lem}


Let $c$ as in the previous lemma and define $\tilde\t^W:=\{\s\in\hat\t^W\mid \frac{|\s|}{b_\t}-c\epsilon_\t\notin 2\Z\}$.
\begin{prop}\label{DefiningEquationGammaGeneralCaseProposition}
Let $L=L_\t^{w_h}$ and $M=M_{L,0}$. The dense open subscheme $\Gamma_\t\cap U_M^h$ of $\Gamma_\t$ is isomorphic to the closed subscheme of $\G_{\t,w_h}\times\A^1_k$ given by
\[Y^{D_\t}=\prod_{\s\in\tilde\t^W}(X-\ch{u_{w_\s w_h}})\cdot\ch{f_{\t,w_h}^W}(X),\]
where $\ch{u_{w_\s w_h}}=\frac{w_\s-w_h}{\pi^{\rho_\t}}\mod\pi$, and $w_\s$ is any rational centre of $\s$.
\proof
The proposition follows from Lemma \ref{DefiningEquationXFLemma} and the definition of $\G_{\t,w_h}$.
\endproof
\end{prop}

We conclude this subsection describing how the glueing morphism (\ref{GlueingMapEquation}) restricts to the special fibre. Suppose $\t\supseteq\s_l\wedge\s_h$ for $l\neq h$ and let $M$ be a matrix associated to $\t$. Consider the glueing map $U_M^h\rightarrow U_M^l$ explicitly defined in \S\ref{GlueingSubsection}.

Suppose first $M=M_{V,j}$ with $V=V_\t^{w_l}$. By Lemma \ref{TMSpecialFibreCaseVLemma} the glueing morphism restricts to the identity on $X_V\times \P_V$.

Suppose $M=M_{L,i}$ with $L=L_\t^{w_l}$. Note that $\tilde m_{12}=0$ from Lemma \ref{MatricesLemma}. Recall the open subscheme $\Gamma_\t^h$ of $\mathring X_{F_\t^{w_h}}$ defined in Lemma \ref{DefiningEquationXFLemma}. Then, Lemma \ref{TMSpecialFibreCaseLLemma} implies that the glueing map restricts to an isomorphism $\Gamma_\t^h\mapsto\Gamma_\t^l$ induced by the ring homomorphism sending $X\mapsto X+\ch{u_{w_hw_l}}$, where $\ch{u_{w_hw_l}}=\tfrac{w_h-w_l}{\pi^{\rho_\t}}\mod\pi$. Similarly, it restricts to an isomorphism $X_{L_\t^{w_h}}^W\times \P_{L_\t^{w_h}}\rightarrow X_{L_\t^{w_l}}^W\times \P_{L_\t^{w_l}}$, where $\P_{L_\t^{w_h}}\rightarrow \P_{L_\t^{w_l}}$ is the identity and $X_{L_\t^{w_h}}^W\rightarrow X_{L_\t^{w_l}}^W$ is induced by the ring homomorphism sending $X\mapsto X+\ch{u_{w_hw_l}}$.

\subsection{Regularity}
Let $w_h\in W$. We want to show that if $\Sigma=\Sigma_C^\mathrm{min}$, and $C$ has an almost rational cluster picture and is $y$-regular, then $\mathring{\mathcal{C}}_{\Delta}^{w_h}$ is a regular scheme. 


\begin{lem}\label{SingularPointsModelGeneralCaseLemma}
Consider the model $\mathcal{C}_\Delta^{w_h}/O_K$ and let $f_h(x)=f(x+w_h)$. Suppose $\Sigma=\{\s_1,\dots,\s_m\}=\Sigma_C^\mathrm{min}$, and $C$ has an almost rational cluster picture and is $y$-regular. If $P$ is a singular point of $\mathcal{C}_\Delta^{w_h}$ then
\[P\in\Spec{\dfrac{O_K[X^{\pm 1},Y,Z]}{(\mathcal{F}_M^h(X,Y,Z), X+ u_{hl},Y^{\tilde m_{23}}Z^{\tilde m_{33}},\pi)}}\subset\mathcal{C}_\Delta^{w_h}\cap X_M,\]
for some $l\neq h$, where $M= M_{L_{\s_h\wedge\s_l}^{w_h},i}$ for $i=0,\dots,r_{L_{\s_h\wedge\s_l}^{w_h}}$. 
\proof
Denote by $m_\alpha(X)\in O_K[X]$ a lift of the minimal polynomial in $k[X]$ of $\bar\alpha\in k^\mathrm{s}$. By Lemma \ref{MultipleRootsfhLemma}, we only need to show that if
$P\in\mathcal{C}_\Delta^{w_h}$ is a singular point then 
\begin{equation}\label{SingularPointsEquation}
    P\in\Spec{\dfrac{O_K[X^{\pm 1},Y,Z]}{(\mathcal{F}_{M_{L,i}}^h(X,Y,Z),m_\alpha(X),Y^{\tilde m_{23}}Z^{\tilde m_{33}},\pi)}},
\end{equation}
for some $v$-edge $L=L_\t^{w_h}$ of $\Delta^{w_h}$, and some multiple root $\bar \alpha$ of $\ch{f_h|_L}$.
We study the polynomial $\mathcal{F}_M^h$ to varying of the matrix $M$, using \cite[\S 4.5]{Dok}. Let $g_h(x,y)=y^2-f_h(x)$. Let $L=L_\t^{w_h}$ and $M=M_{L,i}$. Note that $\ch{g_h|_L}=-\ch{f_h|_L}$. We have $\mathcal{F}_M^h(X,0,Z)=\ch{g_h|_L}(X)$ for any $i$. On the other hand, $\mathcal{F}_M^h(X,Y,0)=\ch{g_h|_L}(X)$ if $i>0$ and $\mathcal{F}_M^h(X,Y,0)=\ch{g_h|_F}(X,Y)$ if $i=0$. From the description given in Lemma \ref{DefiningEquationXFLemma}, we conclude that for these matrices $M$ the points in (\ref{SingularPointsEquation}) are the only possibly singular points of $\mathcal{C}_{\Delta}^{w_h}\cap X_M$. In particular, this proves that for any $v$-face $F$ of $\Delta^{w_h}$, the points in $X_F$ are non-singular in $\mathcal{C}_{\Delta}^{w_h}$.

Let $V=V_\t^{w_h}$ or $V=V_0^{w_h}$ and $M=M_{V,j}$. Since $C$ is $y$-regular, $p\nmid\deg(\ch{g_h|_V})$ by Lemma \ref{DsInteriorPointsLemma}. By \cite[\S 4.5]{Dok} and the fact that the points in $X_F$ are non-singular for all $v$-faces $F$, we conclude that $\mathcal{C}_\Delta^{w_h}$ has no singular point on $X_M$ for these matrices $M$, as required.
\endproof
\end{lem}

\begin{prop}\label{RegularityModelProposition}
Suppose $\Sigma=\Sigma_C^\mathrm{min}$, and $C$ has an almost rational cluster picture and is $y$-regular, then  $\mathcal{C}$ is a regular scheme. 
\proof
Lemmas \ref{SingularPointsModelGeneralCaseLemma} and \ref{SpecialFibreCaseLLemma} show that $\mathring{\mathcal{C}}_{\Delta}^{w_h}$ is regular for every $h$. Thus their glueing $\mathcal{C}$ is regular as well. 
\endproof
\end{prop}

\subsection{Separatedness}
It remains to prove that $\mathcal{C}$ is a proper scheme. We first show it is separated. Clearly it suffices to prove that $\mathcal{X}/O_K$ is separated.
Since the schemes $X_\Delta^h$ are separated, then the open subschemes $\mathring{X}_\Delta^{h}$ are separated as well by \cite[Proposition 3.3.9]{Liu}. 
Consider the open cover $\{V_M^h\}_{h,M}$ of $\mathcal{X}$. Let $h,l=1,\dots,m$ and let $M_h$ and $M_l$ be matrices associated to proper clusters $\t_h\in\Sigma_C^{w_h}$ and $\t_l\in\Sigma_C^{w_l}$ respectively. By \cite[Proposition 3.3.6]{Liu} 
we want to show
\begin{enumerate}[label=(\roman*)]
\item $V_{M_h}^h\cap V_{M_l}^l$ is affine,\label{i}
\item The canonical homomorphism
\[O_{\mathcal{X}}(V_{M_h}^h)\otimes_\Z O_{\mathcal{X}}(V_{M_l}^l)\longrightarrow O_{\mathcal{X}}(V_{M_h}^h\cap V_{M_l}^l)\]
is surjective.\label{ii}
\end{enumerate}
The definition of the glueing map (\ref{GlueingMapEquation}) implies \ref{i}. If $h=l$, or $\s_l\subseteq\t_h$, or $\s_h\subseteq\t_l$, then \ref{ii} follows from the separatedness of $\mathring{X}_\Delta^{h}$ and $\mathring{X}_\Delta^{l}$. So assume $l\neq h$, and $\t_h,\t_l\subsetneq\s_h\wedge\s_l$. 
Consider the Moebius transformation
\[\psi_l:\quad x\mapsto \frac{x}{xw_{hl}^{-1}+1},\quad y\mapsto\frac{y}{(xw_{hl}^{-1}+1)^{g+1}}.\]
It sends the curve $C^{w_l}$ to the isomorphic hyperelliptic curve \[C_l^h:y^2=(xw_{hl}^{-1}+1)^{2g+2}f\lb x(xw_{hl}^{-1}+1)^{-1}+w_l\rb.\]
As \begin{align*}
    f_l^h(x):&=(xw_{hl}^{-1}+1)^{2g+2}f\lb x(xw_{hl}^{-1}+1)^{-1}+w_l\rb\\&=c_fw_{hl}^{|\roots|}(xw_{hl}^{-1}+1)^{2g+2-|\roots|}\prod_{r\in\roots\smallsetminus\{w_h\}}\frac{r-w_h}{w_{lh}}\lb xw_{hl}^{-1}+\frac{r-w_l}{r-w_h}\rb,
\end{align*}
every cluster $\s\in\Sigma_C^{w_l}$ such that $\s\subsetneq\s_h\wedge\s_l$, corresponds to a unique cluster $\s^h\in\Sigma_{C_l^h}^0$ of same size, same radius and rational centre $0$. Moreover,
\[\epsilon_{\s^h}=v(c_{f_l^h})+\sum_{r'\in\s^h}\rho_{\s^h}+\sum_{r'\notin\s^h}v(r')=\epsilon_{\s}.\]
Call $\t_l^h$ the cluster in $\Sigma_{C_l^h}^0$ corresponding to $\t_l$.
Let $\Delta^{lh}$ and $\Delta_v^{lh}$ be the Newton polytopes attached to $y^2-f_l^h(x)$ and let $X_\Delta^{lh}$ be the associated toric scheme (defined in \cite[\S4.2]{Dok}). Since $\t_l\subsetneq\s_h\wedge\s_l$, the $v$-faces $F_{\t_l}$ of $\Delta^{w_l}$ and $F_{\t_l^h}$ of $\Delta^{lh}$ are identical by Lemma \ref{DescriptionNewtonPolytopeClustersLemma}. Furthermore, note that if $\t_l<\s_h\wedge\s_l$, then $\rho_{P(\t_l^h)}\leq\rho_{hl}=\rho_{P(\t_l)}$ and so $s_2^{V^0}\leq s_2^{V}$, where $V^0=V_{\t_l^h}^0$ and $V=V_{\t_l}^{w_l}$. Therefore the matrix $M:=M_l$ is also associated to $\t_l^h$. 

For every $o=1,\dots,m$, with $o\neq l$, define \[w_{hlo}=\begin{cases}\frac{w_{hl}w_{lo}}{w_{ho}}&\mbox{if } o\neq h,\\
w_{hl}&\mbox{if }o=h,\end{cases}\]
and write $w_{hlo}=u_{hlo}\pi^{\rho_{hlo}}$, where $u_{hlo}\in O_K^\times$ and $\rho_{hlo}\in\Z$, i.e.\ \[u_{hlo}=\begin{cases}\frac{u_{hl}u_{lo}}{u_{ho}}&\mbox{if } o\neq h,\\u_{hl}&\mbox{if }o=h,\end{cases}\quad\mbox{ and }\quad\rho_{hlo}=\begin{cases}\rho_{hl}+\rho_{lo}-\rho_{ho}&\mbox{if } o\neq h,\\\rho_{hl}&\mbox{if }o=h.\end{cases}\] 
Define 
\[\tilde T_M^{hlo}(X,Y,Z):=\begin{cases}1+u_{hlo}X^{\rho_{hlo}\tilde m_{13}-\tilde m_{11}}Y^{\rho_{hlo}\tilde m_{23}-\tilde m_{21}}Z^{\rho_{hlo}\tilde m_{33}-\tilde m_{31}}&\mbox{if }\t_l\supseteq\s_o,\\
u_{hlo}^{-1}X^{\tilde m_{11}-\rho_{hlo}\tilde m_{13}}Y^{\tilde m_{21}-\rho_{hlo}\tilde m_{23}}Z^{\tilde m_{31}-\rho_{hlo}\tilde m_{33}}+1&\mbox{if }\t_l\not\supseteq\s_o. \end{cases}\]
We want to show $\tilde T_M^{hlo}(X,Y,Z)\in R$. If $o=h$ then \[\tilde T_M^{hlo}(X,Y,Z)=T_M^{hl}(X,Y,Z)\in R.\] So assume $o\neq h$. If $\s_o\subseteq\t_l$, then it follows from Lemma \ref{NonnegativeExponentsMatricesLemma} as $\s_l\wedge\s_o\subsetneq\s_l\wedge\s_h$ and so $\rho_{hlo}=\rho_{lo}$. On the other hand, if $\s_o\not \subseteq\t_l$, then it follows from Lemma \ref{negativeexponentsLemma} as $\tilde m_{23},\tilde m_{33}>0$ and $\rho_{hlo}\leq\max\{\rho_{hl},\rho_{lo}\}$.
Let \[\tilde T_M^{hl}(X,Y,Z):=\prod_{o\neq l}\tilde T_M^{hlo}(X,Y,Z).\]
The Moebius transformation
\[K[x^{\pm 1}, y^{\pm 1},\tprod_{o\neq l}(x+w_{lo})^{-1}]\stackrel{\psi_l}{\longrightarrow} K[x^{\pm 1},y^{\pm 1},\tprod_{o\neq l}\lb x+w_{hlo}\rb^{-1}]\]
considered above induces an isomorphism
\[R[T_M^l(X,Y,Z)^{-1}]\xrightarrow{M^{-1}\circ\psi_l\circ M} R[\tilde T_M^{hl}(X,Y,Z)^{-1}],\]
sending
\begin{align*}
X&\mapsto X\cdot T_M^{hl}(X,Y,Z)^{-m_{11}-(g+1) m_{21}},\\ 
Y&\mapsto Y\cdot T_M^{hl}(X,Y,Z)^{-m_{12}-(g+1) m_{22}},\\
Z&\mapsto Z\cdot T_M^{hl}(X,Y,Z)^{-m_{13}-(g+1) m_{23}}.
\end{align*}
Then
\[\tilde V_M^{lh}:=\Spec R[\tilde T_M^{hl}(X,Y,Z)^{-1}]\]
is an open subscheme of $X_\Delta^{lh}$, isomorphic to $V_M^l$. We can clearly carry out similar constructions for $t_h$, $M_h$.

By comparing the Newton polytopes $\Delta_v^{lh}$ and $\Delta_v^{hl}$, we see that the Moebius transformation $x\mapsto w_{hl}/(w_{lh}^{-1}x)$, $y\mapsto y/(w_{lh}^{-1}x)^{g+1}$ gives an isomorphism
\[
   \psi: K[x^{\pm 1},y^{\pm 1}, \prod_{o\neq l}(x+w_{hlo})^{-1}]\longrightarrow K[x^{\pm 1},y^{\pm 1}, \prod_{o\neq h}(x+w_{lho})^{-1}]
\]
which induces a birational map $X_\Delta^{hl}\-->X_\Delta^{lh}$, defined on the open set $\tilde V_{M_h}^{hl}$ of $X_\Delta^{hl}$. In particular, there exists an open set $\tilde V_{M_h}^{lh}$ of $X_\Delta^{lh}$, isomorphic to $V_{M_h}^h$ via the map induced by $\psi_h^{-1}\circ\psi$. 

Recall the definition of $\phi$ in (\ref{GlueingMap(-1)equation}), which induces the glueing map between $V_{M_l}^l$ and $V_{M_h}^h$. Since the following diagram
\[\begin{tikzcd}
K[x^{\pm 1}, y^{\pm 1},\prod_{o\neq l}(x+w_{lo})^{-1}]\arrow[r,"\phi"]\arrow[d,"\psi_l"]&K[x^{\pm 1}, y^{\pm 1},\prod_{o\neq h}(x+w_{ho})^{-1}]\arrow[d, "\psi_h"]\\
K[x^{\pm 1},y^{\pm 1}, \prod_{o\neq l}(x+w_{hlo})^{-1}]\arrow[r, "\psi"] &K[x^{\pm 1},y^{\pm 1}, \prod_{o\neq h}(x+w_{lho})^{-1}]
\end{tikzcd}\]
is commutative, then the surjectivity of 
\[O_{\mathcal{X}}(V_{M_h}^h)\otimes_\Z O_{\mathcal{X}}(V_{M_l}^l)\longrightarrow O_{\mathcal{X}}(V_{M_h}^h\cap V_{M_l}^l)\]
follows from the separatedness of $X_\Delta^{lh}$.

\subsection{Properness}
By \cite[IV.15.7.10]{EGA}, it remains to show that $\mathcal{C}_s$ is proper. From \cite[Exercise 3.3.11]{Liu}, we only need to prove that the $1$-dimensional subscheme $\Gamma_\t$ is proper for every $\t=\s_h\wedge\s_l$. Indeed every other component is entirely contained in a model $\mathcal{C}_\Delta^{w_h}$, which is proper (see \S\ref{SpecialFibreSubsection}). Let $\t=\s_h\wedge\s_l$ for some $h,l=1,\dots,m$, with $h\neq l$. For any $o=1,\dots, m$ such that $\s_o\subset\t$, let $\t_o$ be the unique child of $\t$ with $\s_o\subseteq\t_o<\t$. Then $\Gamma_\t$ is equal to the glueing of the schemes
\[\Spec\frac{R[T_M^{o}(X,Y,Z)^{-1}]}{\lb\mathcal{F}_M^o(X,Y,Z),Z,\pi\rb},\quad M=M_{L_\t^{w_o},0},M_{V_\t^{w_o},0},\]
and
\[\Spec\frac{R[T_M^{o}(X,Y,Z)^{-1}]}{\lb\mathcal{F}_M^o(X,Y,Z),Y,\pi\rb},\quad M=M_{V_{\t_o}^{w_o},r_{V_{\t_o}^{w_o}}},\]
for all $o$ such that $\s_o\subset\t$, through the isomorphism (\ref{GlueingMapEquation}) and the glueing maps in the definition of $\mathcal{C}_\Delta^{w_o}$.
In particular, for any $o$ as above there exists a natural birational map $s_o:\Gamma_\t\-->\bar X_{F_\t^{w_o}}$ which is defined as the identity morphism on the dense open $\mathring X_{F_\t^{w_o}}=\Gamma_\t\cap\mathring{\mathcal{C}}_\Delta^{w_o}$. 

Let $D/k$ be a normal curve, let $P\in D$ and let $D\smallsetminus\{P\}\stackrel{g}{\longrightarrow} \Gamma_\t$ be a non-constant morphism of curves. We want to show that $g$ extends to $D$. For every $o$ as above, $\bar X_{F_\t^{w_o}}$ is proper, so the birational map 
\[g_o:=s_o\circ g:D\smallsetminus\{P\}\--> \bar X_{F_\t^{w_o}}\]
extends to a morphism $\bar g_o: D\longrightarrow \bar X_{F_\t^{w_o}}$. 
If \[P_o:=\bar g_o(P)\in\lb\bar X_{F_\t^{w_o}}\cap\mathring{\mathcal{C}}_\Delta^{w_o}\rb= s_o\lb\Gamma_\t\cap\mathring{\mathcal{C}}_\Delta^{w_o}\rb\] for some $o$ such that $\s_o\subset\t$ (we will later show this is always the case), then there exists an open neighbourhood $U$ of $P_o$ such that $U\subseteq\lb\bar X_{F_\t^{w_o}}\cap\mathring{\mathcal{C}}_\Delta^{w_o}\rb$ and so $s_o|_{s_o^{-1}(U)}^U$ is an isomorphism. Since $P\in \bar g_o^{-1}(U)$, the map
\[\bar g_o^{-1}(U)\xrightarrow{\bar g_o|_{\bar g_o^{-1}(U)}^U}U\xrightarrow{\big( s_o|_{s_o^{-1}(U)}^U\big)^{-1}}s_o^{-1}(U)\hookrightarrow\Gamma_\t,\]
induces an extension $D\longrightarrow \Gamma_\t$ of $g$. 

Suppose that $P_o\notin\bar X_{F_\t^{w_o}}\cap\mathring{\mathcal{C}}_\Delta^{w_o}$ for any $o$ such that $\s_o\subset\t$. From \S \ref{SpecialFibreSubsection} we have
\begin{equation}\label{PropernessSingularPointEquation}
    P_o\in S_M=\Spec{\dfrac{R}{(\mathcal{F}_M^o(X,Y,Z),\prod_l\lb X+ u_{ol}\rb,Z,\pi)}},
\end{equation}
where $M=M_{L_\t^{w_o},0}$, and the product runs over all $l\neq o$ such that $\t=\s_o\wedge\s_l$. In particular $P_o$ is a point of each irreducible component of $\bar X_{F_\t^{w_o}}$ by Lemma \ref{DefiningEquationXFLemma}.
Let $h\neq o$ such that $X+u_{oh}$ vanishes at $P_o$. Let $\xi$ be the generic point of $D$ and let $\xi_o=g_o(\xi)$, $\xi_h=g_h(\xi)$ be generic points of $\bar X_{F_\t^{w_o}}$ and $\bar X_{F_\t^{w_h}}$ respectively. Then the birational maps $s_o$ and $s_h$ give
\[\begin{tikzcd}[row sep=10pt]
&&\bar X_{F_\t^{w_o}}\\
D\smallsetminus\{P\}\arrow[r, "g"]&\Gamma_\t\arrow[ur, dashrightarrow, "s_o"]\arrow[dr, dashrightarrow, "s_h"]&\\
&&\bar X_{F_\t^{w_h}}
\end{tikzcd}
\quad\Longrightarrow\quad
\begin{tikzcd}[row sep=10pt]
&&
k\lb\ch{\xi_o}\rb
\arrow[dll, "\phi_{g_o}"',bend left=10]\arrow[dd, "\simeq", bend right=20]\\
k(D)&&\\
&&
k\lb\ch{\xi_h}\rb
\arrow[ull, "\phi_{g_h}",bend right=10]
\end{tikzcd}\]
where we denote by $\phi_{g_o}$ and $\phi_{g_h}$ the homomorphisms between function fields induced by $g_o$ and $g_h$.
The vertical isomorphism is induced by the map
\[\frac{R[T_M^o(X,Y,Z)^{-1}]}{\lb \mathcal{F}_M^o(X,Y,Z),Z\rb}\longrightarrow\frac{R[T_M^h(X,Y,Z)^{-1}]}{\lb \mathcal{F}_M^h(X,Y,Z),Z\rb}\]
which sends (see \S \ref{GlueingSubsection} and Lemma \ref{TMSpecialFibreCaseLLemma})
\[X+u_{oh}\mapsto X\cdot T_M^{ho}(X,Y,Z)^{m_{11}}+u_{oh}=X\lb1+u_{ho}X^{-1}\rb+u_{oh}=X.\]
But the rational function $X+u_{oh}$ vanishes at $P_o$, while $X$ does not vanish at $P_h$ by (\ref{PropernessSingularPointEquation}). This gives a contradiction, as $\bar g_o(P)=P_o$ and $\bar g_h(P)=P_h$.

\subsection{Genus}\label{GenusSubsection}
Suppose $\Sigma=\{\s_1,\dots,\s_m\}=\Sigma_C^\mathrm{min}$, and $C$ has an almost rational cluster picture and is $y$-regular. In the previous subsections we proved that $\mathcal{C}/O_K$ is a proper regular model of $C$. Let $\t\in\Sigma_C^{w_h}$ be a proper cluster. 

\begin{prop}\label{DefiningEquationGammaProposition}
Let $\t\in\Sigma_C^{w_h}$. Then $\Gamma_\t$ is isomorphic to the smooth projective $1$-dimensional scheme given by 
\[Y^{D_\t}=\prod_{\s\in\tilde\t^W}(X-\ch{u_{w_\s w_h}})\ch{f_{\t,w_h}}(X)\]
where $\ch{u_{w_\s w_h}}=\frac{w_\s-w_h}{\pi^{\rho_\t}}\mod\pi$, and $w_\s$ is any rational centre of $\s$.

In particular, 
\begin{enumerate}
    \item if $D_\t=1$, then $\Gamma_\t\simeq\P^1_k$;
    \item if $D_\t=2$ and $\t$ is \"{u}bereven, then $\Gamma_\t$ is the disjoint union of two $\P^1$s over some quadratic extension of $k$;
    \item in all other cases, $\Gamma_\t$ is a hyperelliptic curve of genus $g(\t)$.
\end{enumerate}
\proof
The first part of the proposition follows from Proposition \ref{DefiningEquationGammaGeneralCaseProposition}.

For the second part of the statement note that if $D_\t=1$ then the result follows. Suppose $D_\t=2$. Then $p\neq 2$ as $C$ is $y$-regular. Note that since $\Sigma=\Sigma_C^\mathrm{min}$, the proper clusters in $\Sigma^W$ correspond to the proper clusters in $\Sigma_C^\mathrm{rat}$. Recall the definition of $\tilde\t$ given in Definition \ref{GenussDefinition}. Let $h(X)=\prod_{\s\in\tilde\t^W}(X-\ch{u_{w_\s w_h}})\ch{f_{\t,w_h}}(X)$.

Suppose $\t$ is \"{u}bereven. Then all its children are (proper) rational cluster by Lemma \ref{UberevenLemma} since they are even and $p\neq 2$. In particular $b_\t=1$ by Lemma \ref{TwoChildIntegralRadiusLemma} and so $\epsilon_\t\in 2\Z$ and $\tilde \t=\tilde\t^W=\varnothing$ since it equals the set of odd rational children. 
Moreover, $\t=\bigcup_{\s<\t,\,\s\,\text{proper}}\s$, and so $\ch{f_{\t,w_h}}\in k$. Thus $h(X)\in k$.

Now suppose $h(X)\in k$. Then $\tilde\t^W=\varnothing$ and $\t=\bigcup_{\s<\t}\s$, where $\s$ runs through all children $\s\in\Sigma^W$ of $\t$. The non-proper clusters in $\Sigma^W$ are of the form $\{w_l\}$ for some $l=1,\dots,m$. If $\{w_l\}<\t$, then $\t=\s_l$, but in that case $\t$ would not equal the union of its children in $\Sigma^W$. Hence $\t$ has no non-proper children. It follows that $\tilde\t=\tilde\t^W$ and $\t$ equals the union of its proper rational children. In particular, $\t$ has two or more children in $\Sigma_C^\mathrm{rat}$, so $b_\t=1$, by Lemma \ref{TwoChildIntegralRadiusLemma}. But then $\tilde\t$ is the set of odd children of $\t$ as $\epsilon_\t\in 2\Z$, and so all rational children of $\t$ are even.

It only remains to prove that if $h(x)\notin k$, then the genus of $\Gamma_\t$ is $g(\t)$. Since $h(X)$ is a separable polynomial, we need to show that \[\deg h=\frac{|\t|-\sum_{\s\in\Sigma_C^\mathrm{rat},\,\s<\t}|\s|}{b_\t}+\tilde\t.\] It suffices to prove that if $\s\in\Sigma_C^\mathrm{rat}$ is a non-proper rational child of $\t$ different from $\{w_h\}$, then $b_\t=1$ and $\s\in\tilde\t$. Suppose $\s=\{r\}$ is such a rational cluster. Since $r\in\t$, we have $v(r-w_h)\geq \rho_\t$. Suppose $v(r-w_h)> \rho_\t$. Then $\s\in\mathring\Sigma_C^{w_h}$, as $\s<\t$ and $r\neq w_h$. But this contradicts our choice of $W$. Then $\rho_\t=v(r-w_h)\in\Z$ and so $b_\t=1$. It follows that $\tilde\t$ is the set of odd children of $\t$. Thus $\s\in\tilde \t$.
\endproof
\end{prop}

\subsection{Minimal regular NC model} 
Suppose the base extended curve $C_{K^{nr}}$ is $y$-regular and has an almost rational cluster picture. Consider the model $\mathcal{C}/O_{K^{nr}}$ constructed before with $\Sigma=\Sigma_{C_{K^{nr}}}^\mathrm{min}$. We want to see what components of $\mathcal{C}_s$ should be blown down to obtain the minimal regular model with normal crossings. Recall \cite[\S 5]{Dok}. Let $\Sigma_{K^{nr}}=\Sigma_{C_{K^{nr}}}^\mathrm{rat}$ and fix a proper cluster $\t\in\Sigma_{C_{K^{nr}}}^{w_h}$. 

Suppose first $\t\neq\s_h\wedge\s_l$ for all $l=1,\dots,m$ with $l\neq h$. Equivalently, $\t$ has at most one proper child in $\Sigma_{K^{nr}}$. Then $\Gamma_\t\simeq\bar X_{F_\t^{w_h}}$ and can be seen entirely in $\mathring{\mathcal{C}}_{\Delta}^{w_h}$. In particular, if $\Gamma_\t$ can be blown down then $F_\t^{w_h}$ is a removable or contractible $v$-face (see \cite[Theorem 5.7]{Dok}). By Lemma \ref{DescriptionNewtonPolytopeClustersLemma}, we find
\begin{itemize}
    \item $F_\t^{w_h}$ is removable if and only if $\t=\roots$ with a child in $\Sigma_{K^{nr}}$ of size $2g+1$.
    \item $F_\t^{w_h}$ is contractible if and only if either $|\t|=2$ and $\frac{\epsilon_\t}{2}-\rho_\t\in\Z$  or $\t$ has a proper rational child $\s\in\Sigma_{K^{nr}}$, of size $2g$, and $\frac{\epsilon_{\t}}{2}-g\rho_{\t}\in\Z$.
\end{itemize}
Recall Definition \ref{RemovableContractibleDefinition}. Note that $F_\t^{w_h}$ is removable if and only if $\t$ is removable. In this case, $F_\t^{w_h}$ can be ignored for the construction of $\mathcal{C}_\Delta^{w_h}$ (for any $h$ since $\t=\roots$), and so $\t$ can be ignored for the construction of $\mathcal{C}$. 

Assume now $F_\t^{w_h}$ contractible. We want to understand when $\Gamma_\t$ can be blown down. First consider the case $|\t|=2$ and $\frac{\epsilon_\t}{2}-\rho_\t\in\Z$. Then $\Gamma_\t$ intersects other components of $\mathcal{C}_s$ in $2$ points (as $V_\t^{w_h}$ gives two chains of $\P^1$s and the $v$-edges $V_0^{w_h}$ and $L_\t^{w_h}$ give no component in $\mathcal{C}_{\Delta,s}^{w_h}$). To have self-intersection $-1$, $\Gamma_\t$ has to have multiplicity $>1$. It follows from Lemma \ref{MultiplicityLemma} that $\rho_\t\notin\Z$, as $\frac{\epsilon_\t}{2}-\rho_\t\in\Z$. Moreover, by Lemma \ref{RadiusTimesSizeIntegerLemma}, one has $\rho_\t\in\frac{1}{2}\Z$. Therefore $\epsilon_\t$ is odd and the multiplicity of $\Gamma_\t$ is $2$. Let $r:=r_{V_\t^{w_h}}$ and consider
\[\gamma_\t s_\t=\sfrac{n_0}{d_0}>\sfrac{n_1}{d_1}>\dots>\sfrac{n_r}{d_r}>\sfrac{n_{r+1}}{d_{r+1}}=\gamma_\t\lb s_\t-\rho_\t+\rho_{P(\t)}\rb\]
given by $V_\t^{w_h}$. If $\Gamma_\t$ can be blown down then $d_1=1$. Since $\gamma_\t s_\t=-\frac{\epsilon_\t}{2}+2\rho_\t$, we have $d_0=2$. In particular $d_1=1$ if and only if $\rho_\t-\rho_{P(\t)}=\frac{n_0}{d_0}-\frac{n_{r+1}}{d_{r+1}}\geq\frac{1}{2}$ (see also \cite[Remark 3.15]{Dok}). Thus if $|\t|=2$, then $\Gamma_\t$ can be blown down if and only if $\rho_\t\notin\Z$, $\epsilon_\t$ odd, $\rho_{P(\t)}\leq\rho_\t-\frac{1}{2}$. Note that this is case (\ref{1Contr}) of Definition \ref{RemovableContractibleDefinition}.

Second consider the case $|\t|=2g+2$ with a proper rational child $\s$ of size $2g$ and $\tfrac{\epsilon_\t}{2}-g\rho_\t\in\Z$. The argument is very similar to the previous one. If $\Gamma_\t$ can be blown down then it must have multiplicity $>1$ and this implies $\rho_\t\notin\Z$ again by Lemma \ref{MultiplicityLemma}. From Lemma \ref{RadiusTimesSizeIntegerLemma} it follows that $(|\t|-|\s|)\rho_{\t}\in\Z$, so $\rho_\t\in\frac{1}{2}\Z$. Then $m_\t=2$ and 
\[\sfrac{v(c_f)}{2}=\sfrac{\epsilon_{\t}}{2}-(g+1)\rho_{\t}\in\tfrac{1}{2}\Z\setminus\Z,\]
so $v(c_f)$ odd. Let $r:=r_{V_{\s}^{w_h}}$ and consider
\[\gamma_{\s}s_{\s}=\sfrac{n_0}{d_0}>\sfrac{n_1}{d_1}>\dots>\sfrac{n_r}{d_r}>\sfrac{n_{r+1}}{d_{r+1}}=\gamma_{\s}(s_{\s}-\rho_{\s}+\rho_\t)\]
given by $V_{\s}^{w_h}$. If $\Gamma_\t$ can be blown down then $d_r=1$. Recall that $\epsilon_{\s}-|\s|\rho_{\s}=\epsilon_\t-|\s|\rho_\t$. Then $\gamma_{\s}(s_{\s}-\rho_{\s}+\rho_\t)=-\frac{\epsilon_\t}{2}+(g+1)\rho_\t$, so $d_{r+1}=2$. In particular $d_r=1$ if and only if $\rho_{\s}-\rho_\t=\frac{n_0}{d_0}-\frac{n_{r+1}}{d_{r+1}}\geq\frac{1}{2}$.
Thus if $\t$ has size $2g+2$ and has a unique proper rational child $\s\in\Sigma_{K^{nr}}$, then $\Gamma_\t$ can be blown down if and only if $|\s|=2g$, $\rho_\t\notin\Z$, $v(c_f)$ odd, $\rho_{\s}\geq\rho_\t+\frac{1}{2}$. This is case (\ref{2Contr}) of Definition \ref{RemovableContractibleDefinition}. 

Finally, if $|\t|=2g+1$, $\t$ has a proper child $\s\in\Sigma_{K^{nr}}$ of size $2g$ and $\frac{\epsilon_\t}{2}-g\rho_\t\in\Z$, then $\rho_\t\in\Z$, as $(|\t|-|\s|)\rho_\t\in\Z$. It follows that $\epsilon_\t\in\Z$ and so $m_\t=1$. This implies the self-intersection of $\Gamma_\t$ is not $-1$, since it intersects the rest of $\mathcal{C}_\t$ in at least two points as before. Hence in this case $\Gamma_\t$ can never be blown down.

Now assume there exists $l\neq h$ such that $\t=\s_h\wedge\s_l$. Then $\t$ is not minimal. Let $\t_h,\t_l\in\Sigma_{K^{nr}}$ be such that $\s_h\subseteq\t_h<\t$ and $\s_l\subseteq\t_l<\t$. Suppose $\Gamma_\t$ irreducible. If $|\t|\leq 2g$ (or, equivalently, $\t$ is not the largest non-removable cluster), then $\Gamma_\t$ intersects at least other $3$ components of $\mathcal{C}_s$ (given by $\t_h,\t_l$, and $P(\t)$). So it cannot be contracted to obtain a model with normal crossings. 
A similar argument holds if there exists $o\neq l$ such that $\s_o\wedge\s_h=\t$: at least $3$ components (given by $\t_h$, $\t_l$ and $\t_o$) intersect $\Gamma_\t$, so blowing down $\Gamma_\t$ would make the model lose normal crossings. Assume then $|\t|>2g$ and $\s_o\wedge\s_h\neq\t$ for all $o\neq l$. Then $\Gamma_\t$ intersects at least other $2$ components of $\mathcal{C}_s$ given by $V_{\t_h}^{w_h}$ and $V_{\t_l}^{w_l}$. Firstly, if $\Gamma_\t$ can be blown down, then $m_\t>1$. But $\rho_\t=\rho_{hl}\in\Z$. Then $m_\t$ is at most $2$. If $m_\t=2$ then $D_\t=1$, that implies $\epsilon_\t$ odd
and $\Gamma_\t\simeq\P^1$ by Proposition \ref{DefiningEquationGammaProposition}. It also follows $s_\t\in\tfrac{1}{2}\Z\setminus\Z$. If $\t$ is odd then this implies that $V_\t^{w_h}$ gives a $\P^1$ intersecting $\Gamma_\t$. Since that would be a third component intersecting $\Gamma_\t$, the cluster $\t$ has to be even. 
Hence $\t=\roots$ and $|\t|=2g+2$. Then $\epsilon_\t$ is odd if and only if $v(c_f)$ is odd, as $\rho_\t\in\Z$. Now, $L_\t^{w_h}$ gives some $\P^1$s intersecting $\bar X_{F_\t^{w_h}}\subset\mathcal{C}_{\Delta,s}^{w_h}$. 
All these $\P^1$s are not in $\mathring{\mathcal{C}}_{\Delta,s}^{w_h}$ (and so in $\mathcal{C}_s$) if and only if $\t_h\cup\t_l=\t$. 
In particular, $\t_h$ and $\t_l$ are either both even or both odd.
If $\t_h$ is even, then $\gamma_{\t_h}=2$, and so the component given by $V_{\t_h}^{w_h}$ has multiplicity at least $2$. 
The self-intersection of $\Gamma_\t$ could not be $-1$ in this case. Assume $\t_h$ is odd. Let $r:=r_{V_{\t_h}^{w_h}}$ and consider
\[\gamma_{\t_h}s_{\t_h}=\sfrac{n_0}{d_0}>\sfrac{n_1}{d_1}>\dots>\sfrac{n_r}{d_r}>\sfrac{n_{r+1}}{d_{r+1}}=\gamma_{\t_h}\lb s_{\t_h}-\sfrac{\rho_{\t_h}-\rho_\t}{2}\rb\]
given by $V_{\t_h}^{w_h}$. We want $d_r=1$. Since
\[\gamma_{\t_h}\lb s_{\t_h}-\sfrac{\rho_{\t_h}-\rho_\t}{2}\rb=-\sfrac{\epsilon_\t}{2}+\sfrac{|\t_h|-1}{2}\rho_\t\in\tfrac{1}{2}\Z\smallsetminus\Z,\]
we have $d_{r+1}=2$. As before $d_r=1$ if and only if $\frac{\rho_{\t_h}-\rho_\t}{2}=\frac{n_0}{d_0}-\frac{n_{r+1}}{d_{r+1}}\geq\frac{1}{2}$ and similarly for $\t_l$. Thus if $\t$ has two or more rational children and $\Gamma_\t$ is irreducible then it can be blown down if and only if $v(c_f)$ is odd and $\t=\roots$ is union of its $2$ odd rational children $\t_h$ and $\t_l$, satisfying $\rho_{\t_h}\geq\rho_\t+1$, $\rho_{\t_l}\geq\rho_\t+1$. This is case (\ref{3Contr}) of Definition \ref{RemovableContractibleDefinition}. 

Suppose now $\Gamma_\t$ reducible. By Proposition \ref{DefiningEquationGammaProposition} the cluster $\t$ is \"{u}bereven, $\epsilon_\t$ is even and $\Gamma_\t$ is the disjoint union of $\Gamma_\t^-\simeq\P^1$ and $\Gamma_\t^+\simeq\P^1$. As before, both $\Gamma_\t^-$ and $\Gamma_\t^+$ intersect at least other two components (given by the proper children of $\t$). But then neither $\Gamma_\t^-$ nor $\Gamma_\t^+$ has self-intersection $-1$, as $m_\t=1$.

We have showed that, for a rational cluster $\t\in\Sigma_{K^{nr}}$, an irreducible component of $\Gamma_\t$ can be blown down if and only if $\t$ is contractible. Moreover, in this case, $\Gamma_\t$ is irreducible. 
It remains to show that after blowing down all components $\Gamma_\t$ where $\t$ is a contractible cluster, no other component can be blown down. 
First note that if $\t$ is a contractible cluster, then $m_\t=2$ and $\Gamma_\t$ intersects one or two other components of multiplicity $1$ at two points in total. 
If it intersects only one component, then after the blowing down, the latter will have a node and will not be isomorphic to $\P^1$.
If $\Gamma_\t$ intersects two components and those intersect something else in $\mathcal{C}_s$, then they will not have self-intersection $-1$ also when $\Gamma_\t$ is blown down.
Therefore suppose that one of those two does not intersect any other component of $\mathcal{C}_s$. 
If we are in case (\ref{1Contr}) or case (\ref{2Contr}), it is easy to see that this never happens. Indeed, in those cases, $\Gamma_\t$ intersects non-open-ended chains of $\P^1$s.
Then without loss of generality assume to be in case (\ref{3Contr}) and that $\Gamma_{\t_h}$ is the component that can be blown down once $\Gamma_\t$ has been contracted. This implies $\s_h=\t_h$ and $\rho_{\s_h}=\rho_\t+1$. Then $b_{\s_h}=1$ and $\epsilon_{\s_h}=\epsilon_{\t}+|\s_h|$. Since both $\epsilon_\t$ and $\s_h$ are odd, we have $\epsilon_{\s_h}\in 2\Z$. So $D_{\s_h}=2$ and $\tilde\s_h$ is the set of rational children of $\s_h$.
Hence $g(\s_h)=\left\lfloor\frac{|\s_h|-1}{2}\right\rfloor\geq 1$ since $|\s_h|\geq 3$. But then $\Gamma_{\s_h}$ cannot be blown down.

\subsection{Galois action}
Consider the base extended hyperelliptic curve $C_{K^{nr}}/K^{nr}$. The rational clusters of $C_{K^{nr}}$ and their corresponding rational centres are then over $K^{nr}$. Denote $\Sigma_{K^{nr}}=\Sigma_{C_{K^{nr}}}^\mathrm{rat}$. For any proper cluster $\s\in\Sigma_{K^{nr}}$, let $G_\s=\mathrm{Stab}_{G_K}(\s)$, $K_\s=\lb K^\mathrm{s}\rb^{G_\s}$ and $k_\s$ be the residue field of $K_\s$. Let $\Sigma_{C_{K^{nr}}}^\mathrm{min}=\{\s_1,\dots,\s_m\}$ be the set of rationally minimal clusters of $C_{K^{nr}}$. Fix a set $W=\{w_1,\dots,w_m\}\subset K^{nr}$ of corresponding rational centres. By Lemma \ref{RationalCentreTameExtensionLemma}, we can assume this choice to be $G_K$-equivariant, i.e.\ for any $\sigma\in G_K$, one has $\sigma(w_l)=w_h$ if and only if $\sigma(\s_l)=\s_h$. We can also require that $w_h\in\s_h$ if $\s_h\cap K_{\s_h}\neq \varnothing$. 
Similarly, for any proper cluster $\t\in\Sigma_{K^{nr}}\setminus\Sigma_{C_{K^{nr}}}^\mathrm{min}$, fix a rational centre $w_\t$ in such a way that $w_{\sigma(t)}=\sigma(w_\t)$ for any $\sigma\in G_K$. Set $w_{\s_o}:=w_o$ for any $o=1,\dots,m$.

\begin{lem}\label{ModelHypothesisGeneralCaseLemma}
With the choices above, for any $h=1,\dots,m$, the set of proper clusters in $\Sigma_{C_{K^{nr}}}^{w_h}$ coincides with $\mathring{\Sigma}_{C_{K^{nr}}}^{w_h}$.
\proof
Suppose by contradiction that there exists a non-proper cluster $\{r\}=\s\in\Sigma_{C_{K^{nr}}}^{w_h}$, with $r\neq w_h$. Note that $r\in\s_h$ and so $\s<\s_h$. Recall that since $\s$ is a cluster centred at $w_h$, it is cut out by the disc $\mathcal{D}=\{x\in\bar K\mid v(x-w_h)\geq\rho_\s^{w_h}\}$, with $\rho_\s^{w_h}=v(r-w_h)>\rho_{\s_h}$. This implies that $w_h\notin\roots$, otherwise $w_h\in\s$ and $|\s|\geq 2$. In particular, $w_h\notin\s_h$. For our choice of $w_h$, it follows that $\s_h\cap K_{\s_h}=\varnothing$. Therefore $r\notin K_{\s_h}$ and so there exists $\sigma\in G_{\s_h}$ such that $\sigma(r)\neq r$. 
Since $w_h\in K_{\s_h}$ we have
\[v(\sigma(r)-w_h)=v(\sigma(r-w_h))=v(r-w_h)=\rho_\s^{w_h}.\]
But then $\sigma(r)\in\s$, and so $|\s|\geq 2$, a contradiction.
\endproof
\end{lem}

Assume that $C_{K^{nr}}$ is $y$-regular and has an almost rational cluster picture. By the previous lemma, from the set of rational centres $W$ we can construct the proper regular model $\mathcal{C}/O_{K^{nr}}$ of $C_{K^{nr}}$ as previously presented in this section. In this subsection we show how the Galois group $\Gal(K^{nr}/K)$ acts on the $O_{K^{nr}}$-scheme $\mathcal{C}$. Moreover, we describe the induced action of $G_k$ on the special fibre $\mathcal{C}_s/k^\mathrm{s}$, and give defining equations for the principal components of $\mathcal{C}_s$ compatibly with the action.

For any $l=1,\dots,m$, recall the notation $f_l(x)=f(x+w_l)\in K^{nr}[x]$ and $C^{w_l}/K^{nr}:y^2=f_l(x)$. Fix $\sigma\in G_K$. Let $l,h=1,\dots,m$ such that $\sigma(\s_l)=\s_h$. Then $\sigma(f_l)=f_h$. Now, let $\t\in\Sigma_{C_{K^{nr}}}^{w_l}$ be a proper cluster. 
Then $\sigma(\t)\in\Sigma_{C_{K^{nr}}}^{w_h}$ and $\rho_\t=\rho_{\sigma(\t)}$. It follows that most of the quantities attached to $\t$, such as those of Definition \ref{QuantitiesForTheoremsOnModelsDefinition}, are the same for $\sigma(\t)$, e.g.\ $\epsilon_\t=\epsilon_{\sigma(\t)}$. In particular, if $M$ is a matrix associated to $\t$ then $M$ is associated to $\sigma(\t)$ as well. So $\sigma(\mathcal{F}_M^l)=\mathcal{F}_{M}^h$. 
Finally, as $\sigma(\prod_{o\neq l}(x+w_{lo})^{-1})=\prod_{o\neq h}(x+w_{ho})^{-1}$ we also have $\sigma(T_M^l)=T_{M}^h$.

Hence the natural $K^{nr}$-isomorphism $C^{w_h}\xrightarrow{\sigma}C^{w_l}$ induces $O_{K^{nr}}$-isomorphisms of schemes
\begin{equation}\label{GaloisActionOnChartsEquation}
\mathcal{C}_\Delta^{w_h}\xrightarrow{\sigma}\mathcal{C}_\Delta^{w_l},\qquad\mathring{\mathcal{C}}_\Delta^{w_h}\xrightarrow{\sigma}\mathring{\mathcal{C}}_\Delta^{w_l},\qquad U_M^h\xrightarrow{\sigma}U_M^l.
\end{equation}
Via the glueing morphisms (\ref{GlueingMapEquation}), these maps describe the action of $G_K$ on $\mathcal{C}$.

We want to study the action of $G_k$ on the special fibre of $\mathcal{C}$ more in detail. Let $\sigma\in\Gal(K^{nr}/K)$ and let $\bar\sigma\in G_k$ corresponding to $\sigma$ via the canonical isomorphism $\Gal(K^{nr}/K)\simeq G_k$. 
Let $l,h$ and $\t$ as above. In \S\ref{ComponentsSubsection} we described closed $1$-dimensional subschemes composing $\mathring{\mathcal{C}}_{\Delta,s}^{w_l}$ and the morphisms induced by the glueing maps. Recall the polynomials introduced in Definition \ref{PolynomialsGeneralDefinition}.
From (\ref{GaloisActionOnChartsEquation}) we get
\[\bar\sigma(\ch{g_{\s_l,w_l}^0})=\ch{g_{\s_h, w_h}^0},\quad\bar\sigma(\ch{g_{\t,w_l}})=\ch{g_{\sigma(\t),w_h}},\quad\bar\sigma(\ch{g_l|_{L_{\t}^{w_l}}})=\ch{g_h|_{L_{\sigma(\t)}^{w_h}}}.\]
From the equality (\ref{ghftEquality}) we obtain $\bar\sigma(f_{\t,w_l})=f_{\sigma(\t),w_h}$. Note that the previous relations can also be recovered directly from the definitions.


\begin{lem}\label{GaloisActionOnPolyLemma}
Let $w_\t$ be the rational centre of $\t$ fixed above. Then
\begin{enumerate}[label=(\roman*)]
    \item $\overline{g_{\t,w_\t}},\overline{f_{\t,w_\t}}\in k_\t[X]$; 
    \item 
    $\overline{g_{\t,w_\t}}=\overline{g_{\t,w_l}}$ and $   \overline{f_{\t,w_\t}}(X)=\overline{f_{\t,w_l}}(X+\overline{u_{w_\t w_l}})$ where $\overline{u_{w_\t w_l}}=\tfrac{w_\t-w_l}{\pi^{\rho_\t}}\mod\pi$; 
\end{enumerate}
\proof
For any rational centre $w$ of $\t$, let $u_{\t,w}=c_f\tprod_{r\in\roots\setminus\t} (w-r)$ as in Definition \ref{PolynomialsGeneralDefinition}. Note that $u_{\t,w}/\pi^{v(u_{\t,w})}$ is independent of $w$ since 
\[v((w_\t-r)-(w-r))=v(w_\t-w)\geq\rho_\t>v(w_\t-r)\] for any $r\in\roots\setminus\t$. Then $\overline{g_{\t,w_\t}}=\overline{g_{\t,w_l}}$. 
If $\bar\sigma\in \Gal(k^\mathrm{s}/k_\t)$, i.e.\ $\sigma\in\Gal(K^{nr}/K_\t)$, then \[\bar\sigma(\overline{g_{\t,w_\t}})=\bar\sigma(\overline{g_{\t,w_l}})=\overline{g_{\t,w_h}}=\overline{g_{\t,w_{\t}}}.\]
In particular $\overline{g_{\t,w_\t}}\in k_\t[X]$.

Since $u_{\t,w}/\pi^{v(u_{\t,w})}$ is independent of $w$ we also have
\[\overline{f_{\t,w_\t}}(X^{b_\t})=
\overline{f_{\t,w_l}}((X+\ch{u_{w_\t w_l}})^{b_\t}).\]
Suppose $\rho_\t\in\Z$. Then $b_\t=1$ and so the equality above implies $\overline{f_{\t,w_\t}}(X)=\overline{f_{\t,w_l}}(X+\ch{u_{w_\t w_l}})$. Suppose $\rho\notin\Z$. Then $v(w-w_\t)>\rho_\t$ for any rational centre $w$ of $\t$ as $v(w-w_\t)\in\Z$ and $v(w-w_\t)\geq\rho_\t$. Hence $\ch{u_{w_\t w_l}}=0$. Thus $\overline{f_{\t,w_\t}}(X^{b_\t})=\overline{f_{\t,w_l}}(X^{b_\t})$, which implies $\overline{f_{\t,w_\t}}(X)=\overline{f_{\t,w_l}}(X)=\overline{f_{\t,w_l}}(X+\ch{u_{w_\t w_l}})$. 
If $\bar\sigma\in \Gal(k^\mathrm{s}/k_\t)$, i.e.\ $\sigma\in\Gal(K^{nr}/K_\t)$, then
\[\bar\sigma(\overline{f_{\t,w_\t}})(X)=\bar\sigma(\overline{f_{\t,w_l}})(X+\bar\sigma(\ch{u_{w_\t w_l}}))=\overline{f_{\t,w_h}}(X+\ch{u_{w_\t w_h}})=\overline{f_{\t,w_\t}}(X),\]
and so $\overline{f_{\t,w_\t}}\in k_\t[X]$.
\endproof
\end{lem}

\begin{rem}
Note that the polynomials $\ch{f_{\t,w_\t}}$, $\ch{g_{\t,w_\t}}$ and $\ch{g_{\s_h,w_h}^0}$ equal the polynomials $\ch{f_{\t}}$, $\ch{g_{\t}}$ and $\ch{g_{\s_h}^0}$ given in Definition \ref{PolynomialsDefinition}.
\end{rem}

Let $V=V_\t^{w_l}$ and consider the subscheme $X_V\times\P_V$ of $\mathcal{C}_s$ given by $V$, where $\P_V$ is a chain of $\P^1$s and $X_V:\{\overline{g_{\t,w_l}}=0\}$ over $\G_{m,k^\mathrm{s}}$. 
If $\s_h\subset\t$, then the glueing map $U_M^h\rightarrow U_M^l$ induces the identity $\phi_V^{hl}:X_{V_{\t}^{w_h}}\xrightarrow{=}X_{V_\t^{w_l}}$. Define $X_\t\subseteq \G_{m,k^\mathrm{s}}$ given by $g_{\t,w_\t}=0$. By Lemma \ref{GaloisActionOnPolyLemma}, $\phi_V^o:X_\t\xrightarrow{\simeq} X_{V_\t^{w_o}}$, for $o=h,l$, and this isomorphism is compatible with the Galois action and the glueing maps, i.e.\  $\sigma\circ\phi_V^h=\phi_V^l\circ\sigma$ and $\phi_V^{hl}\circ\phi_V^h=\phi_V^l$ as morphisms on $X_\t$.

When $V_0=V_0^{w_l}$ we can consider the subscheme $X_{V_0}\times\P_{V_0}$ given by $V_0$, where $\P_{V_0}$ is a chain of $\P^1$s and $X_{V_0}:\{\overline{g_{\s_l,w_l}}=0\}$ over $\G_{m,k^\mathrm{s}}$.
Since $X_{V_0}\times\P_{V_0}$ is never glued to any other component there is no need to choose a different model for it.

Let $L=L_\t^{w_l}$ and consider the subscheme $X_{L}^W\times\P_{L}$ given by $L$, where $\P_{L}$ is a chain of $\P^1$s and $X_{L}^W:\{\overline{f_{\t,w_l}}=0\}$ over $\A^1_{k^\mathrm{s}}$. If $\s_h\subset\t$, then the isomorphism $\phi_L^{hl}:X_{L_\t^{w_h}}^W\xrightarrow{\simeq} X_{L_\t^{w_l}}^W$ given by the glueing map $U_M^h\rightarrow U_M^l$ is induced by the ring isomorphism $k^\mathrm{s}[X]\rightarrow k^\mathrm{s}[X]$, sending $X\mapsto X+\overline{u_{w_hw_l}}$, where $\overline{u_{w_hw_l}}=\frac{w_h-w_l}{\pi^{\rho_\t}}\mod\pi$. Define $X_\t^W\subseteq \A^1_{k^\mathrm{s}}$ given by $\ch{f_{\t,w_\t}}=0$. By Lemma \ref{GaloisActionOnPolyLemma}, the map $X\mapsto X+\overline{u_{w_\t w_l}}$ induces an isomorphism $\phi_L^o:X_\t^W\xrightarrow{\simeq} X_{L_\t^{w_o}}^W$, for $o=h,l$, compatible with the Galois action and the glueing morphisms, i.e.\ $\sigma\circ\phi_L^h=\phi_L^l\circ\sigma$ and $\phi_L^{hl}\circ\phi_L^h=\phi_L^l$ as morphisms on $X_\t^W$.

Recall the definitions of $\hat\t^W$ and $\G_{\t,w_l}\subseteq\A_{k^{\mathrm{s}}}^1$ given in Definition \ref{hattGtDefinition} and the definition of $\mathring\t$ given in Definition \ref{PolynomialsDefinition}. Note that by Lemma \ref{ModelHypothesisGeneralCaseLemma},
\[\hat\t^W=\{\s\in\Sigma_{K^{nr}}\cup\{\varnothing\}\mid \s<\t\}\setminus\{\{r\}\in\Sigma_{K^{nr}}\mid r\notin W\}.\]
Fix $c=0,\dots,b_\t-1$ such that $1/b_\t-c\rho_\t\in\Z$. For any rational centre $w\in K^{nr}$ of $\t$ define $\hat f_{\t,w}\in k^\mathrm{s}[X,Y]$ by \[\hat f_{\t,w}(X)=\prod_{\s\in\hat\t^W}(X-\ch{u_{w_\s w}})^{\tfrac{|\s|}{b_\t}-c\epsilon_\t}\ch{f_{\t,w}}(X),\] where $\ch{u_{w_\s w}}=\tfrac{w_\s-w}{\pi^{\rho_\t}}\mod \pi$ 
($w_\s=w_l$ if $\s=\varnothing$).
Let $L=L_\t^{w_l}$, $F=F_\t^{w_l}$ and $M=M_{L,0}$. It follows from Lemma \ref{DefiningEquationXFLemma} that the scheme $\Gamma_\t^{w_l}=\mathring{X}_F\cap U_M^{l}$ is given by $Y^{D_\t}=\hat f_{\t,w_l}(X)$ as a subscheme of $\G_{\t,w_l}\times \A_{k^{\mathrm{s}}}^1$. We then obtain $\bar\sigma(\hat f_{\t,w_l})=\hat f_{\sigma(\t),w_h}$ 
from the action (\ref{GaloisActionOnChartsEquation}) of $\sigma$ on $U_M^l$.

\begin{lem}\label{HatfLemma}
With the notation above,
\begin{enumerate}[label=(\roman*)]
    \item $\hat f_{\t,w_\t}\in k_\t[X]$;
    \item $\hat f_{\t,w_\t}(X)=\hat f_{\t,w_l}(X+\ch{u_{w_\t w_l}})$ where $\overline{u_{w_\t w_l}}=\tfrac{w_\t-w_l}{\pi^{\rho_\t}}\mod\pi$; 
\end{enumerate}
\proof
If $\s\in\mathring\t$, then $\sigma(\s)\in\mathring{(\sigma(\t))}$ and $\bar\sigma(\ch{u_{w_\s w}})=\ch{u_{w_{\sigma(\s)}\sigma(w)}}$ for any rational centre $w$ of $\t$. Hence $\hat f_{\t,w_\t}\in k_\t[X]$ and $\bar\sigma(\hat f_{\t,w_l})=\hat f_{\sigma(\t),w_h}$ by Lemma \ref{GaloisActionOnPolyLemma}(i),(iii). Since $\ch{u_{w_\s w_\t}}=\ch{u_{w_\s w_l}}-\ch{u_{w_\t w_l}}$,
Lemma \ref{GaloisActionOnPolyLemma}(ii) implies $\hat f_{\t,w_\t}(X)=\hat f_{\t,w_l}(X+\ch{u_{w_\t w_l}})$.
\endproof
\end{lem}

Define $\Gamma_\t^{w_\t}\subset\G_{\t,w_\t}\times\A^1_{k^\mathrm{s}}$ given by $Y^{D_\t}=\hat f_{\t,w_\t}$. Suppose $\s_h\subset\t$, and let $\phi_\t^{hl}:\Gamma_\t^{w_h}\simeq\Gamma_\t^{w_l}$ be the isomorphism coming from the glueing map $U_M^h\rightarrow U_M^l$ induced by the ring homomorphism $X\mapsto X+\overline{u_{w_hw_l}}$. By Lemma \ref{HatfLemma}, the map $X\mapsto X+\overline{u_{w_\t w_l}}$ induces an isomorphism $\phi_\t^o:\Gamma_\t^{w_\t}\simeq \Gamma_\t^{w_o}$, for $o=h,l$, which is compatible with the Galois action and the glueing maps, i.e.\ $\sigma\circ\phi_\t^h=\phi_\t^l\circ\sigma$ and $\phi_\t^{hl}\circ\phi_\t^h=\phi_\t^l$ as morphisms on $\Gamma_\t^{w_\t}$. Therefore $\Gamma_\t$ is isomorphic to the smooth completion of $\Gamma_\t^{w_\t}$, and so it is given by $Y^{D_\t}=\tilde f_\t(X)$, where $\tilde f_\t(X)=\prod_{\s\in\mathring\t}(X-\ch{u_{w_\s w_\t}})\ch{f_{\t,w_\t}}(X)$ is the polynomial given in Definition \ref{PolynomialsDefinition}.

\section{Integral differentials}\label{IntegralDifferentialsSection}
Let $C$ be a hyperelliptic curve of genus $g\geq 1$ defined over $K$ by a Weierstrass equation $y^2=f(x)$. It is well-known that the $K$-vector space of global sections of the sheaf of differentials of $C$, namely $H^0(C,\Omega^1_{C/K})$, is spanned by the basis \[\underline{\omega}=\left\{\sfrac{dx}{2y},x\sfrac{dx}{2y},\dots,x^{g-1}\sfrac{dx}{2y}\right\}.\]
Let $\mathcal{C}$ be a regular model of $C$ over $O_K$ and consider its canonical (or dualising) sheaf $\omega_{\mathcal{C}/O_K}$. The free $O_K$-module of its global sections $H^0(\mathcal{C},\omega_{\mathcal{C}/O_K})$ can be viewed as an $O_K$-lattice in $H^0(C,\Omega^1_{C/K})$ (see \cite[Corollary 9.2.25(a)]{Liu}). The aim of this section is to present a basis of $H^0(\mathcal{C},\omega_{\mathcal{C}/O_K})$ as an $O_K$-linear combination of the elements in $\underline{\omega}$ under the assumptions of Theorem \ref{MinimalRegularSNCModelTheorem}. Note that by \cite[Corollary 9.2.25(b)]{Liu} the problem is independent of the choice of model $\mathcal{C}$ but it does depend on the choice of the equation $y^2=f(x)$ since the basis $\underline{\omega}$ does. Throughout this section let $C$ and $\mathcal{C}/O_K$ be as above. 

If $C$ is $\Delta_v$-regular, \cite[Theorem 8.12]{Dok} gives an $O_K$-basis of $H^0(\mathcal{C},\omega_{\mathcal{C}/O_K})$, as required. We rephrase it in terms of rational cluster invariants, by using results of \S \ref{ClustersSection} and Lemma \ref{Regularity-RegularityLemma}.


\begin{thm}\label{DifferentialsNestedTheorem}
Suppose $C$ has an almost rational cluster picture and is $y$-regular, and all proper clusters $\s\in\Sigma_C$ have same rational centre $w\in K$. Let $\s_1\subset\dots\subset \s_n=\roots$ be the proper clusters in $\Sigma_C^\mathrm{rat}$. For every $j=0,\dots,g-1$, define \[i_j:=\min\{i\in\{1,\dots,n\}\mid 2(j+1)<|\s_i|\}\] and 
\[e_j:=\tfrac{1}{2}\epsilon_{\s_{i_j}}-(j+1)\rho_{\s_{i_j}}.\] Then the differentials
\[\mu_j=\pi^{\lfloor e_j \rfloor} (x-w)^j\sfrac{dx}{2y}\qquad j=0,\dots,g-1,\]
form an $O_K$-basis of $H^0(\mathcal{C},\omega_{\mathcal{C}/O_K})$.
\proof
Let $C^w:y^2=f(x+w)$ be the hyperelliptic curve isomorphic to $C$ through the change of variable $y\mapsto y,$ $x\mapsto x+w$. By Corollary \ref{RegularityAfterTranslationCorollary} and Lemma \ref{Regularity-RegularityLemma}, the curve $C^w$ is $\Delta_v$-regular. Since $\mathring\Sigma_C^\mathrm{rat}$ consists of the proper clusters in $\Sigma_C^w$, Lemma \ref{DescriptionNewtonPolytopeClustersLemma} and \cite[Theorem 8.12]{Dok} implies that
\[\mu_j=\pi^{\lfloor e_j \rfloor} x^j\sfrac{dx}{2y}\qquad j=0,\dots,g-1,\]
form an $O_K$-basis of $H^0(\mathcal{C},\omega_{\mathcal{C}/O_K})$ as a lattice in $H^0(C^w,\Omega^1_{C^w/K})$ (that is if $\mathcal{C}$ is regarded as a model of $C^w$). Changing variables concludes the proof.
\endproof
\end{thm}

Suppose now $C$ has an almost rational cluster picture and is $y$-regular. Let $\Sigma_C^\mathrm{min}$ be the set of rationally minimal clusters and let $W=\{w_\s\,\mid\,\s\in\Sigma_C^\mathrm{min}\}$ be a corresponding set of rational centres, such that all clusters in $\mathring\Sigma_C^{w_\s}$ are proper. For every proper cluster $\t\in\Sigma_C^\mathrm{rat}$, choose a minimal cluster $\s\subseteq\t$ and set $w_\t:=w_\s$. Consider the regular model $\mathcal{C}/O_K$ of $C$ of Theorem
\ref{ConstructionProperModelGeneralCaseTheorem}, constructed in \S \ref{ConstructionModelsSection} by glueing the open subschemes $\mathring{\mathcal{C}}_\Delta^{w}$ of $\mathcal{C}_\Delta^w$ for $w\in W$. We want to describe the canonical morphism $C\rightarrow\mathcal{C}$. Write $W=\{w_1,\dots,w_m\}$ as in \S\ref{ConstructionModelsSection}. For any $h=1,\dots,m$, let $\t\in\Sigma_C^{w_h}$ be a proper cluster and let $M$ be a matrix associated to $\t$. Let $C^{w_h}:y^2=f(x+w_h)$ and 
\[y^2-f(x+w_h)\stackrel{M}{=}Y^{n_Y}Z^{n_Z}\mathcal{F}_M^h(X,Y,Z).\]
Then, on the affine chart $X_M$ the projection $C\rightarrow\mathcal{C}_\Delta^{w_h}$ is induced by
\[\frac{R}{\lb\mathcal{F}_M^h(X,Y,Z)\rb}\xrightarrow{M}\frac{K[(x')^{\pm 1},(y')^{\pm 1}]}{\lb(y')^2-f(x'+w_h)\rb}\xrightarrow{\simeq}\frac{K[x^{\pm 1},y^{\pm 1}]}{\lb y^2-f(x)\rb},\]
where $(X,Y,Z)=(x',y',\pi)\bullet M$ and $(x',y')=(x-w_h,y)$. In particular it follows that $(X,Y,Z)=(x-w_h,y,z)\bullet M$ and so 
\[\begin{pmatrix}
     x-w_h  \\
     y  \\
     \pi
\end{pmatrix}=
\begin{pmatrix}
X^{\tilde m_{11}}Y^{\tilde m_{21}}Z^{\tilde m_{31}}\\
X^{\tilde m_{12}}Y^{\tilde m_{22}}Z^{\tilde m_{32}}\\
X^{\tilde m_{13}}Y^{\tilde m_{23}}Z^{\tilde m_{33}}
\end{pmatrix}
=\begin{pmatrix}
     X  \\
     Y  \\
     Z
\end{pmatrix}\bullet M^{-1}.
\]

For a proper cluster $\t\in\Sigma_C^\mathrm{rat}$ recall the definitions of $\Gamma_\t$ and $m_\t$.

\begin{prop}\label{OrdDifferentialsProposition}
Let $\t\in\Sigma_C^\mathrm{rat}$ be a proper cluster. Then\footnote{If $\Gamma_\t$ is reducible, say $\Gamma_\t= \Gamma_\t^-\cup\Gamma_\t^+$, with $\ord_{\Gamma_\t}(\cdot)$ we mean $\min\{\ord_{\Gamma_\t^-}(\cdot),\ord_{\Gamma_\t^+}(\cdot)\}$}
\begin{align*}
     &\ord_{\Gamma_\t}(x-w_\s)=m_{\t}\rho_\t,\\
     &\ord_{\Gamma_\t}\tfrac{dx}{2y}=-m_\t\left(\tfrac{1}{2}\epsilon_\t-\rho_\t-1\right)-1.
\end{align*}
for every proper cluster $\s\in\Sigma_C^\mathrm{rat}$, $\s\subseteq \t$.
\proof
Let $g(x,y)=y^2-f(x)$. Let $W=\{w_1,\dots,w_m\}$ as above. Let $h=1,\dots,m$ such that $w_h=w_\s$. 
Let $F=F_\t^{w_h}$, $V=V_\t^{w_h}$, $M=M_{V,0}$ and let $X,Y,Z$ be the transformed variables $(X,Y,Z)=(x-w_\s,y,\pi)\bullet M$. 
Define $\mathcal{H}(X,Y,Z)=\pi - X^{\tilde m_{13}}Y^{\tilde m_{23}}Z^{\tilde m_{33}}$, and
$\mathcal{G}(X,Y,Z)=g((X,Y,Z)\bullet M^{-1})$. We have \[\mathcal{F}_M^h(X,Y,Z)=Y^{-n_Y}Z^{-n_Z}\mathcal{G}(X,Y,Z),\]  where 
$n_Z=m_\t\epsilon_\t$, since $\mathrm{ord}_Z(y^2)=m_\t\epsilon_\t$ for Lemma \ref{MatricesLemma}. Write $\mathcal{F}=\mathcal{F}_M^h$.

Note that $d(x-w_\s)=dx$ and $(g_{w_\s})_x'(x-w_\s)=g_x'(x)$, where $g_{w_\s}(x,y)=g(x+w_\s,y)$. Then, by \cite[8.7]{Dok},
\[
\begin{cases}
(x-w_\s)g_x'=m_{11}X\mathcal{G}_X'+m_{12}Y\mathcal{G}_Y'+m_{13}Z\mathcal{G}_Z'\\
yg_y'=m_{21}X\mathcal{G}_X'+m_{22}Y\mathcal{G}_Y'+m_{23}Z\mathcal{G}_Z'
\end{cases}\]
from which it follows that
\begin{align*}m_{11}yg_y'-m_{21}(x-w_\s)g_x'&=(m_{11}m_{22}-m_{21}m_{12})Y\mathcal{G}_Y'-(m_{21}m_{13}-m_{11}m_{23})Z\mathcal{G}_Z'\\
&=\tilde m_{33}Y\mathcal{G}_Y'-\tilde m_{23}Z\mathcal{G}_Z'.
\end{align*}
We will show later that this quantity is non-zero. Moreover,
\[\tilde m_{33}Y\mathcal{G}_Y'-\tilde m_{23}Z\mathcal{G}_Z'=Y^{n_Y}Z^{n_Z}\lb\tilde m_{33}Y\mathcal{F}_Y'-\tilde m_{23}Z\mathcal{F}_Z'+(n_Y+n_Z)\mathcal{F}\rb.\]
Recall that $X=(x-w_\s)^{m_{11}}y^{m_{21}}\pi^{m_{31}}$. Then $\frac{dX}{X}=m_{11}\tfrac{dx}{x-w_\s}+m_{21}\tfrac{dy}{y}$. Furthermore, as $0=dg=g_x'dx+g_y'dy$ in $\Omega_{C/K}$, we have
\[\sfrac{dX}{X}=\lb\sfrac{m_{11}}{x-w_\s}-\sfrac{m_{21}}{y}\sfrac{g_x'}{g_y'}\rb dx=\sfrac{dx}{(x-w_\s)yg_y'}\lb m_{11}yg_y'-m_{21}(x-w_\s) g_x'\rb.\]
Therefore
\begin{equation}\label{BasicDifferentialOrderEquation}
    \frac{dx}{2(x-w_\s)y^2}=
    \frac{dX}{XY^{n_Y}Z^{n_Z}\lb\tilde m_{33}Y\mathcal{F}_Y'-\tilde m_{23}Z\mathcal{F}_Z'+(n_Y+n_Z)\mathcal{F}\rb}.
\end{equation}


Let $S=\Spec O_K$. 
Considering $X^{-1}$ as an independent variable, the scheme 
\[U=\Spec\frac{O_K[Y,Z,X^{-1},X]}{(\mathcal{F},\mathcal{H},X\cdot X^{-1}-1)}\] defines a complete intersection in $\A_S^4$. Furthermore, $U$ is an open subscheme of $\mathcal{C}_\Delta^{w_h}\cap X_M$ that restricted to $\A_S^4\setminus\{T_M^h(X,Y,Z)=0\}$ equals $U_M^h$.
In particular, $U$ is integral. From \S\ref{SpecialFibreSubsection} it follows that $U_\t=U\cap\{Z=0\}$ is a dense open subset of $\mathring X_{F}$. 
Recall that $\mathring X_{F}$ is an open subscheme of $\Gamma_\t$. Hence it suffices to prove the proposition for $U_\t$ instead of $\Gamma_\t$ (\cite[Lemma 9.2.17(a)]{Liu}).
Since $X$ and $Y$ are units and $Z$ vanishes to order $1$ on $U_\t$, Lemma \ref{MatricesLemma} implies that
\begin{equation}\label{OrderxandyEquation}
\ord_{U_\t}(x-w_\s)=\tilde m_{31}=m_\t\rho_\t\quad\mbox{ and }\quad \ord_{U_\t}y=\tilde m_{32}=m_\t\tfrac{\epsilon_\t}{2}.
\end{equation}

Recall that $U$ is integral and that $U_\eta$ is isomorphic to an open subscheme of $C$. Then $U_\eta$ is smooth. Hence, by \cite[Corollary 6.4.14(b)]{Liu}, the sheaf $\omega_{\mathcal{C}/O_K}$ is generated on $U$ by $\mathcal{E}^{-1}dX$ where \[\mathcal{E}:=\left|\begin{matrix}\mathcal{F}'_Y&\mathcal{F}'_Z&\mathcal{F}_{X^{-1}}'\cr H'_Y&H'_Z&\mathcal{F}_{X^{-1}}'\cr
0&0&X\end{matrix}\right|=-\pi X Y^{-1}Z^{-1}\lb \tilde m_{33}Y\mathcal{F}_Y'-\tilde m_{23}Z\mathcal{F}_Z'\rb,\] if $\mathcal{E}$ is non-zero. 
Suppose it is; we are going to prove it later. Thus, as $\mathcal{F}=0$ on $U$, from (\ref{BasicDifferentialOrderEquation}) and (\ref{OrderxandyEquation}) we have
\[\ord_{U_\t}\sfrac{dx}{2y}=m_\t\left(\tfrac{1}{2}\epsilon_\t+\rho_\t\right)+\tilde m_{33}-n_Z-1=m_\t\lb-\tfrac{1}{2}\epsilon_\t+\rho_\t+1\rb-1.\]

It remains to show that $\mathcal{E}$ does not equal $0$ on $U$. Suppose it does. Then from the computations above, it follows that $m_{11}y g_y'-m_{21}(x-w_\s)g_x'=0$ in $K(C)$. 
Since $m_{21}$ equals either $1$ or $2$ by Lemma \ref{MatricesLemma}, 
it follows that there exists a non-zero $c\in K$, such that 
\[m_{11} y g_y'-m_{21}(x-w_\s)g_x'+c g=0\]
($c\in K$ from degree analysis). Then $cf(x)=m_{21}(x-w_\s)f'(x)$. Note that $m_{21}$ is non-zero as $\mathrm{char}(K)\neq 2$. But then a contradiction follows since $f$ is a separable polynomial of degree $\geq 3$.
\endproof
\end{prop}

In the following two theorems we describe a basis of integral differentials of $C$.
We use Definitions/Notations \ref{ClusterDepthDefinition}, \ref{ParentChildWedgeDefinition},  \ref{ClusterPictureDefinition}, \ref{RadiusRationalCentreDefinition}, \ref{RationalClusterPictureDefinition}, \ref{AlmostRationalDefinition}, \ref{QuantitiesForTheoremsOnModelsDefinition}, \ref{yRegularDefinition} in the statements.

\begin{thm}\label{DifferentialTheorem}
Let $C/K$ be a hyperelliptic curve of genus $g\geq 1$ defined by the Weierstrass equation $y^2=f(x)$ and let $\mathcal{C}/O_K$ be a regular model of $C$. Suppose $C$ has an almost rational cluster picture and is $y$-regular. For $i=0,\dots,g-1$ choose inductively proper clusters $\s_i\in\Sigma_C^\mathrm{rat}$ so that
\[e_i:=\frac{\epsilon_{\s_i}}{2}-\sum_{j=0}^i\rho_{\s_j\wedge\s_i}=\max_{\t\in\Sigma_C^\mathrm{rat}}\bigg\{\frac{\epsilon_{\t}}{2}-\rho_\t-\sum_{j=0}^{i-1}\rho_{\s_j\wedge\t}\bigg\},\]
where if $\s$ and $\s'$ are two possible choices for $\s_i$  satisfying $\s'\subset\s$, then choose $\s_i=\s$.
Then the differentials
\[\mu_i=\pi^{\lfloor e_i\rfloor}\prod_{j=0}^{i-1}(x-w_{\s_j})\sfrac{dx}{2y},\qquad i=0,\dots,g-1,\]
form an $O_K$-basis of $H^0(\mathcal{C},\omega_{\mathcal{C}/O_K})$.
\proof
Since $H^0(\mathcal{C},\omega_{\mathcal{C}/O_K})$ is independent of the choice of regular model, we consider $\mathcal{C}$ to be the model described in Theorem
\ref{ConstructionProperModelGeneralCaseTheorem}
and constructed in \S\ref{ConstructionModelsSection}. 

We first show that the differentials $\mu_i$ are global sections of $\omega_{\mathcal{C}/O_K}$. It suffices to prove they are regular along all components $\Gamma_\t$, where $\t\in\Sigma_C^\mathrm{rat}$ proper. Indeed for the construction of $\mathcal{C}$ and the definition of the $e_i$'s, the differentials $\mu_i$ are regular along all other components of $\mathcal{C}_s$ by Theorem \ref{DifferentialsNestedTheorem}.
Fix $i=1,\dots,g-1$ and let $j=0,\dots,i-1$. Let $\t\in\Sigma_C^\mathrm{rat}$ be a proper cluster. If $\s_j\subseteq\t$ then 
\[\ord_{\Gamma_\t}(x-w_{\s_j})=m_\t\rho_\t=m_\t\rho_{\s_j\wedge \t},\]
by Proposition \ref{OrdDifferentialsProposition}. If $\t\subsetneq\s_j$ then $w_\t$ is a rational centre of $\s_j$. Hence
\[ v(w_\t-w_{\s_j})\geq\min_{r\in\t}\min\{v( r-w_\t),v( r-w_{\s_j})\}\geq\min\{\rho_\t,\rho_{\s_j}\}=\rho_{\s_j}=\rho_{\s_j\wedge\t}.\]
Therefore Proposition \ref{OrdDifferentialsProposition} implies 
\begin{align*}
\ord_{\Gamma_\t}(x-w_{\s_j})&\geq\min\{\ord_{\Gamma_\t}(x-w_\t),\ord_{\Gamma_\t}(w_\t-w_{\s_j})\}\\
&\geq\min\{m_\t\rho_\t,m_\t\rho_{\s_j\wedge\t}\}=m_\t\rho_{\s_j\wedge \t}.
\end{align*}
If $\s_j\nsubseteq\t$ and $\t\nsubseteq\s_j$ then from Lemma \ref{TwoChildIntegralRadiusLemma} it follows that
\[
\ord_{\Gamma_\t}(x-w_{\s_j})=\min\{m_\t\rho_\t,m_\t\rho_{\s_j\wedge\t}\}=m_\t\rho_{\s_j\wedge \t}.
\]
as $\rho_\t>\rho_{\s_j\wedge\t}$.
Thus we have proved that
\begin{equation}\label{equation1}
\ord_{\Gamma_\t}(x-w_{\s_j})\geq m_\t\rho_{\s_j\wedge \t},\qquad\text{where the equality holds if $\t\not\subset\s_j$.}
\end{equation}

Hence it follows from Proposition \ref{OrdDifferentialsProposition} that
\[
\ord_{\Gamma_\t}\mu_i\geq m_\t\bigg(\lfloor e_i\rfloor+\sum_{j=0}^{i-1}\rho_{\s_j\wedge\t}-\frac{\epsilon_\t}{2}+\rho_t+1\bigg)-1.
\]
But 
\[\lfloor e_i\rfloor\geq\bigg\lfloor \frac{\epsilon_{\t}}{2}-\rho_\t-\sum_{j=0}^{i-1}\rho_{\s_j\wedge\t}\bigg\rfloor>\frac{\epsilon_{\t}}{2}-\rho_\t-\sum_{j=0}^{i-1}\rho_{\s_j\wedge\t}-1,\]
then $\ord_{\Gamma_\t}\mu_i>-1$, that implies $\ord_{\Gamma_\t}\mu_i\geq 0$, as required.

Now we need to show that the differentials $\mu_i$ span $H^0(\mathcal{C},\omega_{\mathcal{C}/O_K})$, i.e.\ the lattice they span is saturated in the global sections of $\omega_{\mathcal{C}/O_K}$. Suppose not. Then there exist $I\subseteq\{0,\dots,g-1\}$ and $u_i\in O_K^\times$ for $i\in I$ such that the differential
$\frac{1}{\pi}\sum_{i\in I}u_i\mu_i$
is regular along $\Gamma_\t$, for every proper cluster $\t\in\Sigma_C^\mathrm{rat}$.
First we want to show that for any $i_1,i_2=0,\dots,g-1$ with $i_1<i_2$, one has $\s_{i_2}\not\subset\s_{i_1}$. Suppose by contradiction that $\s_{i_2}\subsetneq\s_{i_1}$. Then
\footnotesize{\begin{align*}
e_{i_2}&\geq\frac{\epsilon_{\s_{i_1}}}{2}-\rho_{\s_{i_1}}-\sum_{j=0}^{i_2-1}\rho_{\s_j\wedge\s_{i_1}}= e_{i_1}-\rho_{\s_{i_1}}-\sum_{j=i_1+1}^{i_2-1}\rho_{\s_j\wedge\s_{i_1}}\geq e_{i_1}-\rho_{\s_{i_1}}-\sum_{j=i_1+1}^{i_2-1}\rho_{\s_j\wedge\s_{i_2}}\\
&\geq\frac{\epsilon_{\s_{i_2}}}{2}-\rho_{\s_{i_2}}-\sum_{j=0}^{i_1-1}\rho_{\s_j\wedge\s_{i_2}}-\rho_{\s_{i_1}}-\sum_{j=i_1+1}^{i_2-1}\rho_{\s_j\wedge\s_{i_2}}=\frac{\epsilon_{\s_{i_2}}}{2}-\sum_{j=0}^{i_2}\rho_{\s_j\wedge\s_{i_2}}=e_{i_2}.
\end{align*}}\normalsize
Therefore 
\[\max_{\t\in\Sigma_C^\mathrm{rat}}\bigg\{\frac{\epsilon_{\t}}{2}-\rho_\t-\sum_{j=0}^{i_2-1}\rho_{\s_j\wedge\t}\bigg\}=e_{i_2}=\frac{\epsilon_{\s_{i_1}}}{2}-\rho_{\s_{i_1}}-\sum_{j=0}^{i_2-1}\rho_{\s_j\wedge\s_{i_1}},\]
and this means that $\s_{i_1}$ is a possible choice for the $i_2$-th cluster $\s_{i_2}$. But $\s_{i_2}\subsetneq\s_{i_1}$, so the $i_2$-th cluster should have been $\s_{i_1}$, a contradiction.

Let $I_0\subseteq I$ be the set of indices $i$ such that $\gamma_i:=e_i-\lfloor e_i\rfloor$ is maximal. Let $i_0=\min I_0$ and let $\Gamma_0=\Gamma_{\s_{i_0}}$.
Since $\s_{i_0}\not\subset\s_j$, for all $j=0,\dots,i_0-1$, from (\ref{equation1}) it follows that
\begin{align*}
m:=\ord_{\Gamma_0}\frac{1}{\pi}\mu_{i_0}&=-m_{\s_{i_0}}\gamma_{i_0}+m_{\s_{i_0}}\bigg(e_{i_0}-\frac{\epsilon_{\s_{i_0}}}{2}+\rho_{\s_{i_0}}+\sum_{j=0}^{i_0-1}\rho_{\s_j\wedge\s_{i_0}}\bigg)-1\\
&=-m_{\s_{i_0}}\gamma_{i_0}-1<0.
\end{align*}
Furthermore, 
\begin{align*}
\ord_{\Gamma_0}\frac{1}{\pi}\mu_{i}&\geq-m_{\s_{i_0}}\gamma_{i}+m_{\s_{i_0}}\bigg(e_{i}-\frac{\epsilon_{\s_{i_0}}}{2}+\rho_{\s_{i_0}}+\sum_{j=0}^{i-1}\rho_{\s_j\wedge\s_{i_0}}\bigg)-1\\
&\geq-m_{\s_{i_0}}\gamma_{i}-1\geq-m_{\s_{i_0}}\gamma_{i_0}-1=m,
\end{align*}
for all $i\in I$. 
Let $J:=\{i\in I\mid \ord_{\Gamma_0}\frac{1}{\pi}\mu_i=m\}$. Then $J\neq\varnothing$ since $i_0\in J$. 
The order of the differential $\frac{1}{\pi}\sum_{i\in J} u_i\mu_i$ along $\Gamma_0$ must be $>m$. Let $i\in I$. 
From the computations above $i\in J$ if and only if
\begin{enumerate}[label=(\roman*)]
    \item $\ord_{\Gamma_0}(x-w_{\s_j})=m_{\s_{i_0}}\rho_{\s_{i_0}\wedge\s_j}$ for all $j=0,\dots,i-1$. Equivalently, if $\s_j\supsetneq\s_{i_0}$ for some $j<i$, then $v(w_{\s_{i_0}}-w_{\s_j})=\rho_{\s_{i_0}\wedge\s_j}$.
    \item $e_i=\frac{\epsilon_{\s_{i_0}}}{2}-\rho_{\s_{i_0}}-\sum_{j=0}^{i-1}\rho_{\s_j\wedge\s_{i_0}}$. In particular, if $\s_i\subseteq\s_{i_0}$, then $\s_i=\s_{i_0}$.
    \item $\gamma_i=\gamma_{i_0}$. Equivalently, $i\in I_0$.
\end{enumerate}
Therefore $J\subseteq I_0$, $i_0=\min J$ and 
\[\lfloor e_i\rfloor-\lfloor e_{i_0}\rfloor=e_i-\gamma_i-e_{i_0}+\gamma_{i_0}=e_i-e_{i_0}=-\sum_{j=i_0}^{i-1}\rho_{\s_j\wedge\s_{i_0}},\]
for all $i\in J$. Hence
\[\frac{1}{\pi}\sum_{i\in {J}}u_i\mu_i=\frac{1}{\pi}\mu_{i_0}\bigg(\sum_{i\in {J}}\frac{u_i}{\pi^{\sum_{j=i_0}^{i-1}\rho_{\s_j\wedge\s_{i_0}}}}\prod_{j=i_0}^{i-1}(x-w_{\s_j})\bigg),
\]
and since $\ord_{\Gamma_0}\frac{1}{\pi}\mu_{i_0}=m<0$ we must have \begin{equation}\label{equationord}
\ord_{\Gamma_0}\bigg(\sum_{i\in {J}}\frac{u_i}{\pi^{\sum_{j=i_0}^{i-1}\rho_{\s_j\wedge\s_{i_0}}}}\prod_{j=i_0}^{i-1}(x-w_{\s_j})\bigg)>0.
\end{equation}
For any $j<i\in J$, with $i_0\leq j$ we have $\s_j\not\subset\s_{i_0}$. Therefore either $\s_j=\s_{i_0}$ or $\s_j\wedge\s_{i_0}\supsetneq\s_{i_0}$. 
In the latter case, \[\ord_{\Gamma_0}(x-w_{\s_{i_0}})=m_{\s_{i_0}}\rho_{\s_{i_0}}>m_{\s_{i_0}}\rho_{\s_j\wedge\s_{i_0}}=\ord_{\Gamma_0}(x-w_{\s_j}).\]
It follows from (\ref{equationord}) that
\[\ord_{\Gamma_0}\bigg(\sum_{i\in {J}}v_i\frac{(x-w_{\s_{i_0}})^{\beta_i}}{{\pi^{\beta_i\rho_{\s_{i_0}}}}}\bigg)>0,\]
where $J_i=\{j\in I\mid i_0\leq j<i \mbox{ and } \s_j\neq\s_{i_0}\}$, $v_i=u_i\prod_{j\in J_i}\frac{w_{\s_{i_0}}-w_{\s_{j}}}{\pi^{\rho_{\s_j\wedge\s_{i_0}}}}\in O_K^\times$, and $\beta_i=|\{i_0,\dots,i-1\}\setminus J_i|$.

To find a contradiction, we will use the explicit description of a dense open affine subset of $\Gamma_0$.
Let $W=\{w_1,\dots,w_m\}$ be the set of rational centres of the rationally minimal clusters for $C$ fixed at the beginning of the section. Let $w_h\in W$ such that $w_h=w_{\s_{i_0}}$, and let $L=L_{\s_{i_0}}^{w_h}$, $M=M_{L,0}$, and consider
\[U_M^h\cap\{Z=0\}=\Spec\frac{R[T_M^h(X,Y,Z)^{-1}]}{\lb\mathcal{F}_M^h(X,Y,Z),Z\rb}\subset\Gamma_\t,\]
dense open subscheme of $\Gamma_\t$. From Lemma \ref{MatricesLemma},
\[\sum_{i\in {J}}v_i\frac{(x-w_h)^{\beta_i}}{{\pi^{\beta_i\rho_{\s_{i_0}}}}}
=\sum_{i\in {J}}v_i X^{\beta_i/b_{\s_{i_0}}},
\]
which is a unit since the polynomial $\mathcal{F}_M^h(X,Y,Z)$ in $\{Z=0\}$ is of the form $Y^2-G(X)$ or $Y-G(X)$ for some non-constant $G(X)\in K[X]$ (for more details see Lemma \ref{DefiningEquationXFLemma}). This gives a contradiction and concludes the proof.
\endproof
\end{thm}

Assume now $C_{K^{nr}}$ has an almost rational cluster picture and is $y$-regular as in Theorem \ref{MinimalRegularSNCModelTheorem}. Since $|\Sigma_C|$ is finite, there exists a finite unramified extension $F/K$ such that $C_F$ has an almost rational cluster picture and is $y$-regular. Denote by $O_F$ the ring of integers of $F$.
Let $\Sigma_F=\Sigma_{C_F}^\mathrm{rat}$. Fix a rational centre $w_\s\in F$ for every rationally minimal cluster $\s\in\Sigma_F$. For all non-minimal proper clusters $\t\in\Sigma_F$ choose a rational centre $w_\t=w_\s$ for some rationally minimal cluster $\s\subseteq\t$.
In this setting the next theorem gives a basis of integral differentials of $C$.

\begin{thm}\label{DifferentialsGeneraclCaseTheorem}
Let $C/K$ be a hyperelliptic curve of genus $g\geq 1$ defined by the Weierstrass equation $y^2=f(x)$ and let $\mathcal{C}/O_K$ be a regular model of $C$. Suppose there exists a finite unramified extension $F/K$ such that $C_F$ has an almost rational cluster picture and is $y$-regular. For $i=0,\dots,g-1$ choose inductively proper clusters $\s_i\in\Sigma_F$ so that
\[e_i:=\frac{\epsilon_{\s_i}}{2}-\sum_{j=0}^i\rho_{\s_j\wedge\s_i}=\max_{\t\in\Sigma_F}\bigg\{\frac{\epsilon_{\t}}{2}-\rho_\t-\sum_{j=0}^{i-1}\rho_{\s_j\wedge\t}\bigg\},\]
where if $\s$ and $\s'$ are two possible choices for $\s_i$  satisfying $\s'\subset\s$, then choose $\s_i=\s$. Let $\beta\in O_F^\times$ such that $\mathrm{Tr}_{F/K}(\beta)\in O_K^\times$.
Then the differentials
\[\mu_i=\pi^{\lfloor e_i\rfloor}\cdot\mathrm{Tr}_{F/K}\bigg(\beta\prod_{j=0}^{i-1}(x-w_{\s_j})\bigg)\sfrac{dx}{2y},\qquad i=0,\dots,g-1,\]
form an $O_K$-basis of $H^0(\mathcal{C},\omega_{\mathcal{C}/O_K})$.
\proof
First note that without loss of generality we can suppose $F/K$ Galois. Moreover, since $F/K$ is unramified, $\Gal(F/K)\simeq\Gal(\mathfrak{f}/k)$, where $\mathfrak{f}$ is the residue field of $F$, and so the existence of $\beta$ is guaranteed by the surjectivity of $\mathrm{Tr}_{\mathfrak{f}/k}$. Let $\mathcal{C}$ be the minimal regular model of $C$ over $O_K$. By \cite[Proposition 10.1.17]{Liu}, the base extended scheme $\mathcal{C}_{O_F}$ is the minimal regular model of $C_F$ over $O_F$. 
Let $\mu_0^F,\dots,\mu_{g-1}^F$ be the basis of integral differentials of $C_F$ given by Theorem \ref{DifferentialTheorem}.

Suppose $\mu_0',\dots,\mu_{g-1}'$ is a basis of integral differentials of $C_F$ that, for any $\sigma\in\Gal(F/K)$ and any $j=0,\dots,g-1$, satisfies
\begin{equation}\label{DifferentialsCondition}
\sigma(\mu_j')=\mu_j'+\sum_{0\leq i< j}\lambda_{ij}\mu_i',
\end{equation}
for some $\lambda_{ij}\in O_F$ (depending on $\sigma$). Note that $\mu_0^F,\dots,\mu_{g-1}^F$ is in fact such a basis. We want to prove that, for any $j=0,\dots,g-1$, the differentials 
\begin{equation}\label{DifferentialsTraceEquation}
    \mu_0',\dots,\mu_{j-1}',\mathrm{Tr}_{F/K}(\beta \mu_j'),\mu_{j+1}',\dots, \mu_{g-1}'
\end{equation}
still form a basis of $H^0(\mathcal{C}_F,\omega_{\mathcal{C}_F/O_F})$ satisfying condition (\ref{DifferentialsCondition}). 
From equation (\ref{DifferentialsCondition}) it follows that
\[\mathrm{Tr}_{F/K}(\beta\mu_j')=\sum_{\sigma\in\Gal(F/K)}\sigma(\beta)\sigma(\mu_j')=\mathrm{Tr}_{F/K}(\beta)\mu_j'+\sum_{i<j}\lambda_{ij}'\mu_i',\]
for some $\lambda_{ij}'\in O_F$.
Since $\mathrm{Tr}_{F/K}(\beta)\in O_K^\times$, the differentials in (\ref{DifferentialsTraceEquation}) form a basis of $H^0(\mathcal{C}_F,\omega_{\mathcal{C}_F/O_F})$ satisfying condition (\ref{DifferentialsCondition}).

Since $\mu_0^F,\dots,\mu_{g-1}^F$ satisfies (\ref{DifferentialsCondition}), by induction it follows that 
\[\mathrm{Tr}_{F/K}(\beta\mu_0^F),\dots,\mathrm{Tr}_{F/K}(\beta\mu_{g-1}^F)\] 
is a basis of $H^0(\mathcal{C}_F,\omega_{\mathcal{C}_F/O_F})$. Proposition \ref{DualisingSheafProposition} concludes the proof.
\endproof
\end{thm}

We conclude this section with an application of Theorem \ref{DifferentialTheorem}.

\begin{exa}\label{IntegralDifferentialsExample}
Let $p$ be a prime number and let $a\in\Z_p$, $b\in\Z_p^\times$ such that the polynomial $x^2+ax+b$ is not a square modulo $p$. Let $C$ be the hyperelliptic curve over $\Q_p$ of genus $4$ described by the equation $y^2=f(x)$, where $f(x)=(x^6+ap^4x^3+bp^8)((x-p)^3-p^{11})$. We have already shown in Examples \ref{AlmostRationalExample} and \ref{MinimalRegularModelExample} that $C$ satisfies the hypothesis of Theorem \ref{DifferentialTheorem} and has rational cluster picture
\begin{center}
    \includegraphics[trim=6cm 22.2cm 6cm 4.5cm,clip]{Clusters/Cluster5.pdf}
    \end{center}
We choose rational centres for the minimal clusters $\t_3$ and $\t_4$: $w_{\t_3}=0$ and $w_{\t_4}=p$. Since $\roots=\t_3\wedge\t_4$, we can set either $w_\roots=w_{\t_3}$ or $w_\roots=w_{\t_4}$. Let us fix $w_\roots=w_{\t_3}=0$. Then to choose $\s_0,\s_1,\s_2,\s_3$ as in Theorem \ref{DifferentialTheorem} we draw the following table:

\cellspacetoplimit5pt
\cellspacebottomlimit5pt
\begin{center}
    \small\begin{tabular}{|Sc|Sc|Sc!{\vline width 1pt}Sc|Sc|Sc|Sc|}
\hline
     &  $\rho_\t$&$\epsilon_\t $& $\dfrac{\epsilon_{\t}}{2}-\rho_\t$&$\dfrac{\epsilon_{\t}}{2}-\rho_\t-\rho_{\s_0\wedge\t}$&$\dfrac{\epsilon_{\t}}{2}-\rho_\t-\displaystyle\sum_{j=0}^1\rho_{\s_j\wedge\t}$&$\dfrac{\epsilon_{\t}}{2}-\rho_\t-\displaystyle\sum_{j=0}^2\rho_{\s_j\wedge\t}$\\
     \Xhline{1pt}
     $\t_3$& $\dfrac{4}{3}$ & $11$ & $\dfrac{25}{6}$ & $\color{red}\dfrac{19}{6}$ & $\color{red}\dfrac{11}{6}$ & $\dfrac{1}{2}$\\
     \hline
     $\t_4$&  $\dfrac{11}{3}$ & $17$ & $\color{red}\dfrac{29}{6}$ & $\dfrac{7}{6}$ & $\dfrac{1}{6}$ & $-\dfrac{5}{6}$\\
     \hline
    $\roots$& $1$ & $9$ & $\dfrac{7}{2}$ & $\dfrac{5}{2}$ & $\dfrac{3}{2}$ & $\color{red}\dfrac{1}{2}$\\ \hline
\end{tabular}
\end{center}
$\,$\\
\noindent The numbers in red indicate that $\s_0=\t_4$, $\s_1=\s_2=\t_3$ and $\s_3=\roots$. Thus the differentials
\[\mu_0=p^4\cdot\sfrac{dx}{2y},\quad\mu_1=p^3\cdot(x-p)\sfrac{dx}{2y},\quad\mu_2=p\cdot(x-p)x\sfrac{dx}{2y},\quad\mu_3=(x-p)x^2\sfrac{dx}{2y}\]
form a $\Z_p$-basis of $H^0(\mathcal{C},\omega_{\mathcal{C}/\Z_p})$, for any regular model $\mathcal{C}/\Z_p$ of $C$.
\end{exa}

\appendix
\section{Rational centres over tame extensions}

Let $C/K$ be a hyperelliptic curve given by $y^2=f(x)$.

\begin{lem}\label{RationalCentreTameExtensionLemma}
Let $L/K$ be a field extension. Consider the base extended curve $C_L/L$ and its associated cluster picture $\Sigma_{C_L}$. Let $\s\in\Sigma_{C_L}$ be a proper cluster $G_\s=\mathrm{Stab}_{G_K}(\s)$, and $K_\s=\lb K^\mathrm{s}\rb^{G_\s}$. If $L/L\cap K_\s$ is tamely ramified, then $\s$ has a rational centre $w_\s\in L\cap K_\s$.
\proof
This proof takes ideas from \cite[Lemma B.1]{D2M2}. 
Let $w_\s\in L$ be a rational centre of $\s$ and let $\rho_\s=\max_{w\in L}\min_{r\in\s}v(r-w)$ be its radius. Recall the rationalisation $\s^\mathrm{rat}\in \Sigma_{C_L}^\mathrm{rat}$ of $\s$ (Definition \ref{RationalisationDefinition}).
Denote $\t=\s^\mathrm{rat}$ and define $G_{\t}=\mathrm{Stab}_{G_K}(\t)$. Since $\s\subseteq\t$ we have $G_\s\subseteq G_{\t}$. Furthermore, $\Gal(K^\mathrm{s}/L)\subseteq G_\t$.
Let $F_\s=L\cap K_\s$. Then $\Gal(K^\mathrm{s}/F_\s)\subseteq G_\t$. Since $L/F_\s$ is tamely ramified, we can consider a maximal tamely ramified extension $F_\s^\mathrm{t}$ of $F_\s$ extending $L$. Write $F_\s^{nr}$ for the maximal unramified extension of $F_\s$ in $F_\s^t$. Fix a uniformiser $\pi_\s$ of $F_\s$. Since $L/F_\s$ is tamely ramified and $w_\s\in L$, for a sufficiently large $b$ fix a choice of  $\sqrt[b]{\pi_\s}$ such that $w_\s\in F_\s^{nr}(\sqrt[b]{\pi_\s})$. 
Write the $v$-adic expansion of $w_\s$ as
\[w_\s=u_t\sqrt[b]{\pi_\s}^t+u_{t+1}\sqrt[b]{\pi_\s}^{t+1}+\dots\]
for a suitable $t\in\Z$, with $u_l\in F_\s^{nr}$ canonical representatives of elements in $k^\mathrm{s}$.
Define
\[w=\sum_{l<e_{F_\s/K}b\rho_\s}u_l\sqrt[b]{\pi_\s}^l.\]
We first show that $w\in F_\s^\mathrm{t}$. It trivially follows if $w=0$. Suppose $0\neq w\notin F_\s^\mathrm{t}$, and that $u_{l_0}\sqrt[b]{\pi_\s}^{l_0}$ is the lowest valuation term of the expansion which is not in $F_\s^\mathrm{t}$. Let $w'=\sum_{l<l_0}u_l\sqrt[b]{\pi_\s}^l$. Note that $w'\in F_\s^\mathrm{t}$ for our assumption on $l_0$. As $v(w-w_\s)\geq\rho_\s$, we have $v(w_\s-w')=v(w-w')=l_0/e_{F_\s/K}b$. Since $L\subseteq F_\s^\mathrm{t}$, we have $w_\s-w'\in F_\s^\mathrm{t}$ and so the denominator of $l_0/b$ is not divisible by $p$. But then $u_{l_0}\sqrt[b]{\pi_\s}^{l_0}\in F_\s^\mathrm{t}$ as $u_{l_0}\in F_\s^{nr}\subseteq F_\s^\mathrm{t}$ and $\sqrt[b]{\pi_\s}^{l_0}\in F_\s^\mathrm{t}$.

Let $\mathcal{D}_\t=\{x\in K^\mathrm{s}\mid v(x-w_\s)\geq\rho_\s\}$ be the smallest disc in $K^\mathrm{s}$ cutting out $\t$. Note that $\mathrm{Stab}_{G_K}(\mathcal{D}_\t)=G_\t$. Since $w\in\mathcal{D}_\t$, for $\sigma\in \Gal(K^\mathrm{s}/F_\s)\subseteq G_\t$ we have $\sigma(w)\in\mathcal{D}_\t$ and so $v(\sigma(w)-w_\s)\geq\rho_\s$. Therefore the terms in the $v$-adic expansions of $\sigma(w)$ and $w$ agree up to $\sqrt[b]{\pi_\s}^{e_{F_\s/K}b\rho_\s}$ (excluded). Furthermore, if $w\in L$, then $w$ is a rational centre of $\s$. Indeed, for any $r\in\s$ one has \[v(r-w)\geq \min\{v(r-w_\s),v(w-w_\s)\}\geq\rho_\s.\]

We showed $w\in F_\s^\mathrm{t}$. It remains to prove that $w\in F_\s$, i.e.\ it is $\Gal(K^\mathrm{s}/F_\s)$-invariant. Suppose not, and that $u_l\sqrt[b]{\pi_\s}^l$ is the lowest valuation term of the expansion which is not $\Gal(K^\mathrm{s}/F_\s)$-invariant. Note that the denominator of $l/b$ is not divisible by $p$ since $w\in F_\s^\mathrm{t}$. If $b\nmid l$, then there is some element $\sigma$ of tame inertia of $F_\s$ which fixes $u_l\in F_\s^{nr}$ and maps $\sqrt[b]{\pi_\s}^l$ to $\zeta\sqrt[b]{\pi_\s}^l$, where $\zeta\neq 1$ is a root of unity; this contradicts the fact that $\sigma(w)\equiv w\mod\sqrt[b]{\pi_\s}^{e_{F_\s/K}b\rho_\s}$. If $b\mid l$, then 
we must have $u_l\notin F_\s$. Then there exists some element $\sigma\in\Gal(F_\s^{nr}/F_\s)$ so that $\sigma(u_l)\neq u_l$; this contradicts $\sigma(w)\equiv w\mod \sqrt[b]{\pi_\s}^{e_{F_\s/K}b\rho_\s}$ similarly to before.
\endproof
\end{lem}

\section{Dualising sheaf under base extensions}
Let $F/K$ be a finite Galois extension and let $O_F$ be the ring of integers of $F$.

\begin{lem}
Let $M$ be a flat $O_K$-module and $M_F:=M\otimes_{O_K}O_F$. Then \[M\simeq M_F^{\Gal(F/K)}=\{m\in M_F\mid \sigma(m)=m\mbox{ for every }\sigma\in\Gal(F/K)\}.\]
\proof
As $M$ is flat, the functor $M\otimes_{O_K}-$ is (left) exact. From the isomorphism $O_K\simeq O_F^{\Gal(F/K)}$ it follows that \[M\otimes_{O_K}O_K\simeq M\otimes_{O_K}O_F^{\Gal(F/K)},\] that is $M\simeq M_F^{\Gal(F/K)}$, as required.
\endproof
\end{lem}

\begin{prop}\label{DualisingSheafProposition}
Let $C$ be a smooth projective curve of genus $g\geq 1$ and let $\mathcal{C}$ be a regular model of $C$ over $O_K$. Denote by $C_F$ and $\mathcal{C}_{O_F}$ the base extended schemes. Then
$H^0(\mathcal{C}_F,\omega_{\mathcal{C}_F/O_F})\simeq H^0(\mathcal{C},\omega_{\mathcal{C}/O_K})\otimes_{O_K}O_F$ and\[ H^0(\mathcal{C},\omega_{\mathcal{C}/O_K})\simeq H^0(\mathcal{C}_F,\omega_{\mathcal{C}_F/O_F})^{\Gal(F/K)}.\]
\proof
The 
lemma follows from the following results: \cite[Proposition 10.1.17]{Liu}, \cite[Theorem 6.4.9(b)]{Liu}, \cite[Exercise 6.4.6]{Liu}, \cite[Corollary 5.2.27]{Liu} and the previous lemma. 
\endproof
\end{prop}



\end{document}